\RequirePackage{ifpdf}
\ifpdf 
\documentclass[pdftex]{sigma}
\else
\documentclass{sigma}
\fi

\usepackage{rotating}
\usepackage{xypic}
\usepackage{pst-node,pstricks,multido,pst-plot,pst-text,pst-3d}%
\usepackage{tikz}
\usepackage{tkz-graph}
\usepackage{tkz-berge}
\usetikzlibrary{arrows}
\usetikzlibrary{calc,3d}
\usepackage{soul}

\definecolor{orange}{rgb}{1,0.5,0}

\def\S{{\mathfrak S}}
\def\cal#1{{\mathfrak #1}}

\def\<{\langle}
\def\>{\rangle}
\def\R{{\mathbb R}}
\def\C{{\mathbb C}}
\def\Z{{\mathbb Z}}
\def\N{{\mathbb N}}

\def\CT{{\rm C.T.}}

\def\goth{\mathfrak}

\def\ashuff#1#2#3{
\kern 1pt \vrule height#1 \overline{\vrule height#3 width 0pt
\hskip#2} \rule{.3pt}{#1}\overline{\vrule height#3 width 0pt
\hskip#2} \rule{.3pt}{#1} \kern 1pt }

\def\CT{{\rm CT}}

\numberwithin{equation}{section}

\newtheorem{Theorem}{Theorem}[section]
\newtheorem{Lemma}[Theorem]{Lemma}
\newtheorem{Proposition}[Theorem]{Proposition}
\newtheorem{Corollary}[Theorem]{Corollary}
{\theoremstyle{definition}
\newtheorem{Definition}[Theorem]{Definition}
\newtheorem{Example}[Theorem]{Example}
}

\begin{document}

\allowdisplaybreaks

\renewcommand{\PaperNumber}{026}

\FirstPageHeading

\ShortArticleName{Vector-Valued Jack Polynomials from Scratch}

\ArticleName{Vector-Valued Jack Polynomials from Scratch}
\Author{Charles F. DUNKL~$^\dag$ and Jean-Gabriel LUQUE ~$^\ddag$}

\AuthorNameForHeading{C.F. Dunkl and J.-G. Luque}

\Address{$^\dag$~Dept. of Mathematics, University of Virginia, Charlottesville VA 22904-4137, USA}
\EmailD{\href{mailto:cfd5z@virginia.edu}{cfd5z@virginia.edu}}
\URLaddressD{\url{http://people.virginia.edu/~cfd5z/}}
\Address{$^\ddag$~Universit\'e de Rouen, LITIS
Saint-Etienne du Rouvray, France}
\EmailD{\href{mailto:jean-gabriel.luque@univ-rouen.fr}{jean-gabriel.luque@univ-rouen.fr}}
\URLaddressD{\url{http://www-igm.univ-mlv.fr/~luque/}}

\ArticleDates{Received September 21, 2010, in f\/inal form March 11, 2011;  Published online March 16, 2011}

\Abstract{Vector-valued Jack polynomials associated to the symmetric group $\S_N$
are polynomials with
multiplicities in an irreducible module of $\S_N$ and which are
simultaneous eigenfunctions of the Cherednik--Dunkl operators with some
additional properties concerning the leading monomial.
These polynomials were introduced by Grif\/feth in the general setting of
the complex ref\/lections groups $G(r,p,N)$ and studied by one of the
authors (C.~Dunkl) in the specialization $r=p=1$ (i.e.\ for the
symmetric group).
By adapting a construction due to Lascoux, we describe an algorithm
allowing us to compute explicitly the Jack polynomials following a
Yang--Baxter graph. We recover some properties already studied by C.~Dunkl
and restate them in terms of graphs together with additional new
results. In particular, we investigate normalization, symmetrization and
antisymmetrization, polynomials with minimal degree, restriction
etc. We give also a shifted version of the construction and we
discuss  vanishing properties of the associated
polynomials.}


\Classification{05E05; 16T25; 05C25; 33C80}

\section{Introduction}

 The Yang--Baxter graphs are used
to study Jack polynomials. In particular, these objects have been
investigated in this context by Lascoux \cite{Las} (see also
\cite{Macdo,MacdoHecke} for general properties about Jack and Macdonald
polynomials).  Vector-valued Jack polynomials
 are associated with irreducible
representations of the symmetric group $\mathfrak{S}_{N}$, that is, to
partitions of $N$. A Yang--Baxter graph is a directed graph with no loops
and a
unique root, whose edges are labeled by generators of a certain subsemigroup
of the extended af\/f\/ine symmetric group. In this paper the vertices are
labeled
by a pair consisting of a  weight in $\N^N$ and the content vector of a standard
tableau. The  weights  are the labels of
monomials which are the leading
terms of polynomials, and the tableaux all have the same shape. There is a
vector-valued Jack polynomial associated with each vertex. These polynomials
are special cases of the polynomials introduced by Grif\/feth~\cite{Gr,Gr2}
for the
family of complex ref\/lection groups denoted by $G\left( r,p,N\right) $
(where~$p|r$). This is the group of unitary $N\times N$ matrices such that
their nonzero entries are $r^{\rm th}$ roots of unity, the product of the
nonzero
entries is a $\left( r/p\right) ^{\rm th}$ root of unity, and there is exactly
one nonzero entry in each row and each column. The symmetric group is the
special case $G\left( 1,1,N\right) $. The vector space in which the Jack
polynomials take their values is equipped with the nonnormalized basis
described by Young, namely, the simultaneous eigenvectors of the
Jucys--Murphy
elements.

The labels on edges denote transformations to be applied to the objects at a
vertex. Vector-valued Jack polynomials are uniquely determined by their
spectral vector, the vector of eigenvalues under the (pairwise commuting)
Cherednik--Dunkl operators. This serves to demonstrate the claim that
dif\/ferent
paths from one vertex to another produce the same result, a situation
which is
linked to the braid or Yang--Baxter relations. These refer to the
transformations.

Following Lascoux \cite{Las} we def\/ine the monoid
$\widehat{\mathfrak{S}}_{N}$, a subsemigroup of the af\/f\/ine symmetric group, with generators $\left\{
s_{1},s_{2},\ldots,s_{N-1},\Psi\right\} $ and relations:%
\begin{gather*}
s_{i}s_{j}  =s_{j}s_{i},\qquad \left\vert i-j\right\vert >1,\\
s_{i}s_{i+1}s_{i}  =s_{i+1}s_{i}s_{i+1},\qquad 1\leq i<N-1,\\
s_{1}\Psi^{2}  =\Psi^{2}s_{N-1},\\
s_{i}\Psi  =\Psi s_{i-1},\qquad 2\leq i\leq N-1.
\end{gather*}
The relations $s_{i}^{2}=1$ do not appear in this list because the graph has no loops.

The main objects of our study are polynomials in $x=\left( x_{1},\ldots,x_{N}\right) \in\mathbb{R}^{N}$ with coef\/f\/icients in $\mathbb{Q}\left(
\alpha\right) $, where $\alpha$ is a transcendental (indeterminate), and
with
values in the $\mathfrak{S}_{N}$-module corresponding to a partition
$\lambda$
of $N$. The Yang--Baxter graph $G_{\lambda}$ is a pictorial representation of
the algorithms which produce the Jack polynomials starting with
constants. The
generators of $\widehat{\mathfrak{S}}_{N}$ correspond to transformations
taking a Jack polynomial to an adjacent one. At each vertex there is such a
polynomial, and a 4-tuple which identif\/ies it. The 4-tuple consists of a~standard tableau denoting a basis element of the
$\mathfrak{S}_{N}$-module, a
 weight  (multi-index) descri\-bing the
leading monomial of the polynomial, a
spectral vector, and a permutation, essentially the rank function of the
 weight. The spectral vector and
permutation are determined by the f\/irst
two elements. For technical reasons the standard tableaux are actually
reversed, that is, the entries decrease in each row and each column. This
convention avoids the use of a reversing permutation, in contrast to
Grif\/feth's paper~\cite{Gr2} where the standard tableaux have the usual
ordering.

The symmetric and antisymmetric Jack polynomials are constructed in terms of
certain subgraphs of~$G_{\lambda}$. Furthermore the graph technique leads to
the def\/inition and construction of shifted inhomogeneous vector-valued
Jack polynomials.

Here is an outline of the contents of each section.

Section~\ref{s2YB} contains the basic def\/initions and construction of
the graph~$G_{\lambda}$. The presentation is in terms of the 4-tuples mentioned above.
It is important to note that not every possible label need appear on edges
pointing away from a given vertex: if the e
weight  at the vertex is
$v\in\mathbb{N}^{N}$ then the transposition $\left( i,i+1\right) $ (labeled
by $s_{i}$) can be applied only when $v\left[ i\right] \leq v\left[
i+1\right] $, that is, when the resulting
weight is greater than or
equal to $v$ in the dominance order. The action of the af\/f\/ine element $\Psi$
is given by $v\Psi=\left( v\left[ 2\right] ,v\left[ 3\right]
,\ldots,v\left[ N\right] ,v\left[ 1\right] +1\right) $.

The Murphy basis for the irreducible representation of $\mathfrak{S}_{N}$
along with the def\/inition of the action of the simple ref\/lections $\left(
i,i+1\right) $ on the basis is presented in Section~\ref{s3VectorValued}
Also the
vector-valued polynomials, their partial ordering, and the Cherednik--Dunkl
operators are introduced here.

Section~\ref{s4NonSym} is the detailed development of Jack polynomials.
Each edge of the
graph $G_{\lambda}$ determines a transformation that takes the Jack
polynomial
associated with the beginning vertex to the one at the ending vertex of the
edge. There is a canonical pairing def\/ined for the vector-valued
polynomials;
the pairing is nonsingular for generic~$\alpha$ and the Cherednik--Dunkl
operators are self-adjoint. The Jack polynomials are pairwise orthogonal for
this pairing and the squared norm of each polynomial can be found by use of
the graph.

In Section~\ref{s5Sym} we investigate the symmetric and
antisymmetric vector-valued Jack polynomials in relation with
connectivity of the Yang--Baxter graph whose af\/f\/ine edges have been removed. Also, in
this section one f\/inds the method of producing coef\/f\/icients so that the
corresponding sum of Jack polynomials is symmetric or antisymmetric.
The idea is explained in terms of certain subgraphs of~$G_{\lambda}$.

 Vertices of
$G_{\lambda}$ satisfying certain conditions may be mapped to vertices of a
graph related to~$\mathfrak{S}_{M}$, $M<N$, by a restriction map. This
topic
is the subject of Section~\ref{s6Rest}. This section also describes the restriction map on the
Jack polynomials.

In Section~\ref{s7shift} the shifted vector-valued Jack polynomials are
presented.
These are inhomogeneous and the parts of highest degree coincide with the
homogeneous Jack polynomials of the previous section. The construction again
uses the Yang--Baxter graph $G_{\lambda}$; it is only necessary to change the
operations associated with the edges.

Throughout the paper there are numerous f\/igures to concretely illustrate the
structure of the graphs.

\section[Yang-Baxter type graph associated to a partition]{Yang--Baxter type graph associated to a partition}\label{s2YB}

\subsection{Sorting a vector}\label{ss2.1vector}

Consider a vector $v\in\N^N$, we want to compute the unique decreasing
partition $v^{+}$, which is in the orbit of $v$ for the action of the
symmetric group $\S_N$ acting on  the  right on the position,
using the minimal number of elementary transpositions $s_i=(i\,i+1)$.

If $v$ is a vector we will denote by $v[i]$ its $i$th component. Each
$\sigma\in\S_N$ will be associated to the vector of its images
$[\sigma(1),\dots,\sigma(N)]$.
Let $\sigma$ be a permutation, we will denote
$\ell(\sigma)=\min\{k:\sigma=s_{i_1}\cdots s_{i_k}\}$ the length of the
permutation. By a straightforward induction one f\/inds:

\begin{Lemma}\label{unicity_sigma_v}
Let $v\in\N^N$ be a vector, there exists a unique permutation $\sigma_v$
such that $v=v^{+}\sigma_v$ with $\ell(\sigma_v)$ minimal.
\end{Lemma}
The permutation $\sigma_v$ is obtained by a standardization process: we
label with integer from $1$ to $N$ the positions in $v$ from the largest
entries to the smallest one and from left to right.
\begin{Example}\rm
Let $v=[2,3,3,1,5,4,6,6,1]$, the construction gives:
\[
\begin{array}{@{}r@{\,\,}c@{\,\,}ccccccccccc}\sigma_v&=&[&7&5&6&8&3&4&1&2&9&] \\
v&=&[&2&3&3&1&5&4&6&6&1&]
\end{array}\]
We verify that $v\sigma_v^{-1}=[6,6,5,4,3,3,2,1,1]=v^{+}$.
\end{Example}

The def\/inition  of $\sigma_v$ is
compatible with the action of $\S_N$ in the following sense:

\begin{Proposition}\label{si2sigma_v} \qquad
\begin{enumerate}\itemsep=0pt
\item
  $ \sigma_{vs_i}=\left\{\begin{array}{ll}\sigma_v&\mbox{ if }v=vs_i,\\
\sigma_vs_i&\mbox{ otherwise.} \end{array}
\right.$
\item If $v[i]=v[i+1]$ then $\sigma_vs_i=s_{\sigma_v[i]}\sigma_v$.
\end{enumerate}
\end{Proposition}

This can be easily obtained from the construction.

Def\/ine the af\/f\/ine operation $\Psi$ acting on a vector by
\[
[v_1,\dots,v_N]\Psi=[v_2,\dots,v_N,v_1+1],
\]
{ and more generally let $\Psi^\alpha$ by
\[
[v_1,\dots,v_N]\Psi^\alpha=[v_2,\dots,v_N,v_1+\alpha].
\]
Denote also by $\theta:=\Psi^0$ the circular permutation
$[2,\dots,N,1]$.
}

Again, one can prove easily that the computation of $\sigma_v$ is
compatible (in a certain sense) with the action of $\Psi$:

\begin{Proposition}\label{Psi2sigma_v}
\[
\sigma_{v\Psi}=\sigma_v \theta.
\]
\end{Proposition}

\begin{Example}
Consider $v=[2,3,3,2,5,4,6,6,1]$, one has
\[\sigma_v=[7,5,6,8,3,4,1,2,9]\] and \[\sigma_v{
\theta}=[5,6,8,3,4,1,2,9,7].\]
But \[v':=v[2,3,4,5,6,7,8,9,1]=[3,3,2,5,4,6,6,1,2]\] and
$\sigma_{v'}=[5,6,\underline{7},3,4, 1,2,9,\underline{8}]$; here
 an
underlined integer
 means that there is a dif\/ference with the same position in
$\sigma_v{\theta}$ and in {$\sigma_{v'}$}. This is due to
the fact that $v[1]$ is the f\/irst occurrence of $2$ in $v$ while
$v'[9]$ is the last occurrence of $2$ in $v'$. Adding $1$ to $v'[9]$ one
obtains
\[
 v\Psi=[3,3,2,5,4,6,6,1,3].
\]
The last occurrence of $2$ becomes the last occurrence of $3$ (that is
the number of the f\/irst occurrence of $2$ minus $1$). Hence,
\[
 \sigma_{v\Psi}=[5,6,8,3,4,1,2,9,7]=\sigma_v{\theta}.
\]
\end{Example}

\subsection{Construction and basic properties of the
graph}\label{ss2.2basic}

\begin{Definition}
A {\it tableau} of shape $\lambda$ is a
f\/illing with integers weakly increasing in each row and in each column.
In the sequel
{\it row-strict} means increasing in each row and {\it column-strict}
means increasing in
each column.

A {\it reverse standard tableau} (RST) is obtained by f\/illing the shape
$\lambda$ with integers $1,\dots, N$ and with the conditions of strictly
decreasing in the  row and the column. We
will denote by ${\rm Tab}_\lambda$, the set of the RST with shape
$\lambda$.

Let $\tau$ be a RST, we def\/ine the {\it vector of contents} of $\tau$ as
the vector $\CT_\tau$ such that
$\CT_\tau[i]$ is the {\it content} of $i$ in $\tau$ (that is the number
of the diagonal in which $i$ appears; the number of the main diagonal is~0, and the numbers decrease from down to up). {In other words, if
$i$ appears in the box $[{\rm col},{\rm row}]$ then $\CT_\tau[i]={\rm col}-{\rm row}$.}
\end{Definition}

\begin{Example}
 Consider the tableau $\tau=\begin{array}{ccc}
2\\
5&4\\
6&3&1\\
\end{array}$, we obtain the vector of contents by labeling the
numbers of the diagonals $\begin{array}{rcc}
-2\\
-1&0\\
0&1&2\\
\end{array}$. So one obtains,
\[
\CT_{\mbox{\tiny $\begin{array}{ccc}
2\\
5&4\\
6&3&1\\
\end{array}$}}=[2,-2,1,0,-1,0].
\]
\end{Example}

We construct a Yang--Baxter-type graph with vertices labeled by 4-tuples
$(\tau, \zeta, v,\sigma)$, where $\tau$ is a RST, $\zeta$ is a vector
of length $N$ with entries in $\Z[\alpha]$ ($\zeta$ will
be called the {\it spectral vector}), $v\in\N^N$ and $\sigma\in\S_N$, as
follows: First, consider a RST of shape $\lambda$ and write a vertex
labeled by the 4-tuple $(\tau,\CT_\tau, 0^N, [1,\dots,N])$. Now, we
consider the action of the elementary transposition of $\S_N$ on the
4-tuple given by
\[(\tau,\zeta,v,\sigma)s_i=\left\{\begin{array}{ll} (\tau,\zeta
s_i,vs_i,\sigma s_i)&\mbox{ if }v[i+1]\neq v[i],\vspace{1mm} \\
\big(\tau^{(\sigma[i],\sigma[i+1])},\zeta s_i,v,\sigma\big) &\mbox{ if }
v[i]=v[i+1]\mbox{ and }\tau^{(\sigma[i],\sigma[i+1])}\in{\rm
Tab}_\lambda,\vspace{1mm}\\
(\tau,\zeta,v,\sigma) &\mbox{ otherwise,}\end{array}\right.
\]
where $\tau^{(i,j)}$ denotes the f\/illing obtained by permuting the
values $i$ and $j$ in $\tau$.

 Consider also the af\/f\/ine action given by
\[ (\tau,\zeta,v,\sigma)\Psi=
(\tau,\zeta\Psi^\alpha,v\Psi,\sigma[2,\dots,N,1])
 =(\tau,[\zeta_2,\dots,],[v_2,\dots],[\sigma_2,\dots]).\]
\begin{Example}\quad\small
\begin{enumerate}\itemsep=0pt
\item \ \\ $\begin{array}{r}
\left({31\ \atop
542},[1,0,2\alpha,\alpha+2,\alpha-1],[0,0,2,1,1],[45123]\right)s_2=\ \ \
\ \ \ \ \ \ \ \ \ \ \ \ \ \ \ \ \ \ \ \ \ \\\left({31\ \atop
542},[1,2\alpha,0,\alpha+2,\alpha-1],[0,2,0,1,1],[41523]\right)\end{array}
 $
\item \ \\$\begin{array}{r}
\left({31\ \atop
542},[1,0,2\alpha,\alpha+2,\alpha-1],[0,0,2,1,1],[45123]\right)s_4=\ \ \
\ \ \ \ \ \ \ \ \ \ \ \ \ \ \ \ \ \ \ \ \ \\\left({21\ \atop
543},[1,0,2\alpha,\alpha-1,\alpha+
2],[0,0,2,1,1],[45123]\right)\end{array}
 $
\item \ \\$\begin{array}{r}
\left({31\ \atop
542},[1,0,2\alpha,\alpha+2,\alpha-1],[0,0,2,1,1],[45123]\right)s_1=\ \ \
\ \ \ \ \ \ \ \ \ \ \ \ \ \ \ \ \ \ \ \ \ \\\left({31\ \atop
542},[1,0,2\alpha,\alpha+2,\alpha-1],[0,0,2,1,1],[45123]\right)\end{array}
 $
\item \ \\$\begin{array}{r}
\left({31\ \atop
542},[1,0,2\alpha,\alpha+2,\alpha-1],[0,0,2,1,1],[45123]\right)\Psi=\ \
\ \ \ \ \ \ \ \ \ \ \ \ \ \ \ \ \ \ \ \ \ \ \\\left({31\ \atop
542},[0,2\alpha,\alpha+2,\alpha-1,\alpha+1],[0,2,1,1,1],[51234]\right)\end{array}
 $
\end{enumerate}
\end{Example}

\begin{Definition} \label{Glambda}
If $\lambda$ is a partition, denote by $\tau_\lambda$ the tableau
obtained by f\/illing the shape $\lambda$ from bottom to top and left to
right by the integers $\{1,\dots,N\}$ in the decreasing order.

The graph $G_\lambda$ is an inf\/inite directed graph constructed from the
4-tuple \[(\tau_\lambda,\CT_{\tau_\lambda},[0^N],[1,2,\dots,N]),\]
called the {\it root}   and adding
vertices and edges following the rules
\begin{enumerate}\itemsep=0pt
\item We add an arrow labeled by $s_i$ from the vertex
$(\tau,\zeta,v,\sigma)$ to $(\tau',\zeta',v',\sigma')$ if
$(\tau,\zeta,v,\sigma)s_i=(\tau',\zeta',v',\sigma')$ and $v[i]<v[i+1]$
or $v[i]=v[i+1]$ and $\tau$ is obtained from $\tau'$ by interchanging
the position of two integers $k<\ell$ such that $k$ is at the south-east
of $\ell$ (i.e.\ $\CT_\tau(k)\geq \CT_\tau(\ell)+2$).
\item We add an arrow labeled by $\Psi$ from the vertex
$(\tau,\zeta,v,\sigma)$ to $(\tau',\zeta',v',\sigma')$ if
$(\tau,\zeta,v,\sigma)\Psi=(\tau',\zeta',v',\sigma')$.
\item We add an arrow $s_i$ from the vertex $(\tau,\zeta,v,\sigma)$ to
$\varnothing$ if $(\tau,\zeta,v,\sigma)s_i=(\tau,\zeta,v,\sigma)$.
 \end{enumerate}
An arrow of the form
\begin{center}
\begin{tikzpicture}
\GraphInit[vstyle=Shade]
    \tikzstyle{VertexStyle}=[shape = rectangle,
draw
]
\SetUpEdge[lw = 1.5pt,
color = orange,
 labelcolor = gray!30,
 style={post},
labelstyle={sloped}
]
\tikzset{LabelStyle/.style = {draw,
                                     fill = white,
                                     text = black}}
\tikzset{EdgeStyle/.style={post}}
\Vertex[x=0, y=0,
 L={\tiny$(\tau,\zeta,v,\sigma)$}]{x}
\Vertex[x=4, y=0,
 L={\tiny$(\tau,\zeta',v',\sigma')$}]{y}
\Edge[label={\tiny$s_i$ or $\Psi$}](x)(y)
\end{tikzpicture}
\end{center}
will be called a {\it step}. The other arrows will be called {\it
jumps}, and in particular an arrow
\begin{center}
\begin{tikzpicture}
\GraphInit[vstyle=Shade]
    \tikzstyle{VertexStyle}=[shape = rectangle,
draw
]
\SetUpEdge[lw = 1.5pt,
color = orange,
 labelcolor = gray!30,
 style={post},
labelstyle={sloped}
]
\tikzset{LabelStyle/.style = {draw,
                                     fill = white,
                                     text = black}}
\tikzset{EdgeStyle/.style={post}}
\Vertex[x=0, y=0,
 L={\tiny$(\tau,\zeta,v,\sigma)$}]{x}
\Vertex[x=3, y=0,
 L={\tiny$\varnothing$}]{y}
\Edge[label={\tiny$s_i$}](x)(y)
\end{tikzpicture}
\end{center}
 will be called a {\it fall}; the other jumps will be called {\it
correct jumps}.

As usual a {\it path} is a
 sequence of consecutive arrows in $G_\lambda$ starting
 from the root and is denoted by the sequence if the labels of its
arrows.
Two paths ${\cal P}_1=(a_1,\dots,a_k)$ and ${\cal
P}_2=(b_1,\dots,b_\ell)$ are said to be {\it equivalent} (denoted by
${\cal P}_1\equiv{\cal P}_2$) if they lead to the same vertex.
\end{Definition}

We remark that from Proposition~\ref{si2sigma_v}, in the case $v[i]=v[i+1]$, the part~1 of Def\/inition~\ref{Glambda} is equivalent to the following statement: $\tau'$ is
obtained from $\tau$ by interchanging $\sigma_v[i]$ and
$\sigma_v[i+1]=\sigma_v[i]+1$ where $\sigma_v[i]$ is to the south-east
of $\sigma_v[i]+1$, that is, $\CT_v[\sigma_v[i]]-\CT_v[\sigma_v[i]+1]\geq2$.

\begin{Example}
The following arrow is a correct jump
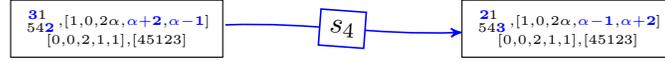
\begin{center}
\begin{tikzpicture}
\GraphInit[vstyle=Shade]
    \tikzstyle{VertexStyle}=[shape = rectangle,
draw
]
\SetUpEdge[lw = 1.5pt,
color = orange,
 labelcolor = gray!30,
 style={post},
labelstyle={sloped}
]
\tikzset{LabelStyle/.style = {draw,
                                     fill = white,
                                     text = black}}
\tikzset{EdgeStyle/.style={post}}
\Vertex[x=0, y=0,
 L={\tiny${{\blue \bf3}1\ \atop 54{\blue \bf2}}, [1,0,2\alpha,{\blue
\bf\alpha+2},{\blue \bf\alpha-1}]\atop [0,0,2,1,1], [45123]
 $}]{x}
\Vertex[x=6, y=0,
 L={\tiny${{\blue \bf2}1\ \atop 54{\blue \bf3}}, [1,0,2\alpha,{\blue
\bf\alpha-1},{\blue \bf\alpha+2}]\atop [0,0,2,1,1], [45123]
 $}]{y}
\Edge[label={$s_4$},style={post,in=182,out=2},color=blue](x)(y)
\end{tikzpicture}
\end{center}
 whilst
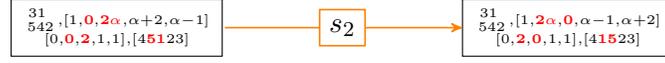
\begin{center}
\begin{tikzpicture}
\GraphInit[vstyle=Shade]
    \tikzstyle{VertexStyle}=[shape = rectangle,
draw
]
\SetUpEdge[lw = 1.5pt,
color = orange,
 labelcolor = gray!30,
 style={post},
labelstyle={sloped}
]
\tikzset{LabelStyle/.style = {draw,
                                     fill = white,
                                     text = black}}
\tikzset{EdgeStyle/.style={post}}
\Vertex[x=0, y=0,
 L={\tiny${31\ \atop 542}, [1,{\red \bf 0},{\red \bf
2\alpha},\alpha+2,\alpha-1]\atop [0,{\red \bf 0},{\red \bf2},1,1],
[4{\red \bf51}23]
 $}]{x}
\Vertex[x=6, y=0,
 L={\tiny${31\ \atop 542}, [1,{\red \bf2\alpha},{\red \bf
0},\alpha-1,\alpha+2]\atop [0,{\red \bf2},{\red \bf0},1,1], [4{\red
\bf15}23]
 $}]{y}
\Edge[label={$s_2$}](x)(y)
\end{tikzpicture}
\end{center}
is a step.

The arrows
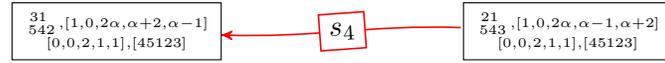
\begin{center}
\begin{tikzpicture}
\GraphInit[vstyle=Shade]
    \tikzstyle{VertexStyle}=[shape = rectangle,
draw
]
\SetUpEdge[lw = 1.5pt,
color = orange,
 labelcolor = gray!30,
 style={post},
labelstyle={sloped}
]
\tikzset{LabelStyle/.style = {draw,
                                     fill = white,
                                     text = black}}
\tikzset{EdgeStyle/.style={post}}
\Vertex[x=0, y=0,
 L={\tiny${{3}1\ \atop 54{2}}, [1,0,2\alpha,{\alpha+2},{\alpha-1}]\atop
[0,0,2,1,1], [45123]
 $}]{x}
\Vertex[x=6, y=0,
 L={\tiny${{2}1\ \atop 54{3}}, [1,0,2\alpha,{\alpha-1},{\alpha+2}]\atop
[0,0,2,1,1], [45123]
 $}]{y}
\Edge[label={$s_4$},style={post,in=-2,out=178},color=red](y)(x)
\end{tikzpicture}
\end{center}
and
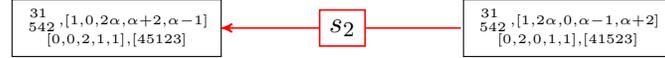
\begin{center}
\begin{tikzpicture}
\GraphInit[vstyle=Shade]
    \tikzstyle{VertexStyle}=[shape = rectangle,
draw
]
\SetUpEdge[lw = 1.5pt,
color = orange,
 labelcolor = gray!30,
 style={post},
labelstyle={sloped}
]
\tikzset{LabelStyle/.style = {draw,
                                     fill = white,
                                     text = black}}
\tikzset{EdgeStyle/.style={post}}
\Vertex[x=0, y=0,
 L={\tiny${31\ \atop 542}, [1,0,2\alpha,\alpha+2,\alpha-1]\atop
[0,0,2,1,1], [45123]
 $}]{x}
\Vertex[x=6, y=0,
 L={\tiny${31\ \atop 542}, [1,2\alpha,0,\alpha-1,\alpha+2]\atop
[0,2,0,1,1], [41523]
 $}]{y}
\Edge[label={$s_2$},color=red](y)(x)
\end{tikzpicture}
\end{center}
are not allowed.

\end{Example}

\begin{Example}\rm
Consider the partition $\lambda=21$, the graph $G_{21}$ in Fig.~\ref{G21} is obtained from the 4-tuple
$\left(\begin{array}{cc}2\\3&1\end{array},(1,-1,0),(0,0,0),1\right)$ by
applying the rules of Def\/inition~\ref{Glambda}. In Fig.~\ref{G21}, the steps
are drawn in orange, the jumps in blue and the falls have been omitted.

%
%

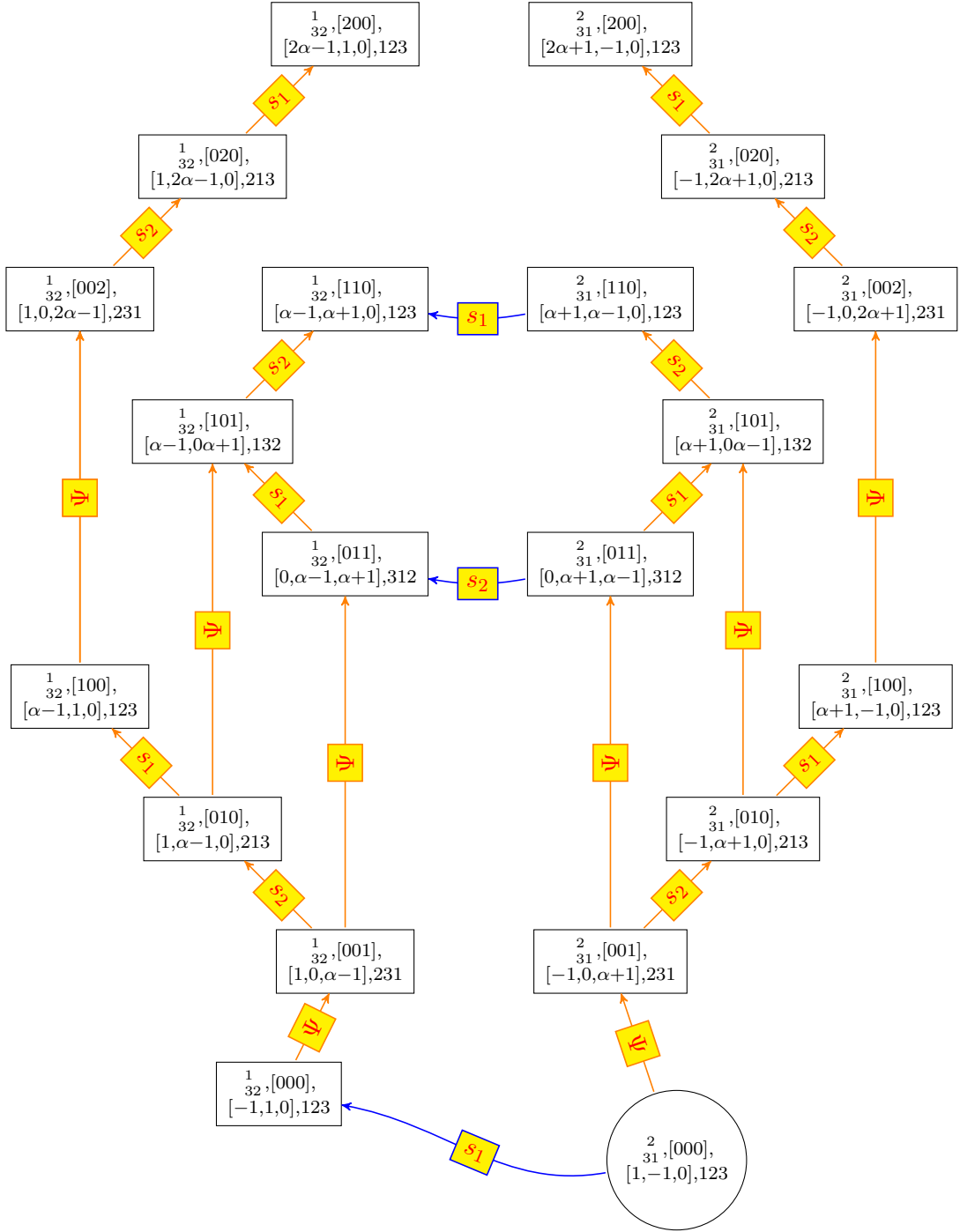
\begin{figure}[h]\centering
\begin{tikzpicture}%
\GraphInit[vstyle=Shade]
    \tikzstyle{VertexStyle}=[shape = rectangle,
                             draw
]

\Vertex[x=3, y=-1,
 L={${2\atop31}, [000],\atop
[1,-1,0],123$},style={shape=circle}
]{x2}
\SetUpEdge[lw = 1.5pt,
color = orange,
 labelcolor = gray!30,
 labelstyle = {draw,sloped},
 style={post}
]

\tikzset{LabelStyle/.style = {draw,
                                     fill = yellow,
                                     text = red}}
\Vertex[x=-3, y=0, L={${1\ \atop32}, [000],\atop [-1,1,0],123$}]{x1}
\Edge[label={$s_1$},style={post,in=-10,out=190},color=blue](x2)(x1)
\Vertex[x=-2, y=2, L={${1\ \atop32}, [001],\atop [1,0,\alpha -1],231$}]{y1}
\Vertex[x=-4, y=4, L={${1\ \atop32}, [010],\atop [1,\alpha -1,0],213$}]{y2}
\Vertex[x=-6, y=6, L={${1\ \atop32}, [100],\atop [\alpha-1,1,0],123$}]{y3}
\Vertex[x=2, y=2, L={${2\ \atop31}, [001],\atop [-1,0,\alpha +1],231$}]{z1}
\Vertex[x=4, y=4, L={${2\ \atop31}, [010],\atop [-1,\alpha +1,0],213$}]{z2}
\Vertex[x=6, y=6, L={${2\ \atop31}, [100],\atop [\alpha+1,-1,0],123$}]{z3}
\Edge[label={$s_2$}](y1)(y2)
\Edge[label={$s_1$}](y2)(y3)
\Edge[label={$s_2$}](z1)(z2)
\Edge[label={$s_1$}](z2)(z3)
\Edge[label={$\Psi$}](x1)(y1)
\Edge[label={$\Psi$}](x2)(z1)

\Vertex[x=-2, y=8, L={${1\ \atop32}, [011],\atop [0,\alpha -1,\alpha
+1],312$}]{xy1}
\Vertex[x=2, y=8, L={${2\ \atop31}, [011],\atop [0,\alpha +1,\alpha
-1],312$}]{xz1}
\Vertex[x=-4, y=10, L={${1\ \atop32}, [101],\atop [\alpha-1,0\alpha
+1],132$}]{xy2}
\Vertex[x=4, y=10, L={${2\ \atop31}, [101],\atop [\alpha+1,0\alpha
-1],132$}]{xz2}
\Vertex[x=-6, y=12, L={${1\ \atop32}, [002],\atop
[1,0,2\alpha-1],231$}]{xy3}
\Vertex[x=6, y=12, L={${2\ \atop31}, [002],\atop
[-1,0,2\alpha+1],231$}]{xz3}
\Edge[label={$\Psi$}](y1)(xy1)
\Edge[label={$\Psi$}](y2)(xy2)
\Edge[label={$\Psi$}](y3)(xy3)
\Edge[label={$\Psi$}](z1)(xz1)
\Edge[label={$\Psi$}](z2)(xz2)
\Edge[label={$\Psi$}](z3)(xz3)

\Edge[label={$s_2$},style={post,in=-10,out=190},color=blue](xz1)(xy1)
\Edge[label={$s_1$}](xy1)(xy2)
\Edge[label={$s_1$}](xz1)(xz2)
\Edge[label={$\Psi$}](y3)(xy3)

\Vertex[x=-2, y=12, L={${1\ \atop32}, [110],\atop [\alpha-1,\alpha
+1,0],123$}]{xy4}
\Vertex[x=2, y=12, L={${2\ \atop31}, [110],\atop [\alpha+1,\alpha
-1,0],123$}]{xz4}
\Edge[label={$s_2$}](xy2)(xy4)
\Edge[label={$s_2$}](xz2)(xz4)
\Edge[label={$s_1$},style={post,in=-10,out=190},color=blue](xz4)(xy4)

\Vertex[x=-4, y=14, L={${1\ \atop32}, [020],\atop
[1,2\alpha-1,0],213$}]{xy5}
\Vertex[x=-2, y=16, L={${1\ \atop32}, [200],\atop
[2\alpha-1,1,0],123$}]{xy6}
\Vertex[x=4, y=14, L={${2\ \atop31}, [020],\atop
[-1,2\alpha+1,0],213$}]{xz5}
\Vertex[x=2, y=16, L={${2\ \atop31}, [200],\atop
[2\alpha+1,-1,0],123$}]{xz6}

\Edge[label={$s_2$}](xz3)(xz5)
\Edge[label={$s_1$}](xz5)(xz6)
\Edge[label={$s_2$}](xy3)(xy5)
\Edge[label={$s_1$}](xy5)(xy6)
\end{tikzpicture}
\caption{\label{G21} The f\/irst vertices of the graph $G_{21}$ where we
omit to write the vertex $\varnothing$ and the associated arrows.}
\end{figure}
\end{Example}

For a reverse standard tableau $\tau$ of shape $\lambda$, a
partition of $N$, let
\def\inv{{\rm inv}}
\def\rw{{\rm rw}}
\def\cl{{\rm cl}}
\[
 \inv(\tau)= \#\{(i,j):1\leq i<j\leq N, \rw(i,\tau)>\rw(j,\tau)\},
\]
where $\rw(i,\tau)$ is the row of $\tau$ containing $i$ (also we denote the column containing $i$
by $\cl(i,\tau)$). Then a correct jump
from $\tau$ to $\tau'$ implies $\inv(\tau')=\inv(\tau)+1$ (the entries
$\sigma[i]$ and $\sigma[i+1]=\sigma[i]+1$ are interchanged in $\tau$ to
produce $\tau'$). Thus the number of correct jumps in a path from the
root to $(\tau,\zeta,v,\sigma)$ equals $\inv(\tau)-\inv(\tau_\lambda)$.

So we consider the number of steps in a path from $0^N$ to
$v$;\footnote{The other components of the label
of the vertices are omitted here.} recall that one step
links~$v$ to~$v'$ where either $v[i] < v[i + 1]$ and $v' = vs_i$ or $v'
= v\Psi$.

For $x\in\Z$ (or $\R$) let $\epsilon(x) := \frac12(|x| + |x + 1|-1)$,
then $\epsilon(x) = x$ for $x\geq0$,
$\epsilon(x) = 0$ for $-1\leq x\leq0$, and $\epsilon(x) =-x-1$ for
$x\leq-1$.

There is a symmetry relation:
$\epsilon(x) = \epsilon(-x-1)$.

\begin{Definition}
For $v\in \N^N$ let $|v|:=\sum\limits_{i=1}^Nv[i]$ and set
\[
 {\rm S}(v):=\sum_{1\leq i<j\leq N}\epsilon(v[i]-v[j]).
\]
\end{Definition}
The above formula can be written as
\[
{\rm S}(v)=\frac12\sum_{1\leq i<j\leq
N}(|v[i]-v[j]|+|v[i]-v[j]+1|)-{N(N-1)\over4}.
\]

\begin{Proposition}\label{nbsteps}
The number of steps in any path joining $0^N$ to $v$ equals $|v|+{\rm
S}(v)$.
\end{Proposition}

\begin{proof}
The base point satisf\/ies $|0^N|=0$ and $|{\rm S}(0^N)|=0$.  Since
jumps do not modify $v$, consider a step
of the form $v' = vs_m$, then $|v'| = |v|$, $v [m + 1]-v [m]\geq1$ and
${\rm S} (v')-{\rm S} (v)$
involves only the pair $(m,m + 1)$ in the sum over all pairs $(i, j)$,
$1\leq i < j\leq N$. Indeed
\begin{gather*}
{\rm S}(v')-{\rm S}(v) = \epsilon(v[m+1]-v[m])-\epsilon(v[m]-v[m+1])\\
\phantom{{\rm S}(v')-{\rm S}(v)}{} =  (v[m+1]-v[m])-(-v[m]+v[m+1]-1)= 1.
\end{gather*}
It remains to show that ${\rm S}(v\Psi)={\rm S}(v)$ (because
$|v\Psi|=|v|+1$). Note $(v\Psi)[N]=v[1]+1$. Then
\begin{gather*}
{\rm S}(v)-{\rm  S}(v\Psi)= \sum_{j=2}^N\epsilon(v[1]-v[j])-\sum_{i=2}^N\epsilon(v[i]-v[1]-1)\\
\phantom{{\rm S}(v)-{\rm  S}(v\Psi)}{} = \sum_{j=2}^N\epsilon(v[1]-v[j])-\epsilon(v[j]-v[1]-1) =0.
\end{gather*}
This completes the proof.
\end{proof}

As a straightforward consequences, Proposition~\ref{nbsteps} implies
\begin{Corollary}\label{lengthpath}
All the paths
 joining two given vertices in $G_\lambda$
have the same length.
\end{Corollary}

This suggests that some properties could be shown by induction on the
common length of all the  all the paths
joining two given vertices.

 For a given 4-tuple $(\tau,\zeta,v,\sigma)$ the values of
$\zeta$ and $\sigma$ are determined by those of $\tau$ and $v$, as shown by
the following proposition.

\begin{Proposition}
If $(\tau,\zeta,v,\sigma)$ is a vertex in $G_\lambda$, then
$\sigma=\sigma_v$ and $\zeta[i]=v_{i}\alpha+CT_\tau[\sigma[i]]$. We will
set $\zeta_{v,\tau}:=\zeta$.
\end{Proposition}

\begin{proof}
We prove the result by induction on the length $k$ of a path
$(a_1,\dots,a_k)$ (from Corollary~\ref{lengthpath} all the paths have
the same length) from the root to $(\tau,\zeta,v,\sigma)$ and set
\[(\tau',\zeta',v',\sigma')=\big(\tau_\lambda,\CT_{\tau_\lambda},0^N,[1,\dots,N]\big)a_{1}\cdots
a_{k-1}.\]
Suppose that $a_k=\Psi$ is the af\/f\/ine operation. More precisely,
$\tau=\tau'$, $\zeta=\zeta'\Psi^\alpha$, $v=v'\Psi$ and
$\sigma=\sigma'[2,\dots,N,1]$. Using the induction hypothesis one has
$\sigma'=\sigma_{v'}$ and
$\zeta'[i]={v'[i]}\alpha+\CT_\tau[\sigma_{v'}[i]] $. Hence,
Proposition~\ref{Psi2sigma_v} gives
$\sigma=\sigma_{v'\Psi}=\sigma_v$. Suppose that $i<N$ then
\[\zeta[i]=\zeta'[i+1]={v'[i+1]}\alpha+\CT_\tau[\sigma_{v'}[i+1]]=
{v[i]}\alpha+\CT_\tau[\sigma_{v}[i]].
\]
If $i=N$ then again
\[\zeta[i]=\zeta'[1]+a=(v'[1]+1)\alpha+\CT_\tau[\sigma_{v'}[1]]=
v[N]\alpha+\CT_\tau[\sigma_{v}[N]].\]
Suppose now that $a_k$ is not an af\/f\/ine operation. Using the induction
hypothesis one has $\sigma'=\sigma_{v'}$ and
$\zeta'[j]={v'[j]}\alpha+\CT_\tau[\sigma_{v'}[j]] $ for each $j$. If
$a_k=s_i$ is a step then
$\tau=\tau'$, $\zeta=\zeta's_i$, $v=v's_i$ and $\sigma=\sigma' s_i$.
Hence, Proposition~\ref{si2sigma_v} gives $\sigma=\sigma_{v's_i}=\sigma_v$. If $j\neq
i,i+1$ then one has $\zeta[j]=\zeta'[j]$, $v[j]=v'[j]$ and
$\sigma[j]=\sigma'[j]$, hence
$\zeta[j]={v[j]}\alpha+\CT_\tau[\sigma_{v}[j]] $. If $j=i$ then one has
$\zeta[i]=\zeta'[i+1]$, $v[i]=v'[i+1]$ and $\sigma[i]=\sigma'[i+1]$, and
again the result is straightforward. And similarly when $j=i+1$ one
f\/inds the correct value for $\zeta[i+1]$.

Suppose now that $a_k=s_i$ is a jump. That is
$\tau=\tau'^{(\sigma[i],\sigma[i+1])}$, $\zeta=\zeta's_i$, $v=v'$ and
$\sigma=\sigma'$. Straightforwardly, $\sigma=\sigma_{v'}=\sigma_v$ and
if $j\neq i,i+1$ then \[\zeta[j]=\zeta'[j]
={v'[j]}\alpha+\CT_{\tau'}[\sigma_{v'}[j]]={v[j]}\alpha+\CT_\tau[\sigma_{v}[j]].\]
Suppose that $j=i$, since $a_k$ is a jump $v'[i]=v'[i+1]$ and
\[\zeta[i]=\zeta'[i+1]
={v'[i+1]}\alpha+\CT_{\tau'}[\sigma_{v'}[i+1]]={v[i]}\alpha+\CT_\tau[\sigma_{v}[i]].\]
Similarly, when $j=i+1$, one obtains the correct value for $\zeta[j]$.
This ends the proof.
 \end{proof}

\begin{Example}
Consider the RST $\tau=\begin{array}{cccc} 3\\7&4&1\\8&6&5&2\end{array}$
and the vector $v=[6,2,4,2,2,3,1,4]$. One has
$\sigma_v=[1,5,2,6,7,4,8,3]$ and $\CT_\tau=[1,3,-2,0,2,1,-1,0]$ and then
\[
\zeta_{v,\tau}=[6\alpha+1,2\alpha+2,4\alpha+3,2\alpha+1,2\alpha-1,3\alpha,\alpha,4\alpha-2].
\]
Hence, the $4$-tuple 
\begin{gather*}
\mbox{\footnotesize $\left(\begin{array}{@{\,}cccc@{}} 3\\7&4&1\\8&6&5&2\end{array},
[6\alpha+1,2\alpha+2,4\alpha+3,2\alpha+1,2\alpha-1,3\alpha,\alpha,4\alpha-2],[6,2,4,2,2,3,1,4],[1,5,2,6,7,4,8,3]\!\right)\!$}
\end{gather*}
labels a vertex of $G_{431}$.
\end{Example}

As a consequence,
\begin{Corollary}
Let $(\tau,v)$ be a pair constituted with a RST $\tau$ of shape
$\lambda$ $($a partition of $N)$ and a vector $v\in\N^N$. Then there
exists a unique vertex in $G_\lambda$ labeled by a $4$-tuple of the form
$(\tau,\zeta,v,\sigma)$. We will denote ${\cal
V}_{\tau,\zeta,v,\sigma}:=(\tau,v)$.
\end{Corollary}
We point out that all the information can be retrieved from the spectral
vector $\zeta$~-- the coef\/f\/icients of $\alpha$ give $v$, the rank
function of $v$ gives $\sigma$, and the constants in the spectral vector
give the content vector which does uniquely determine the RST $\tau$.
\begin{Definition}
We def\/ine the subgraph $G_\tau$ as the graph obtained from $G_\lambda$
by erasing all the vertices labeled by RST other than $\tau$ and the
associated arrows. Such a graph is connected.
\end{Definition}

Note that the graph $G_\lambda$ is the union of the graphs $G_\tau$
connected by jumps. Furthermore, if~$G_\tau$ and~$G_{\tau'}$ are
connected by a succession of jumps then there is no step from~$G_{\tau'}$ to~$G_\tau$.
\begin{Example}

In Fig.~\ref{G21}, the graph $G_{21}$ is constituted with the two graphs
$G_{1\ \atop32}$ and $G_{2\ \atop 31}$ connected by jumps (in blue).
\end{Example}

\section{Vector-valued polynomials}\label{s3VectorValued}

\subsection{About the Young seminormal representation of the symmetric group}\label{ss3.1Young}

We consider the space $V_\lambda$ spanned by reverse tableaux of shape
$\lambda$ and the (Young) action of the symmetric group as def\/ined by Murphy\footnote{The Young seminormal representation was def\/ined in Young's last papers but himself apparently underestimated the importance of the construction. G.~Murphy rediscovered it when reading the Jucys' paper~\cite{Jucys}. See~\cite{OV} for more details about the seminormal representation and its relation with the notions of \emph{Gelfand--Tsetlin basis}.} in~\cite{Murph}
by
\begin{equation}\label{Murphy}
\tau s_i=\left\{
\begin{array}{ll}
  b_\tau[i]\tau&\mbox{if} \ b_\tau[i]^2=1,\\
b_\tau[i]\tau+\tau^{(i,i+1)}&\mbox{if} \ 0<b_\tau[i]\leq \frac12,\\
b_\tau[i]\tau+(1-b_\tau[i]^2)\tau^{(i,i+1)}& \mbox{otherwise},
\end{array}
\right.
\end{equation}
\noindent where $b_\tau[i]:=\frac1{\CT_\tau[i]-\CT_\tau[i+1]}$. Note
that when $|b_\tau[i]|<1$, $\tau^{(i,i+1)}$ is always a reverse standard
tableau when $\tau$ is a reverse standard tableau.

Murphy showed \cite{Murph} that the RST are the simultaneous
eigenfunctions of the Jucys--Murphy elements:
\[
\omega_i=\sum_{j=i+1}^N s_{ij},
\]
where $s_{ij}$ denotes the transposition exchanging $i$ and $j$.
More precisely:
\begin{Proposition}\label{eigenMurph}
\[ \tau\omega_i=\CT_\tau[i]\tau.\]
\end{Proposition}

As usual, a polynomial representation for the Murphy action on the RST
can be computed through the Yang--Baxter graph. We start from
$\tau_\lambda$ and we construct the associated polynomial  in the
variables $t_1,\dots, t_N$:
\[P_{\tau_\lambda}=\prod_{i}\prod_{k>l\atop (i,k),(i,l)\in
\lambda}(t_{\tau_\lambda[i,k]}-t_{\tau_\lambda[i,l]}),\]
where $\tau[i,j]$ denotes the integer belonging at the column $i$ and
the row $j$ in $\tau$.
Such a polynomial is a simultaneous eigenfunction of the Jucys--Murphy
idempotents:
\[
 P_{\tau_\lambda}\omega_i = \CT_{\tau_\lambda}[i]P_{\tau_\lambda}.
\]
Suppose that $P_\tau$ is the polynomial associated to $\tau$. Suppose
also that $0<b_\tau[i]<1$. Hence, the polynomial $P_{\tau^{(i,i+1)}}$ is
obtained from the polynomial $P_\tau$ by acting with $s_i-b_\tau[i]$
(with the standard action of the transposition $s_i$ on the variables
$t_j$).

\begin{Example}
\begin{gather*}
P_{21\atop 43} =P_{31\atop 42}\big(s_2-\tfrac12\big) = (t_3-t_4)(t_1-t_2)\big(s_2-\tfrac12\big)\\
 \hphantom{P_{21\atop 43}}{} =  t_1 t_2-\tfrac12  t_4  t_1-\tfrac12  t_3  t_2+ t_4  t_3-\tfrac12
t_4  t_2-\tfrac12  t_3  t_1
\end{gather*}
\end{Example}

Let us remark that in \cite{YoungLasc}, Lascoux simplif\/ied the Young
construction by having recourse to the covariant algebra (of $\S_N$)
$\C[x_1,\dots,x_N]/{\S{\rm ym}_+}$ where ${\S{\rm ym}_+}$ is the ideal generated
by symmetric functions without constant terms. Note that the covariant
algebra is isomorphic to the regular representation $\C[\S_N]$. In the aim to adapt
his construction to our notations, we replace each polynomial with its
dominant monomial represented by the vectors of its exponents. The
vector associated to the root of the graph is the vector exponent of the
leading monomial in the product of the Vandermonde determinants associated
to each column and is obtained by putting the number of the
 row minus $1$ in the corresponding entry.

\begin{Example}
The vector associated to ${4\ \ \atop {52\ \atop 631}}$ is $[010210]$.
\end{Example}
In fact, the covariant algebra being isomorphic to the regular
representation of $\S_N$, the computation of the polynomials is
completely encoded by the action of the symmetric group on the leading
monomials, as shown in the following example. Observe that we do not replace the representation by the orbit of the leading monomial (since the space generated by the orbit is in general bigger), but we consider the projection which completely determines the elements.
\begin{Example}\rm Consider the RST of shape $221$, one has
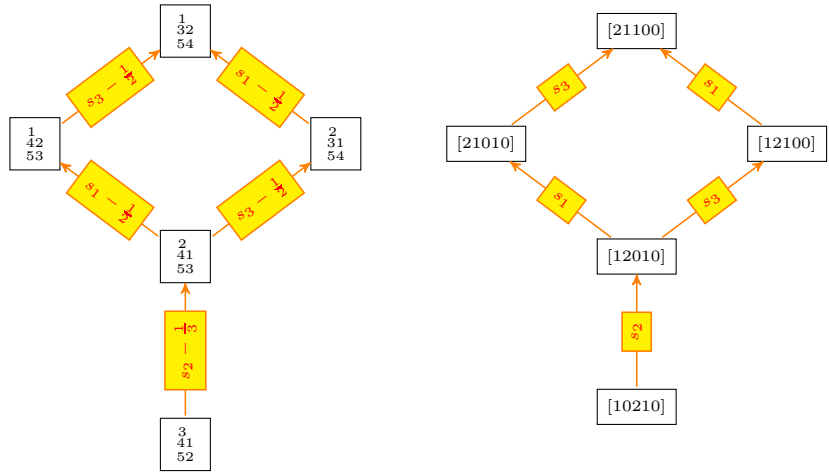
\begin{center}\begin{tikzpicture}%

\GraphInit[vstyle=Shade]
    \tikzstyle{VertexStyle}=[shape = rectangle,
                             draw
]

\SetUpEdge[lw = 1.5pt,
color = orange,
 labelcolor = gray!30,
 labelstyle = {draw,sloped},
 style={post}
]

\tikzset{LabelStyle/.style = {draw,
                                     fill = yellow,
                                     text = red}}

\Vertex[x=0, y=-0.5, L={\tiny $3\ \atop {41\atop 52}$}]{a1}
\Vertex[x=0, y=2, L={\tiny $2\ \atop {41\atop 53}$}]{a2}
\Vertex[x=-2, y=3.5, L={\tiny $1\ \atop {42\atop 53}$}]{a3}
\Vertex[x=2, y=3.5, L={\tiny $2\ \atop {31\atop 54}$}]{a4}
\Vertex[x=0, y=5, L={\tiny $1\ \atop {32\atop 54}$}]{a5}

\Vertex[x=6, y=0, L={\tiny $[10210]$}]{b1}
\Vertex[x=6, y=2, L={\tiny $[12010]$}]{b2}
\Vertex[x=4, y=3.5, L={\tiny $[21010]$}]{b3}
\Vertex[x=8, y=3.5, L={\tiny $[12100]$}]{b4}
\Vertex[x=6, y=5, L={\tiny $[21100]$}]{b5}

\Edge[label={\tiny $s_2-\frac13$}](a1)(a2)
\Edge[label={\tiny$s_1-\frac12$}](a2)(a3)
\Edge[label={\tiny$s_3-\frac12$}](a2)(a4)
\Edge[label={\tiny$s_3-\frac12$}](a3)(a5)
\Edge[label={\tiny$s_1-\frac12$}](a4)(a5)

\Edge[label={\tiny $s_2$}](b1)(b2)
\Edge[label={\tiny$s_1$}](b2)(b3)
\Edge[label={\tiny$s_3$}](b2)(b4)
\Edge[label={\tiny$s_3$}](b3)(b5)
\Edge[label={\tiny$s_1$}](b4)(b5)

\end{tikzpicture}

\end{center}

For instance, one has
\begin{gather*}
P_{1\ \atop{42\atop 53}} =
{\red{t_1^{2} t_4 t_2}}
+\tfrac12 t_4 t_2^{2} t_3
+\tfrac12 t_2^{2}t_5 t_1
-\tfrac12 t_2^{2} t_5 t_3
+\tfrac12 t_4 t_3^{2} t_1
+\tfrac12 t_4^{2} t_1 t_3
 + t_1^{2} t_5 t_3  - t_1^{2} t_5 t_2
 +\tfrac12 t_5^{2} t_4 t_3
 \\
\hphantom{P_{1\ \atop{42\atop 53}} = }{}
 +\tfrac12 t_5 t_4^{2} t_2
 +\tfrac12 t_5^{2} t_1 t_2
-\tfrac12 t_4 t_2^{2} t_1
 -\tfrac12 t_5^{2} t_4 t_2
 +\tfrac12 t_3^{2} t_5 t_2
 -\tfrac12 t_4^{2} t_1 t_2
 -t_1^{2}t_4t_3
 -\tfrac12 t_3^{2} t_5 t_1
  \\
\hphantom{P_{1\ \atop{42\atop 53}} = }{}
 -\tfrac12 t_5 t_4^{2} t_3
 - \tfrac12 t_4 t_3^{2} t_2
 -\tfrac12 t_5^{2} t_1 t_3,
\end{gather*}
whose leading monomial is $t_1^{2}t_2 t_4$.
\end{Example}

 From the construction, the leading monomial
of $P_\tau$ is the product of all the $t_i^{\rw(i,\tau)-1}$.
For example, the leading monomial in $P_{51\ \ \atop{732\ \atop9864}}$
is $t_1^2t_2t_3t_5^2t_7$.

\subsection{Def\/inition and dominance properties of vector-valued
polynomials}\label{ss3.2Def}

Consider the space \[M_N={\rm span}_\C\big\{x_1^{v[1]}\cdots
x_N^{v[N]}\otimes\tau\ :\ v\in\N^N, \tau\in{\rm Tab}_\lambda,
\lambda\vdash N\big\},\] where ${\rm Tab}(N)$ denotes the set of the reverse
standard tableaux on $\{1,\dots,N\}$.

This space splits into a direct sum $M_N=\bigoplus_{\lambda\vdash
N}M_\lambda$, where \[
M_\lambda={\rm span}_\C\big\{x_1^{v[1]}\cdots
x_N^{v[N]}\otimes\tau\,|\, v\in\N^N, \tau\in{\rm Tab}_\lambda\big\}.
\]
The algebra $\C[\S_N]\otimes\C[\S_N]$ acts on these spaces by commuting
the vector of the powers on the variables on the left component and the
action on the tableaux def\/ined by Murphy (equation~(\ref{Murphy})) on the
right component.

\begin{Example}
\[
x_1^3x_2^1\otimes {\begin{array}{cc}2\\3&1\end{array}}(s_2\otimes s_1)
=\frac12x_1^3x_3^1\otimes
{\begin{array}{cc}2\\3&1\end{array}}+x_1^3x_3^1\otimes
{\begin{array}{cc}1\\3&2\end{array}}.
\]
\end{Example}
For simplicity we will denote $x^v=x_1^{v[1]}\cdots x_N^{v[N]}$ and
$x^{v,\tau}:=x^v\otimes \tau\sigma_v$.
By abuse of language~$x^{v,\tau}$ will be referred to as a
polynomial.
Note that the space $M_\lambda$ is spanned by the set of polynomials
 \[{\cal M}_\lambda:=\big\{ x^{v,\tau} : v\in\N^N, \tau\in{\rm
Tab}_\lambda\big\},\]
which can be naturally endowed with the  order
$\lhd$ def\/ined by
\[
x^{v,\tau}\tau\lhd x^{v',\tau'}\mbox{ if\/f } v\lhd v',
\]
with $v\lhd v'$ means that $v^{+}\prec {v'}^{+}$ or
$v^{+}={v'}^{+}$ and $v\prec v'$, where
$\prec$ denotes the classical dominance order on
vectors:
\[
v\preceq v' \quad \mbox{if\/f} \quad \forall\,  i, v[1]+\cdots+v[i]\leq v'[1]+\cdots+v'[i].
\]
\begin{Example}\qquad
\begin{enumerate}\itemsep=0pt
\item $x^{031,{\mbox{\tiny$\begin{array}{cc}2\\3&1\end{array}$}}}\lhd x^{310,{\mbox{\tiny
$\begin{array}{cc}1\\3&2\end{array}$}}}$ since $031\prec 310$.
\item $x^{220,{\mbox{\tiny $\begin{array}{cc}2\\3&1\end{array}$}}}\lhd x^{301,
{\mbox{\tiny $\begin{array}{cc}1\\3&2\end{array}$}}}$ since $220\prec 310$.
\item The polynomials $x^{031,{\mbox{\tiny
$\begin{array}{cc}2\\3&1\end{array}$}}}$ and $x^{031,{\mbox{\tiny
$\begin{array}{cc}1\\3&2\end{array}$}}}$ are not comparable.
\end{enumerate}
\end{Example}

The partial order $\unlhd$ will provide us a relevant dominance notion.
\begin{Definition}
The monomial $x^{v,\tau}$ is the {\it leading monomial} of a polynomial
$P$ if and only if $P$ can be written as
\[
 P=\alpha_v x^{v, \tau}+\sum_{x^{v',\tau'}\lhd
x^{v,\tau}}\alpha_{v',\tau'}x^{v',\tau'}
\]
with $\alpha_v\neq 0$.
\end{Definition}

As in~\cite{KS}, we def\/ine $\Psi:=(\theta \otimes \theta)
x_N$, with $\theta=s_1s_2\cdots s_{N-1}$. The following proposition describes the transformation properties of leading monomials with respect to the $s_i$ and $\Psi$.
\begin{Proposition}\label{dominance}
Suppose that
  $ x^{v,\tau}$ is the leading monomial in $P$ then
\begin{enumerate}\itemsep=0pt
\item[$1.$] If $v[i]< v[i+1]$ then
 $x^{v s_{i}, \tau}$ is the leading monomial in $P(s_{i}\otimes
s_{i})$.
Its leading monomial is $x^{\tilde v}\otimes (\tau.s_{ij})$, where
$x^{\tilde v}$ is the dominant term in $\partial_{ij}x^{\tilde v}$.
\item[$2.$]
$x^{v\Psi,\tau}$ is the leading monomial in $P\Psi$.
\end{enumerate}
\end{Proposition}

\subsection[Dunkl and Cherednik-Dunkl operators for vector-valued polynomials]{Dunkl and Cherednik--Dunkl operators for vector-valued polynomials}\label{ss3.3Dunkl}

We def\/ine the Dunkl operators
\[
{\cal D}_i:={\partial\over\partial x_i}\otimes 1+\frac1\alpha\sum_{i\neq
j}\partial_{ij}\otimes s_{ij},
\]
where $s_{ij}$ denotes the transposition which exchanges $i$ and $j$ and
\[
 \partial_{ij}=(1-s_{ij})\frac1{ x_i-x_j}
\]
is the divided dif\/ference.

 This def\/inition is the same as in \cite{Dunkl}, but our operators
act on their left.
One has
\begin{Lemma}\label{HD}
If ${\cal D}_i$ denotes the Dunkl operator, one has
\[
 (s_i\otimes s_i){\cal D}_i={\cal D}_{i+1}(s_i\otimes s_i).
\]
\end{Lemma}

\begin{proof} Straightforward from the def\/inition of
${\cal D}_i$ and the equalities
\begin{gather*}
s_is_{ij}=s_{i+1,j}s_i, \qquad s_i\partial_{ij}=\partial_{i+1j}s_i \qquad \mbox{and}\qquad
s_i{\partial\over\partial x_i}={\partial\over\partial x_{i+1}}s_i.\tag*{\qed}
\end{gather*}
 \renewcommand{\qed}{}
\end{proof}

The Cherednik--Dunkl operators are pairwise commuting
operators  def\/ined by \cite{Dunkl}
\[
{\cal U}_i:=x_i{\cal D}_i-\frac1\alpha\sum_{j=1}^{i-1}s_{i,j}\otimes
s_{i,j}.
\] We do not repeat the proof of the commutation $[{\cal U}_i,{\cal
U}_j]=0$ which can be found in~\cite{Dunkl}. But, as we will see in the
next section, this property is not used to prove the existence of the
vector-valued Jack polynomials.

One has
\begin{Lemma}\label{HU1}\qquad
\begin{enumerate}\itemsep=0pt
\item[$1.$]
$(s_i\otimes s_i) {\cal U}_i={\cal U}_{i+1} (s_i\otimes s_i)+\frac1\alpha
$.
\item[$2.$] $
 (s_i\otimes s_i) {\cal U}_j={\cal U}_{j} (s_i\otimes s_i)$, $j\neq i,i+1$.
\item[$3.$] $ (s_i\otimes s_i) {\cal U}_{i+1}={\cal U}_{i} (s_i\otimes s_i)-\frac1\alpha$.
\end{enumerate}
\end{Lemma}

\begin{proof}
The three identities are of the same type. We prove only the f\/irst
one which
follows from the equalities
\begin{gather*}
(s_i\otimes s_i){\cal U}_i = (s_i\otimes s_i)x_i{\cal
D}_i-\frac1\alpha\sum_{j=1}^{i-1}s_{i,j}\otimes s_{i,j}\\
\hphantom{(s_i\otimes s_i){\cal U}_i}{} = \left(x_{i+1}{\cal
D}_{i+1}-\frac1\alpha\sum_{j=1}^{i-1}s_{i+1,j}\otimes
s_{i+1,j}\right)(s_i\otimes s_i) = {\cal U}_{i+1}(s_i\otimes s_i)+\frac1\alpha.\tag*{\qed}
\end{gather*}
\renewcommand{\qed}{}
\end{proof}

 The af\/f\/ine operator $\Psi$ has the following commutation properties
with the Dunkl operators:
\begin{Lemma}\label{DU2}\qquad
\begin{enumerate}
 \item[$1.$] ${\cal D}_{i+1}\Psi=\Psi{\cal D}_i+(\theta\otimes \theta)(
s_{i,N}\otimes s_{i,N})$, $i<1$.

 \item[$2.$] ${\cal D}_1\Psi=\Psi{\cal
D}_N+(\theta\otimes\theta)\left(\sum\limits_{j=1}^{N-1}(s_{N,j}\otimes
s_{N,j})-1\right)$.
\end{enumerate}

\end{Lemma}
As a consequence, one f\/inds.
\begin{Lemma}\label{HU2}
\[
\Psi{\cal U_i}=\cal U_{i+1}\Psi,\qquad i\neq N \qquad \mbox{and}\qquad \Psi
{\cal U}_N=\left({\cal U}_1+1\right)\Psi.
\]
\end{Lemma}

The action on the RST is given by

\begin{Lemma}\label{U_itab}
\[
 (1\otimes \tau) {\cal U}_i=\left(1+\frac1\alpha \CT_\tau[i]\right)(1\otimes \tau).
 \]
\end{Lemma}

\begin{proof} One has
\begin{gather*}
(1\otimes \tau) {\cal U}_i = (1\otimes \tau) x_i{\cal
D}_i-\frac1\alpha\sum_{j=1}^{i-1}(1\otimes\tau)( s_{i,j}\otimes s_{i,j})
 = (1\otimes\tau) \left(1+\frac1\alpha1\otimes \omega_i\right),
\end{gather*}
where $\omega_i:=\sum\limits_{j=i+1}^N(i\ j)$ denotes a Jucys--Murphy element.
Since the RST are eigenfunctions of the Jucys--Murphy elements and the
associated eigenvalues are given by the contents, the lemma follows.
\end{proof}

For convenience, def\/ine $\tilde \xi_i:=\alpha{\cal U}_i-\alpha$.
 From the preceding lemmas, one obtains
\begin{Proposition}\label{Cherednik_com}
\begin{gather*}
(s_i\otimes s_i) \tilde\xi_i=\tilde\xi_{i+1}(s_i\otimes
s_i)+1,\\ 
(s_i\otimes s_i)\tilde\xi_{i+1}=\tilde\xi_{i}(s_i\otimes
s_i)-1,\\ 
(s_i\otimes s_i)\tilde\xi_j=\tilde \xi_j (s_i\otimes s_i),\qquad j\neq
i,i+1,\\ 
\Psi\tilde \xi_i=\tilde \xi_{i+1}\Psi,\qquad i\neq N,\label{xiPsi}
\\
\Psi\tilde\xi_N=\big(\tilde\xi_1+\alpha\big)\Psi.
\end{gather*}
\end{Proposition}

\section{Nonsymmetric vector-valued Jack polynomials}\label{s4NonSym}

In this section we recover the construction, due to one of the authors~\cite{Dunkl}, of a basis of vector-valued polynomials $J_{v,\tau}$. This construction belongs to a large family of vector-valued Jack
polynomials associated to the complex ref\/lection groups $G(r,1,n)$
def\/ined by Grif\/feth~\cite{Gr2}. We will denote by $\zeta_{v,\tau}$ their
associated spectral vectors. We will see also that many properties of
this basis can be deduced from the Yang--Baxter structure.

\subsection[Yang-Baxter construction associated to $G_\lambda$]{Yang--Baxter construction associated to $\boldsymbol{G_\lambda}$}\label{ss4.1YB}

Let $\lambda$ be a partition and $G_\lambda$  be  the
associated graph. We construct the set of
polynomials
$\left(J_{\cal P}\right)_{\cal P \ {\rm path \ in} \ G_\lambda}$ using the
following recursive rules:
\begin{enumerate}\itemsep=0pt
 \item $J_{[]}:=(1\otimes\tau_{\lambda})$.
 \item If ${\cal P}=[a_1,\dots,a_{k-1},s_i]$ then
\[J_{\cal P}:=J_{[a_1,\dots,a_{k-1}]}\left(s_i\otimes
s_i+\frac1{\zeta[i+1]-\zeta[i]}\right),\]
where the vector $\zeta$ is def\/ined by
\[(\tau_\lambda,\CT_{\tau_\lambda},0^N,[1,2,\dots,N]) a_1 \dots
a_{k-1}=(\tau,\zeta,v,\sigma).\]
 \item If ${\cal P}=[a_1,\dots,a_{k-1},\Psi]$ then
\[J_{\cal P}=J_{[a_1,\dots,a_{k-1}]}\Psi.\]
\end{enumerate}

One has the following theorem.

\begin{Theorem}\label{thJack}
Let ${\cal P}=[a_0,\dots,a_k]$ be a path in $G_\lambda$ from the root to
$(\tau,\zeta,v,\sigma)$.
The polynomial $J_{\cal P}$ is a simultaneous eigenfunctions of the
operators $\tilde\xi_i$ whose leading monomial is $x^{v,\tau}$. Furthermore,
the eigenvalues of $\tilde\xi_i$ associated to $J_{\cal P}$ are equal to
$\zeta[i]$.

 Consequently $J_{\cal P}$ does not depend on the path, but only on
the end point $(\tau,\zeta,v,\sigma)$, and will be denoted by
$J_{v,\tau}$.
 The family $(J_{v,\tau})_{v,\tau}$ forms a basis of
$M_\lambda$ of simultaneous eigenfunctions of the Cherednik operators.

Furthermore, if $\cal P$ leads to $\varnothing$ then $J_{\cal P}=0$.
\end{Theorem}

\begin{proof}
We will prove the result by induction on the length $k$. If $k=0$ then
the result follows from Proposition~\ref{U_itab}.
Suppose now that $k>0$ and let
\[(\tau',\zeta',v',\sigma_{v'})=\big(\tau_\lambda,\CT_{\tau_\lambda},0^N,[1,\dots,N]\big)a_1\cdots
a_{k-1} .\]
By induction, $J_{[a_1,\dots,a_{k-1}]}$ is a simultaneous eigenfunctions
of the operators $\tilde\xi_i$ such that the associated vector of
eigenvalues is given by
\[J_{[a_1,\dots,a_{k-1}]}\tilde\xi_i=\zeta'[i]J_{[a_1,\dots,a_{k-1}]}\]
and the leading monomial is $x^{v',\tau'}$.

If  $a_k=\Psi$ is an af\/f\/ine arrow, then
$\tau=\tau'$, $\zeta=\zeta'.\Psi^\alpha$, $v=v'\Psi$,
$\sigma_v=\sigma_{v'}[2,\dots,N,1]$ and $J_{\cal
P}=J_{[a_1,\dots,a_{k-1}]}\Psi$. If $i\neq N$
 \begin{gather*}
 J_{\cal P}\tilde\xi_i = J_{[a_1,\dots,a_{k-1}]}\Psi\tilde\xi_i
=J_{[a_1,\dots,a_{k-1}]}\tilde\xi_{i+1}\Psi
= \zeta'[i+1]J_{\cal P}
= \zeta[i] J_{\cal P}.
\end{gather*}

If $i=N$ then,
 \begin{gather*}
 J_{\cal P}\tilde\xi_N =  J_{[a_1,\dots,a_{k-1}]}\Psi\tilde\xi_N
 = J_{[a_1,\dots,a_{k-1}]}\tilde(\xi_{1}+\alpha)\Psi
 = (\zeta'[1]+\alpha)J_{\cal P}
=\zeta[N] J_{\cal P}.
\end{gather*}
The   leading monomial is a consequence of Proposition~\ref{dominance}.

Suppose now that $a_k=s_i$ is a non af\/f\/ine arrow, then
$\zeta=\zeta's_i$, $v=v's_i$ and \[
J_{\cal P}=J_{[a_1,\dots,a_{k-1}]}\left(s_i\otimes
s_i+\frac1{\zeta'[i+1]-\zeta'[i]}\right).
\] If $j\neq i, i+1$ then
\begin{gather*}
 J_{\cal P}\tilde\xi_j =  J_{[a_1,\dots,a_{k-1}]}\left(s_i\otimes
s_i+\frac1{\zeta'[i+1]-\zeta'[i]}\right)\tilde\xi_j \\
\hphantom{J_{\cal P}\tilde\xi_j}{}
 = J_{[a_1,\dots,a_{k-1}]}\tilde\xi_{j}\left(s_i\otimes
s_i+\frac1{\zeta'[i+1]-\zeta'[i]}\right)
 =
\zeta'[j]J_{\cal P}
 = \zeta[j] J_{\cal P}.
\end{gather*}
If $j=i$ then
\begin{gather*}
 J_{\cal P}\tilde\xi_i =  J_{[a_1,\dots,a_{k-1}]}\left(s_i\otimes
s_i+\frac1{\zeta'[i+1]-\zeta'[i]}\right)\tilde\xi_i\\
\hphantom{J_{\cal P}\tilde\xi_i}{}
 = J_{[a_1,\dots,a_{k-1}]}\left(\tilde\xi_{i+1}(s_i\otimes
s_i)+1+\tilde\xi_i\frac1{\zeta'[i+1]-\zeta'[i]}\right)\\
\hphantom{J_{\cal P}\tilde\xi_i}{}
 =
J_{[a_1,\dots,a_{k-1}]}\left(\zeta'[i+1](s_i\otimes
s_i)+1+{\zeta'[i]\over\zeta'[i+1]-\zeta'[i]}\right)\\
\hphantom{J_{\cal P}\tilde\xi_i}{}
 =
\zeta'[i+1]J_{[a_1,\dots,a_{k-1}]}\left(s_i\otimes
s_i+{1\over\zeta'[i+1]-\zeta'[i]}\right) =\zeta[i] J_{\cal P}.
\end{gather*}
If $j=i+1$ then
\begin{gather*}
 J_{\cal P}\tilde\xi_{i+1} =  J_{[a_1,\dots,a_{k-1}]}\left(s_i\otimes
s_i+\frac1{\zeta'[i+1]-\zeta'[i]}\right)\tilde\xi_i\\
\hphantom{J_{\cal P}\tilde\xi_{i+1}}{}
 = J_{[a_1,\dots,a_{k-1}]}\left(\tilde\xi_{i}(s_i\otimes
s_i)-1+\tilde\xi_{i+1}\frac1{\zeta'[i+1]-\zeta'[i]}\right)\\
\hphantom{J_{\cal P}\tilde\xi_{i+1}}{}
=
J_{[a_1,\dots,a_{k-1}]}\left(\zeta'[i](s_i\otimes
s_i)-1+{\zeta'[i+1]\over\zeta'[i+1]-\zeta'[i]}\right)\\
\hphantom{J_{\cal P}\tilde\xi_{i+1}}{}
=
\zeta'[i]J_{[a_1,\dots,a_{k-1}]}\left(s_i\otimes
s_i+{1\over\zeta'[i+1]-\zeta'[i]}\right)=\zeta[i+1] J_{\cal P}.
\end{gather*}
Let us examine the  leading
monomials. First, suppose that $a_k=s_i$ is a step then
$\tau=\tau'$ and $\sigma_v=\sigma_{v'}s_i$. From Proposition~\ref{dominance}, the leading monomial in $J_{\cal P}$ equals the leading
term in $x^{v'\tau'}\left(s_i\otimes
s_i+\frac1{\zeta'[i+1]-\zeta'[i]}\right)$
that is $x^{v's_i}\otimes(\tau'\sigma_{v'}s_i)=x^{v,\tau}$.

Suppose that $a_k=s_i$ is not a step and set
$Q:=x^{v',\tau'}\left(s_i\otimes s_i+{1\over
\zeta'[i+1]-\zeta'|i]}\right)$. One has
\[
 J_{\cal P}=Q+\sum_{x^{v'',\tau''} \lhd
x^{v,\tau}}\alpha_{v'',\tau''}x^{v'',\tau''}.
\]

If $a_k=s_i$ is a jump then
$\tau={\tau'}^{(\sigma_{v'}[i],\sigma_{v'}[i]+1)}$ and $\sigma_v=\sigma_v'$.
 But
\begin{gather*}
Q = x^{v's_i}\otimes(\tau'\sigma_{v'}s_i)+\frac1{\zeta'[i+1]-\zeta'[i]}x^{v'}\otimes(\tau'\sigma_{v'})\\
\phantom{Q}{} = x^{v}\otimes
(\tau's_{\sigma_{v'}[i]}\sigma_{v'})+\frac1{\zeta'[i+1]-\zeta'[i]}x^{v'}\otimes(\tau'\sigma_{v'})\\
\phantom{Q}{} = x^{v}\otimes
(\tau\sigma_v)+\left(b_\tau'[\sigma_{v'}[i]]+\frac1{\zeta'[i+1]-\zeta'[i]}\right)x^{v'}\otimes(\tau'\sigma_{v'})\\
\phantom{Q}{} =
x^{v,\tau}+\left(b_\tau'[\sigma_{v'}[i]]+\frac1{\zeta'[i+1]-\zeta'[i]}\right)x^{v',\tau'}.
\end{gather*}
But $
\zeta'[i]=\CT_{\tau'}[\sigma_{v'}[i]]
$ and $
\zeta'[i+1]=\CT_{\tau'}[\sigma_{v'}[i+1]]=\CT_{\tau'}[\sigma_{v'}[i]+1]
$, hence $b_{\tau'}[\sigma_{v'}[i]]=-\frac1{\zeta'[i+1]-\zeta'[i]}$. And
the leading monomial is $Q=x^{v,\tau}$ as expected.

This proves the f\/irst part of the theorem and that the family
$(J_{v,\tau})_{v,\tau}$ forms a basis of $M_\lambda$ of simultaneous
eigenfunctions of the Cherednik operators.

Finally, if $a_k=s_i$ is a fall, $Q$ is proportional to
$x^{v'}\otimes(\tau'\sigma_{v'})$ and then $J_{\cal P}$ is proportional
to $J_{[a_1,\dots,a_{k-1}]}$. But clearly, the two polynomials are
eigenfunction of the Cherednik operators with dif\/ferent eigenvalues
 from the cases $j=i$ and $j=i+1$. This proves that $J_{\cal
P}=0$.
\end{proof}

 As a consequence, we will consider the family of polynomials
$(J_{v,\tau})_{v,\tau}$ indexed by pairs $(v,\tau)$ where $v\in \N^N$ is
a weight and $\tau$ is a tableau.

\begin{Example}
 For $\lambda=21$,
 the f\/irst polynomials
$J_{v,\tau}$  are
displayed in Fig.~\ref{PG21}. The spectral vectors can be
read on Fig.~\ref{G21}.

\begin{figure}[h]
\centering

\begin{tikzpicture}%
\GraphInit[vstyle=Shade]
    \tikzstyle{VertexStyle}=[shape = rectangle,
                             draw
]

\Vertex[style={
shape=circle},x=3, y=0,
 L={$J_{000,{2\ \atop 31}}$}]{x2}

\SetUpEdge[lw = 1.5pt,
color = orange,
 labelcolor = gray!30,
 labelstyle = {sloped},
 style={post}
]

\tikzset{LabelStyle/.style = {draw,
                                     text = black}}

\Vertex[x=2, y=2, L={\tiny$J_{001,{2\ \atop 31}}$}]{z1}
\Vertex[x=4, y=4, L={\tiny$J_{010,{2\ \atop 31}}$}]{z2}
\Vertex[x=6, y=6, L={\tiny$J_{100,{2\ \atop 31}}$}]{z3}

\Edge[label={\tiny $s_2+\frac1{\alpha+1}$}](z1)(z2)
\Edge[label={\tiny $s_1+\frac1{\alpha+2}$}](z2)(z3)
\Edge[label={$\Psi$}](x2)(z1)

\Vertex[x=2, y=8, L={\tiny$J_{011,{2\ \atop 31}}$}]{xz1}
\Vertex[x=4, y=10, L={\tiny$J_{101,{2\ \atop 31}}$}]{xz2}
\Vertex[x=6, y=12, L={\tiny$J_{002,{2\ \atop 31}}$}]{xz3}
\Edge[label={$\Psi$}](z1)(xz1)
\Edge[label={$\Psi$}](z2)(xz2)
\Edge[label={$\Psi$}](z3)(xz3)

\Edge[label={\tiny $s_1+\frac1{\alpha+1}$}](xz1)(xz2)

\Vertex[x=2, y=12, L={\tiny$J_{110,{2\ \atop 31}}$}]{xz4}
\Edge[label={\tiny $s_2+\frac1{\alpha-1}$}](xz2)(xz4)

\Vertex[x=4, y=14, L={\tiny$J_{020,{2\ \atop 31}}$}]{xz5}
\Vertex[x=2, y=16, L={\tiny$J_{200,{2\ \atop 31}}$}]{xz6}

\Edge[label={\tiny $s_2+\frac1{2\alpha+1}$}](xz3)(xz5)
\Edge[label={\tiny$s_1+\frac1{2\alpha+2}$}](xz5)(xz6)

\Vertex[x=-3, y=0, L={\tiny$J_{000,{1\ \atop32}}$}]{x1}
\Vertex[x=-2, y=2, L={\tiny$J_{001,{1\ \atop32}}$}]{y1}
\Vertex[x=-4, y=4, L={\tiny$J_{010,{1\ \atop32}}$}]{y2}
\Vertex[x=-6, y=6, L={\tiny$J_{100,{1\ \atop32}}$}]{y3}
\Vertex[x=-2, y=8, L={\tiny$J_{011,{1\ \atop32}}$}]{xy1}
\Vertex[x=-4, y=10, L={\tiny$J_{101,{1\ \atop32}}$}]{xy2}
\Vertex[x=-6, y=12, L={\tiny$J_{002,{1\ \atop32}}$}]{xy3}
\Vertex[x=-2, y=12, L={\tiny$J_{110,{1\ \atop32}}$}]{xy4}
\Vertex[x=-4, y=14, L={\tiny$J_{020,{1\ \atop32}}$}]{xy5}
\Vertex[x=-2, y=16, L={\tiny$J_{200,{1\ \atop32}}$}]{xy6}

\Edge[label={\tiny $s_1+\frac12$}](x2)(x1)
\Edge[label={\tiny $s_2+\frac1{\alpha-1}$}](y1)(y2)
\Edge[label={\tiny $s_1+\frac1{\alpha-2}$}](y2)(y3)
\Edge[label={ $\Psi$}](x1)(y1)
\Edge[label={ $\Psi$}](y1)(xy1)
\Edge[label={$\Psi$}](y2)(xy2)
\Edge[label={$\Psi$}](y3)(xy3)
\Edge[label={\tiny $s_2+\frac12$}](xz1)(xy1)
\Edge[label={\tiny $s_1+\frac1{\alpha-1}$}](xy1)(xy2)
\Edge[label={$\Psi$}](y3)(xy3)
\Edge[label={\tiny $s_2+\frac1{\alpha+1}$}](xy2)(xy4)
\Edge[label={\tiny $s_1+\frac1{2}$}](xz4)(xy4)
\Edge[label={\tiny $s_2+\frac1{2\alpha-1}$}](xy3)(xy5)
\Edge[label={\tiny $s_1+\frac1{2\alpha-2}$}](xy5)(xy6)
\end{tikzpicture}
\caption{\label{PG21} First values of the polynomials $J_{v,\tau}$ for
$\lambda=21$ ($s_i$ means $s_i\otimes s_i$).}
\end{figure}
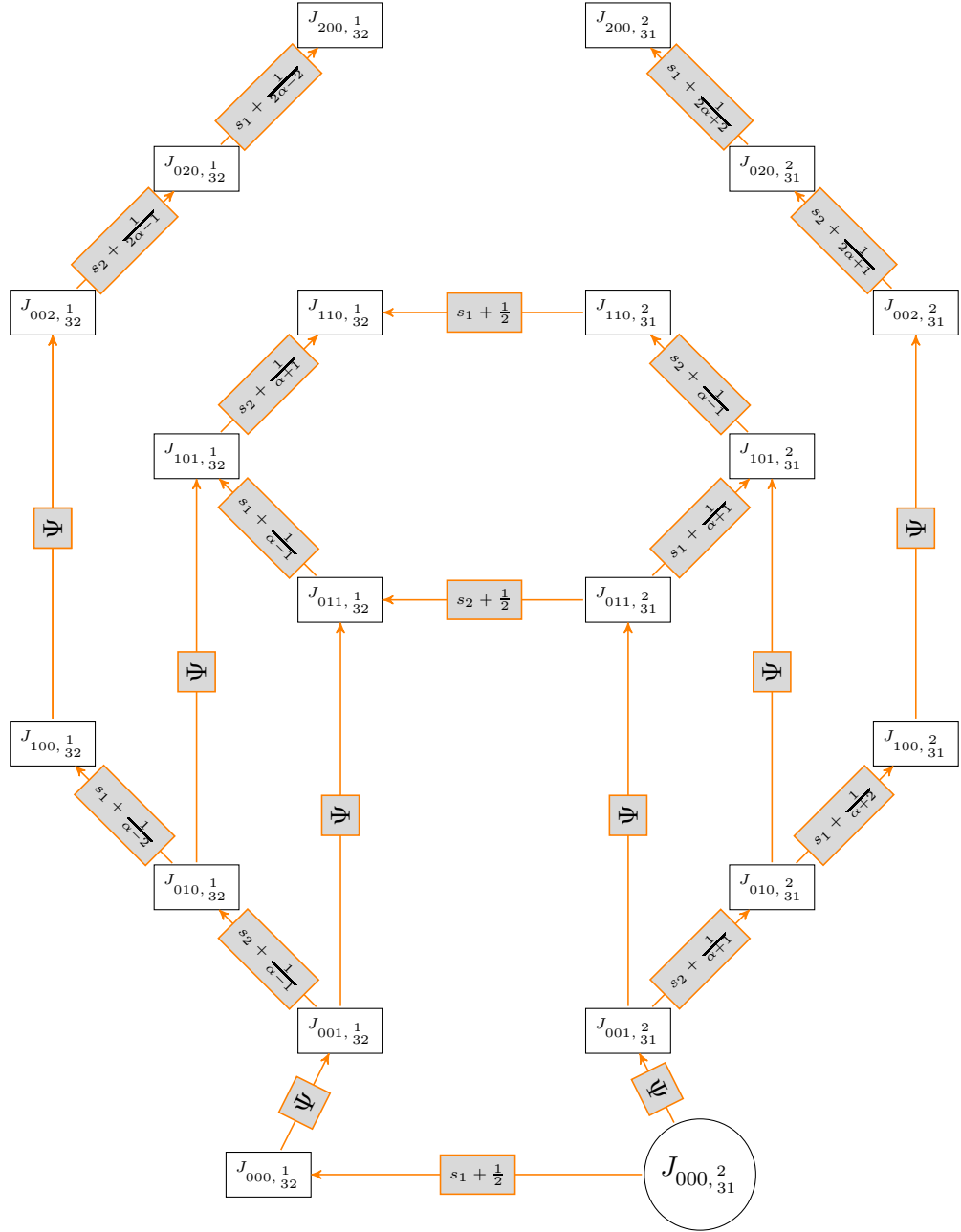
\end{Example}

Note that if $[a_1,\dots,a_{k-1}]$ leads to a vertex other than
$\varnothing$ and $[a_1,\dots,a_{k-1},s_i]$ leads to $\varnothing$, the last
part of Theorem~\ref{thJack} implies that $J_{[a_1,\dots,a_{k-1}]}$ is
symmetric or antisymmetric under the action of $s_i\otimes s_i$.

The recursive rules of this section f\/irst appeared in~\cite{Gr2}. The
Lemma~5.3 and  the Yang--Baxter graph  constitute
essentially what Grif\/feth called {\it calibration graph} in that paper.

\subsection[Partial Yang-Baxter-type construction associated to $G_\tau$]{Partial Yang--Baxter-type construction associated to $\boldsymbol{G_\tau}$}\label{ss4.2Gtau}

To compute an expression for a polynomial $J_{v,\tau}$ it suf\/f\/ices to
f\/ind the good path in the sub\-graph~$G_\tau$ as shown by the following
examples.
\begin{Example}
Consider $\tau=\begin{array}{cc}1\\3&2\end{array}$, Fig.~\ref{KG132}
explains how to obtains the values of $J_{v,{1\atop 32}}$ from the graph
$G_{1\atop32}$.
\begin{figure}[h]
\centering
\begin{tikzpicture}
\GraphInit[vstyle=Shade]
    \tikzstyle{VertexStyle}=[shape = rectangle,
draw
]

\Vertex[style={
shape=circle},x=3, y=0,
 L={$J_{000,{1\ \atop 32}}$}]{x2}
\SetUpEdge[lw = 1.5pt,
color = orange,
 labelcolor = gray!30,
 labelstyle = {sloped},
 style={post}
]

\tikzset{LabelStyle/.style = {draw,
                                     text = black}}
\tikzset{EdgeStyle/.style={post}}
\Vertex[x=2, y=2, L={\tiny$J_{001,{1\ \atop 32}}$}]{z1}
\Vertex[x=4, y=4, L={\tiny $J_{010,{1\ \atop 32}}$}]{z2}
\Vertex[x=6, y=6, L={\tiny$J_{100,{1\ \atop 32}}$}]{z3}

\Edge[label={\tiny $s_2+\frac1{\alpha-1}$}](z1)(z2)
\Edge[label={\tiny $s_1+\frac1{\alpha-2}$}](z2)(z3)
\Edge[label={$\Psi$}](x2)(z1)
\Vertex[x=2, y=8, L={\tiny$J_{011,{1\ \atop 32}}$}]{xz1}
\Vertex[x=4, y=10, L={\tiny$J_{101{1\ \atop 32}}$}]{xz2}
\Vertex[x=6, y=12, L={\tiny$J_{002,{1\ \atop 32}}$}]{xz3}
\Edge[label={$\Psi$}](z1)(xz1)
\Edge[label={$\Psi$}](z2)(xz2)
\Edge[label={$\Psi$}](z3)(xz3)
\Edge[label={\tiny $s_1+\frac1{\alpha-1}$}](xz1)(xz2)

\Vertex[x=2, y=12, L={\tiny$J_{110,{1\ \atop 32}}$}]{xz4}
\Edge[label={\tiny $s_2+\frac1{\alpha+1}$}](xz2)(xz4)
\Vertex[x=4, y=14, L={\tiny$J_{020,{1\ \atop 32}}$}]{xz5}
\Vertex[x=2, y=16, L={\tiny$J_{200,{1\ \atop 32}}$}]{xz6}

\Edge[label={\tiny $s_2+\frac1{2\alpha+1}$}](xz3)(xz5)
\Edge[label={\tiny$s_1+\frac1{2\alpha+2}$}](xz5)(xz6)

\end{tikzpicture}
\caption{\label{KG132} First values of the polynomials $J_{v,{1\ \atop
32}}$.}
\end{figure}

\end{Example}

\begin{Example}
For the trivial representation (i.e., $\lambda$ has a single
part), note that the Cherednik operators (in \cite{Las}) have
the same eigenspaces as the Cherednik--Dunkl operators ${\cal U}_i$ (in
\cite{Dunkl}). In the notations of \cite{Las}, $\xi_i$ reads
\[
\xi_i=\alpha x_i{\partial\over\partial x_i}+ \sum_{j=1\atop j\neq
i}^{N}\overline\pi_{ij}+(1-i),
\]
where \[\overline\pi_{ij}=\left\{\begin{array}{ll}x_i\partial_{ij} & \mbox{if} \ j<i, \\
x_j\partial_{ij} & \mbox{if} \ i<j, \end{array}\right.
\]
 where $\partial_{ij}$ denotes the divided dif\/ference  on the variables $x_i$ and $x_j$.
Noting that $x_i\partial_{ij}=\partial_{ij}x_i-1$,
$x_j\partial_{ij}=\partial_{ij}x_i-(ij)$ and $x_i{\partial\over\partial
x_i}={\partial\over\partial x_i}x_i-1$, one f\/inds
\[
 \xi_i=\alpha{\cal U}_i-(\alpha+N-1)=\tilde\xi_i-\left(
N-1\right).
\]
\end{Example}

\begin{Example}
Consider sign representation associated to the partition $[1^N]$. The set ${\rm Tab}_{[1^N]}$ contains a unique element
$\tau=\begin{array}{c}1\\\vdots\\N\end{array}$. Hence, we can omit $\tau$ when we write the polynomials of~$M_{[1^N]}$. One can see that the
corresponding Jack polynomials are equal to the standard ones for the
coef\/f\/icient $-\alpha$. Indeed, since $\tau s_{ij}=-\tau$ one has
\begin{gather*}
P{\cal D}_i\simeq(P\otimes\tau) {\cal D}_i = (P\otimes \tau)\left({\partial\over\partial x_i}+\frac1\alpha\sum_{i\neq j}\partial_{ij}\otimes s_{i,j}\right)\!
 = (P\otimes \tau)\left({\partial\over\partial x_i}-\frac1\alpha\sum_{i\neq j}\partial_{ij}\otimes 1\right).
\end{gather*}
Hence, the Cherednik--Dunkl operator ${\cal U}_i=x_i{\cal D}_i-\frac1\alpha\sum\limits_{j=1}^{i-1}s_{ij}\otimes s_{ij}$  acts on $M_{[1^N]}$ as the ope\-ra\-tor~${\cal U}_i$ acts on $M_{[N]}$ but for the parameter $-\alpha$.
\end{Example}

\begin{Example}
Let us explain the method on a bigger example: $J_{[0,0,2,1,1,0],\tau}$,
for $\tau:=\begin{array}{@{}c@{\,\,}c@{\,\,}c@{\,\,}c@{\,}} 4&3\\6&5&2&1\end{array}\!$.
First, we obtain the vector $[0,0,2,1,1,0]$ from $[0,0,0,0,0,0]$ by
 the following sequence of opera\-tions:
\begin{gather*}
[0,0,0,0,0,0] \ \mathop\rightarrow^\Psi \ [0,0,0,0,0,1]
\ \mathop\rightarrow^{s_5} \ [0,0,0,0,1,0]
\ \mathop\rightarrow^{s_4} \ [0,0,0,1,0,0]
\ \mathop\rightarrow^{s_3} \ [0,0,1,0,0,0]\\
\quad \mathop\rightarrow^{s_2} \ [0,1,0,0,0,0]
\ \mathop\rightarrow^{s_1} \ [1,0,0,0,0,0]
\ \mathop\rightarrow^\Psi \ [0,0,0,0,0,2]
\ \mathop\rightarrow^{s_5} \ [0,0,0,0,2,0]
\ \mathop\rightarrow^\Psi \ [0,0,0,2,0,1]\\
\quad
  \mathop\rightarrow^{s_5} \ [0,0,0,2,1,0]
\ \mathop\rightarrow^\Psi \ [0,0,2,1,0,1]
\ \mathop\rightarrow^{s_5} \
[0,0,2,1,1,0].
\end{gather*}
Replace $\Psi$ by $\Psi^\alpha$ in the list of the operations, the
associated sequence is
\begin{gather*}
\zeta_{[0,0,0,0,0,0]}= [3,2,0,-1,1,0] \ \mathop\rightarrow^{\Psi^\alpha} \
 \zeta_{[0,0,0,0,0,1]}=[2,0,-1,1,0,\alpha+3] \\
 \quad \mathop\rightarrow^{s_5} \
 \zeta_{[0,0,0,0,1,0]}=[2,0,-1,1,\alpha+3,0] \ \mathop\rightarrow^{s_4} \
 \zeta_{[0,0,0,1,0,0]}=[2,0,-1,\alpha+3,1,0]\\
\quad \mathop\rightarrow^{s_3} \
 \zeta_{[0,0,1,0,0,0]}=[2,0,\alpha+3,-1,1,0] \
 \mathop\rightarrow^{s_2} \
 \zeta_{[0,1,0,0,0,0]}=[2,\alpha+3,0,-1,1,0] \\
\quad \mathop\rightarrow^{s_1} \
 \zeta_{[1,0,0,0,0,0]}=[\alpha+3,2,0,-1,1,0] \
 \mathop\rightarrow^{\Psi^\alpha} \
 \zeta_{[0,0,0,0,0,2]}=[2,0,-1,1,0,2\alpha+3]\\
\quad \mathop\rightarrow^{s_5} \
 \zeta_{[0,0,0,0,2,0]}=[2,0,-1,1,2\alpha+3,0] \
 \mathop\rightarrow^{\Psi^\alpha} \
 \zeta_{[0,0,0,2,0,1]}=[0,-1,1,2\alpha+3,0,\alpha+2] \\
\quad \mathop\rightarrow^{s_5} \
 \zeta_{[0,0,0,2,1,0]}=[0,-1,1,2\alpha+3,\alpha+2,0] \
 \mathop\rightarrow^{\Psi^\alpha} \
 \zeta_{[0,0,2,1,0,1]}=[-1,1,2\alpha+3,\alpha+2,0,\alpha] \\
\quad \mathop\rightarrow^{s_5} \
 \zeta_{[0,0,2,1,1,0]}=[-1,1,2\alpha+3,\alpha+2,\alpha,0].
 \end{gather*}

Now, to obtain the vector-valued Jack polynomial, it suf\/f\/ices to start
from $1\otimes\begin{array}{cccc} 4&3\\6&5&2&1\end{array}$ and act
successively with the af\/f\/ine operator $\Psi$ (when reading
$\Psi^\alpha$) and with $s_i\otimes s_i+\frac1{\zeta[i+1]-\zeta[i]}$
(when reading $s_i$).
\end{Example}

In conclusion,
the computation of vector-valued Jack for a given RST is completely
independent of the computations of the vector-valued Jack indexed by the
other RST with the same shape.

\subsection{Normalization}\label{ss4.3Norm}

The space $V_\lambda$ spanned by the RST $\tau$ of the same shape
$\lambda$ is naturally endowed (up to a multiplicative constant) by
$\S_N$-invariant scalar product $\langle\,,\,\rangle_0$ with respect to which the RST are pairwise
orthogonal. As in \cite{Dunkl}, we set
\[
||\tau||^2=\prod_{1\leq i<j\leq N\atop
\CT_\tau[i]<\CT_\tau[j]-1}{(\CT_\tau[i]-\CT_\tau[j]-1)(\CT_\tau[i]-\CT_\tau[j]+1)\over
(\CT_\tau[i]-\CT_\tau[j])^2}.
\]
As in \cite{Dunkl}, we consider the contravariant form $\langle\ ,\
\rangle$ on the space $M_\lambda$ which is the symmetric
$\S_N$-invariant form extending $\langle\,,\,\rangle_0$ and such that
the Dunkl operator ${\cal D}_i$ is the adjoint to the multiplication by
$x_i$ (see appendix in  \cite{Dunkl} for more details).

The operator $x_i{\cal D}_i$ is self adjoint and the adjoint of $\sigma\in \S_N$ is $\sigma^{-1}$. Since $s_{ij}=s_{ij}^{-1}$ is self adjoint,   ${\cal U}_i$ is self-adjoint for the form $\langle\ ,\ \rangle$ and
the polynomials $J_{v,\tau}$ are pairwise orthogonal.

Let us compute their squared norms $||J_{v,\tau}||^2$ (the bilinear form
is nonsingular for generic $\alpha$ and positive
def\/inite for $\alpha$ in some subset of $\mathbb{R}$ \cite{ES}). The
method is essentially the same as in \cite{Dunkl} and we show that the
result can be read in the Yang--Baxter graph. More precisely, one has

\begin{Proposition}\label{recnorm}\qquad
\begin{enumerate}\itemsep=0pt
\item[$1.$]
$\displaystyle ||J_{(v,\tau)s_i}||^2={(\zeta_{v,\tau}[i+1]-\zeta_{v,\tau}[i]-1)(\zeta_{v,\tau}[i+1]-\zeta_{v,\tau}[i]+1)
\over (\zeta_{v,\tau}[i+1]-\zeta_{v,\tau}[i])^2}||J_{v,\tau}||^2$.
\item[$2.$] $\displaystyle ||J_{(v,\tau)\Psi}||^2=\left(\frac1\alpha
\zeta_{v,\tau}[1]+1\right)||J_{v,\tau}||^2.$
\end{enumerate}
\end{Proposition}

\begin{proof}
1.  Since
\[J_{(v,\tau)s_i}= J_{v,\tau}\left(s_i\otimes
s_i+\frac1{\zeta_{v,\tau}[i+1]-\zeta_{v,\tau}[i]}\right),\]
one obtains
\[ ||J_{v,\tau}(s_i\otimes s_i)||^2=||J_{(v,\tau)s_i}||^2+
\frac1{(\zeta_{v,\tau}[i+1]-\zeta_{v,\tau}[i])^2}||J_{v,\tau}||^2.\]
But $||J_{v,\tau}(s_i\otimes s_i)||^2=||J_{v,\tau}||^2$, which gives the
result.

2.  One has
\begin{gather*}
||J_{(v,\tau)\Psi}||^2 = || J_{v,\tau}(\theta \otimes \theta) x_N||^2
 = \langle J_{v,\tau}(\theta \otimes \theta),
J_{v,\tau}(\theta\otimes \theta) x_N{\cal D}_N \rangle \\
\hphantom{||J_{(v,\tau)\Psi}||^2}{}
 =  \langle J_{v,\tau}(\theta\otimes \theta), J_{v,\tau}x_1 {\cal
D}_1(\theta\otimes \theta)\rangle
 =  \langle J_{v,\tau}, J_{v,\tau}x_1 {\cal D}_1\rangle;
\end{gather*}
 recall that $\theta=s_1s_2\cdots s_{N-1}$.
Since ${\cal U}_1= x_1{\cal D}_1$, one obtains the results.
\end{proof}

\begin{Example}
Let again $\tau=\begin{array}{cc}2\\3&1\end{array}$, we compute the
normalization following the Yang--Baxter graph (see Fig.~\ref{fignorm}).

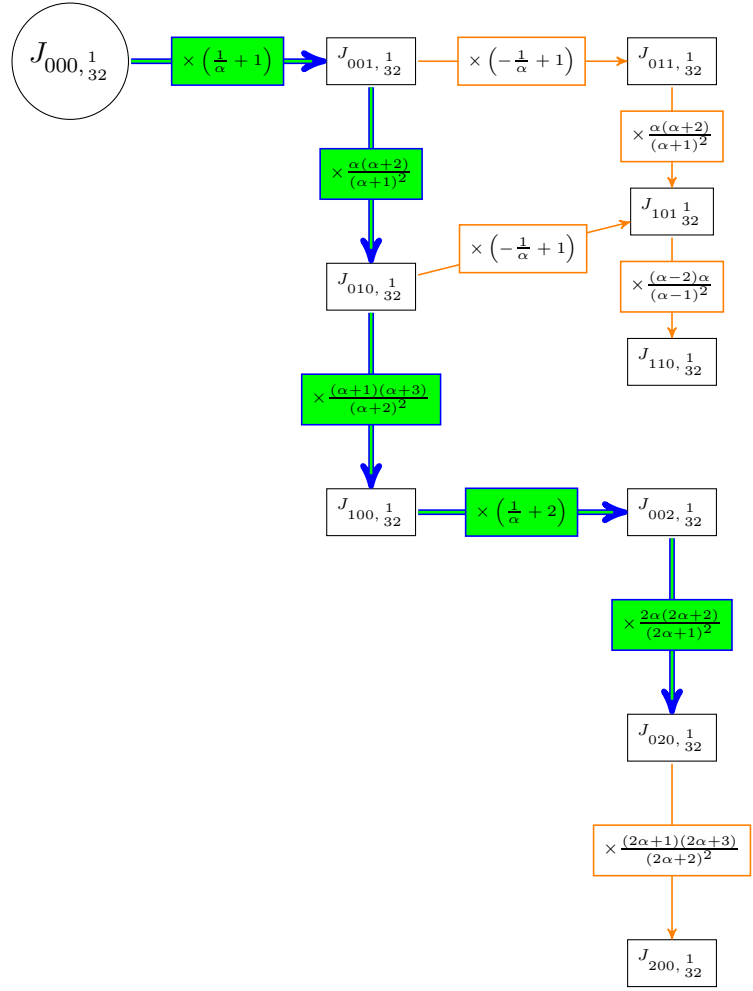
\begin{figure}
\centering
\begin{tikzpicture}
\GraphInit[vstyle=Shade]
    \tikzstyle{VertexStyle}=[shape = rectangle,
draw
]

\Vertex[style={
shape=circle},x=-6, y=0,
 L={$J_{000,{1\ \atop 32}}$}]{x2}

\SetUpEdge[lw = 1.5pt,
color = orange,
 labelcolor = gray!30,
 style={post}
]
 \tikzstyle{TempStyle}=[double = green,
                           double distance = 1pt]

\tikzset{LabelStyle/.style = {draw,
                                     fill = white,
                                     text = black}}

\tikzset{EdgeStyle/.style={post}}
\Vertex[x=-2, y=0, L={\tiny $J_{001,{1\ \atop 32}}$}]{z1}
\Vertex[x=-2, y=-3, L={\tiny $J_{010,{1\atop 32}}$}]{z2}
\Vertex[x=-2, y=-6, L={\tiny $J_{100,{1\ \atop 32}}$}]{z3}
\Edge[style={TempStyle,color=blue,LabelStyle/.style={draw,fill=green}},label={\tiny
$\times {\alpha(\alpha+2)\over(\alpha+1)^2}$}](z1)(z2)
\Edge[style={TempStyle,color=blue,LabelStyle/.style={draw,fill=green}},label={\tiny
$\times {(\alpha+1)(\alpha+3)\over (\alpha+2)^2}$}](z2)(z3)
\Edge[style={TempStyle,color=blue,LabelStyle/.style={draw,fill=green}},label={\tiny
$\times \left(\frac1\alpha+1\right)$}](x2)(z1)
\Vertex[x=2, y=0, L={\tiny $J_{011,{1\ \atop 32}}$}]{xz1}
\Vertex[x=2, y=-2, L={\tiny $ J_{101{1\ \atop 32}}$}]{xz2}
\Vertex[x=2, y=-6, L={\tiny $ J_{002,{1\ \atop 32}}$}]{xz3}
\Edge[label={\tiny $\times\left(-\frac1\alpha+1\right)$}](z1)(xz1)
\Edge[label={\tiny $\times \left(-\frac1\alpha+1\right)$}](z2)(xz2)
\Edge[style={TempStyle,color=blue,LabelStyle/.style={draw,fill=green}},label={\tiny
$\times \left(\frac1\alpha+2\right)$}](z3)(xz3)
\Edge[label={\tiny $\times {\alpha(\alpha+2)\over (\alpha+1)^2}$}](xz1)(xz2)
\Vertex[x=2, y=-4, L={\tiny$ J_{110,{1\ \atop 32}}$}]{xz4}
\Edge[label={\tiny $\times {(\alpha-2)\alpha\over (\alpha-1)^2}$}](xz2)(xz4)
\Vertex[x=2, y=-9, L={\tiny$ J_{020,{1\ \atop 32}}$}]{xz5}
\Vertex[x=2, y=-12, L={\tiny $J_{200,{1\ \atop 32}}$}]{xz6}

\Edge[style={TempStyle,color=blue,LabelStyle/.style={draw,fill=green}},label={\tiny
$\times {2\alpha(2\alpha+2)\over(2\alpha+1)^2}$}](xz3)(xz5)
\Edge[label={\tiny$\times {(2\alpha+1)(2\alpha+3)\over
(2\alpha+2)^2}$}](xz5)(xz6)

\end{tikzpicture}
\caption{\label{fignorm} Computation of $||J_{020,{1\ \atop32}}||^2$
using the graph $G_{1\atop 32}$.}
\end{figure}

 For instance:
\begin{gather*}
||J_{[020],\tau}||^2 = ||\tau||^2\left(1+\frac1\alpha\right)\left(\alpha(\alpha+2)\over(\alpha+1)^2\right)
\left((\alpha+1)(\alpha+3)\over(\alpha+2)^2\right)\left(2+\frac1\alpha\right)\left(2\alpha(2\alpha+2)\over(2\alpha+1)^2\right)\\
\hphantom{||J_{[020],\tau}||^2}{}
 =  {4(\alpha+3)(\alpha+1)\over(2\alpha+1)(\alpha+2)}.
 \end{gather*}
\end{Example}

\section{Symmetrization and antisymmetrization}\label{s5Sym}

In \cite{BF}, Baker and Forrester
investigated the coef\/f\/icients and
the norm of the symmetric Jack polynomials
 by symmetrizing the nonsymmetric Jack
polynomials. 
The symmetrization method was used in~\cite{DG} for the polynomials associated with the complex groups~$G(r,p,N)$.
In this section, we generalize their results and obtain
symmetric and antisymmetric vector-valued Jack polynomials.

\subsection{Non-af\/f\/ine connectivity}\label{ss5.1NonAff}

Let us denote by $H_\lambda$ the graph obtained from $G_\lambda$ by
removing the af\/f\/ine edges, all the falls and the vertex $\varnothing$.
The purpose of this section is to investigate the connected components
of $H_\lambda$.
Recall that $v^+$ is the unique decreasing partition obtained by
permuting the entries of~$v$.

\begin{Definition}
Let $v\in\N^N$ and $\tau\in{\rm Tab}_\lambda$ ($\lambda$ partition). We
def\/ine the f\/illing $T(\tau,v)$ obtained by replacing $i$ by $v^{+}[i]$
in $\tau$ for each~$i$.
\end{Definition}

\begin{Proposition} \label{HT}
Two $4$-tuples $(\tau,\zeta,v,\sigma)$ and $(\tau',\zeta',v',\sigma')$
are in the same connected component of $H_\lambda$ if and only if
$T(\tau,v)=T(\tau',v')$.
\end{Proposition}

\begin{proof}
Remark f\/irst that the steps and correct jumps preserve $T(\tau,v)$.
Indeed steps leave invariant the pairs
$(\tau,v^{+})$ whilst the correct jumps act on the RST by $\tau
s_{\sigma_v[i]}$ where $v[i]=v[i+1]$ (or equivalently by $\tau s_j$
where $v^{+}[j]=v^{+}[j+1]$. Hence, we show that
if $(\tau,\zeta,v,\sigma)$ is connected to $(\tau',\zeta',v',\sigma')$
then $T(\tau,v)=T(\tau',v')$.

Let us prove the converse. Suppose that $T(\tau,v)=T(\tau',v')$.
Since
$(\tau,\zeta,v,\sigma)$ (resp. $(\tau',\zeta'$, $v',\sigma')$) is connected
to $(\tau,\zeta,v^+,{\rm Id})$ (resp. $(\tau',\zeta',{v'}^+,{\it Id})$) by steps, it
suf\/f\/ices to prove the result for $T(\tau,\mu)=T(\tau',\mu)$ when $\mu$
is a (decreasing) partition. Let $\rho\in\S_N$ such that $\tau^\rho$
(the tableau~$\tau$ where the entries have been permuted by $\rho$)
equals $\tau_\lambda$. By construction
$T(\tau,\mu)=T(\tau_\lambda,\mu\rho^{-1})$ and then it suf\/f\/ices to show
the result for
$T(\tau_\lambda,\mu\rho^{-1})=T({\tau'}^\rho,\mu\rho^{-1})$. Again the
connectivity by steps implies that it suf\/f\/ices to prove the result for
$T(\tau_\lambda,\mu)=T(\tau',\mu)$ when $\mu$ is a partition.
 We will show that  if $T(\tau_\lambda,\mu)=T(\tau,\mu)$ then there
exists a series of correct jumps from $(\tau,\zeta',\mu,{\rm Id})$ to
$(\tau_\lambda,\zeta,\mu,{\rm Id})$ when $\mu$ is a partition. We prove
the result by induction on the length of the  shortest
permutation $\omega$ such that $\tau \omega=\tau_\lambda$ for the weak
order. The base point of the induction is straightforward. Now, choose
$i$ such that $\omega[i]>\omega[i+1]$ and $i$ and $i+1$ are
 neither in the same
row nor in the same column in $\tau$
then $\omega=s_i\omega'$ where $\ell(\omega')<\ell(\omega)$.
 Since
$T(\tau,\mu)=T(\tau',\mu)$, this means that $\mu[i]=\mu[i+1]$ and hence,
there is a correct jump from
$(\tau,\zeta',\mu,{\rm Id})$ to $(\tau^{(i,i+1)},\zeta' s_i,\mu,{\rm Id})$. By
the induction hypothesis, this shows the result.
\end{proof}

This shows that the connected components of $H_\lambda$ are indexed by
the $T(\tau,{\mu})$ where $\mu$ is a~partition.

\begin{Definition}
We will denote by $H_T$ the connected component associated to $T$ in
$H_\lambda$. The component $H_T$ will be said to be {\it $1$-compatible} if~$T$ is a column-strict tableau. The compo\-nent~$H_T$ will be said to be {\it
$(-1)$-compatible} if $T$ is a row-strict tableau.
\end{Definition}

\begin{Example}
Let $\mu=[2,1,1,0,0]$ and $\lambda=[3,2]$. There are four connected
components with vertices labeled by permutations of $\mu$ in $H_\lambda$
(see Fig.~\ref{connected}). The possible values of $T(\tau,{\mu})$ are
\[
{12\ \atop 001},\, {02\ \atop 011},\, {01\ \atop 012}\qquad \mbox{and} \qquad {11\
\atop002},
\]squared in red in Fig.~\ref{connected}.
The $1$-compatible components are $H_{12\ \atop 001}$ and $H_{11\
\atop002}$ while there is only one $(-1)$-compatible component $H_{01\
\atop 012}$. The component $H_{02\ \atop 011}$ is neither $1$-compatible
 nor $(-1)$-compatible.

The component $H_{12\ \atop 001}$ contains vertices of $G_{31\ \atop
542}$ and $G_{21\ \atop 543}$ connected by jumps.
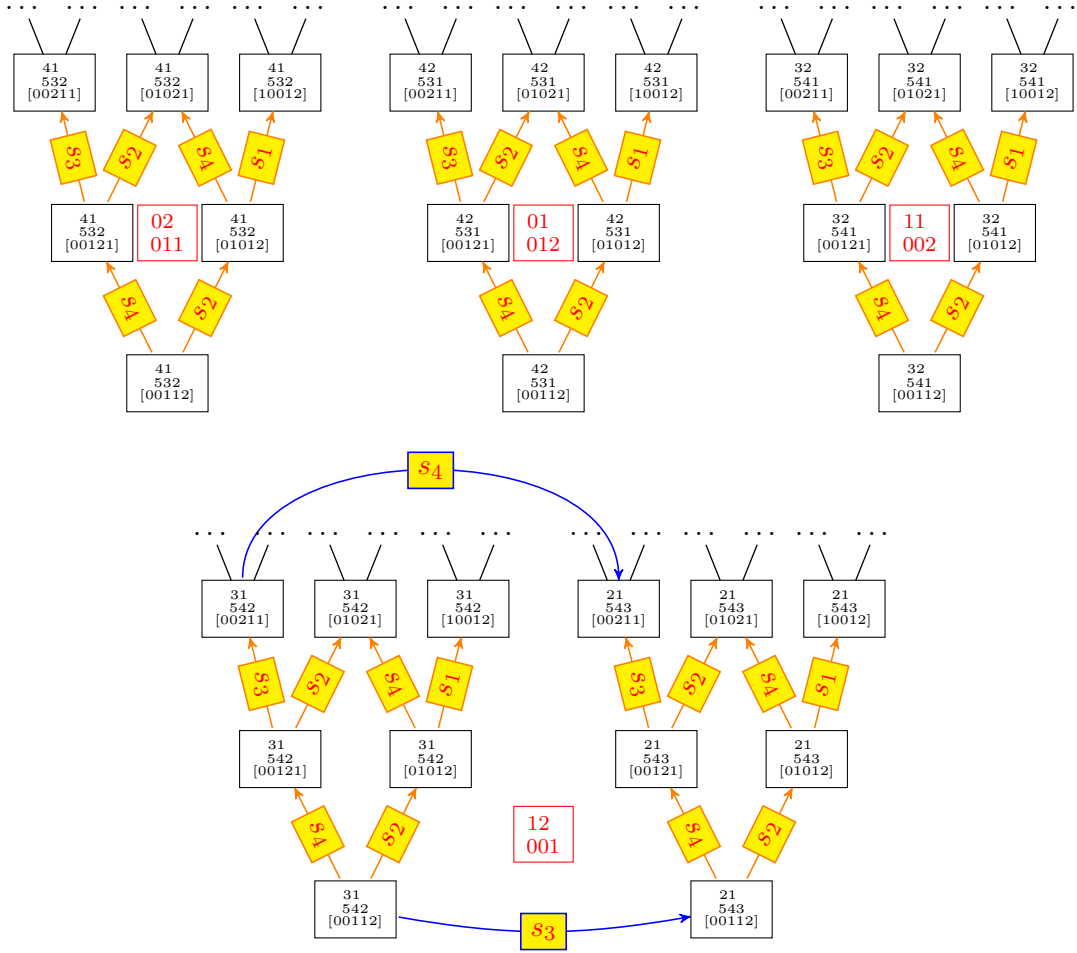
\begin{figure}[t]
\centering
\begin{tikzpicture}%

\GraphInit[vstyle=Shade]
    \tikzstyle{VertexStyle}=[shape = rectangle,
                             draw
]

\SetUpEdge[lw = 1.5pt,
color = orange,
 labelcolor = gray!30,
 labelstyle = {draw,sloped},
 style={post}
]

\tikzset{LabelStyle/.style = {draw,
                                     fill = yellow,
                                     text = red}}

\Vertex[x=0, y=0, L={\tiny ${41\ \atop 532}\atop [00112]$}]{a1}
\Vertex[x=5, y=0, L={\tiny ${42\ \atop 531}\atop [00112]$}]{a2}
\Vertex[x=10, y=0, L={\tiny ${32\ \atop 541}\atop [00112]$}]{a3}

\Vertex[x=-1, y=2, L={\tiny ${41\ \atop 532}\atop [00121]$}]{b11}
\Vertex[x=1, y=2, L={\tiny ${41\ \atop 532}\atop [01012]$}]{b12}

\Vertex[x=0, y=2,
 L={${02\ \atop011}$},style={color=red}
]{x2}

\Vertex[x=-1.5, y=4, L={\tiny ${41\ \atop 532}\atop [00211]$}]{c11}
\Vertex[x=0, y=4, L={\tiny ${41\ \atop 532}\atop [01021]$}]{c12}
\Vertex[x=1.5, y=4, L={\tiny ${41\ \atop 532}\atop [10012]$}]{c13}

\Edge[label={$s_4$}](a1)(b11)
\Edge[label={$s_2$}](a1)(b12)
\Edge[label={$s_3$}](b11)(c11)
\Edge[label={$s_2$}](b11)(c12)
\Edge[label={$s_4$}](b12)(c12)
\Edge[label={$s_1$}](b12)(c13)
\Vertex[x=4, y=2, L={\tiny ${42\ \atop 531}\atop [00121]$}]{b21}
\Vertex[x=6, y=2, L={\tiny ${42\ \atop 531}\atop [01012]$}]{b22}

\Vertex[x=5, y=2,
 L={${01\ \atop012}$},style={color=red}
]{x2}

\Vertex[x=3.5, y=4, L={\tiny ${42\ \atop 531}\atop [00211]$}]{c21}
\Vertex[x=5, y=4, L={\tiny ${42\ \atop 531}\atop [01021]$}]{c22}
\Vertex[x=6.5, y=4, L={\tiny ${42\ \atop 531}\atop [10012]$}]{c23}

\Edge[label={$s_4$}](a2)(b21)
\Edge[label={$s_2$}](a2)(b22)
\Edge[label={$s_3$}](b21)(c21)
\Edge[label={$s_2$}](b21)(c22)
\Edge[label={$s_4$}](b22)(c22)
\Edge[label={$s_1$}](b22)(c23)

\Vertex[x=10, y=2,
 L={${11\ \atop002}$},style={color=red}
]{x2}

\Vertex[x=9, y=2, L={\tiny ${32\ \atop 541}\atop [00121]$}]{b31}
\Vertex[x=11, y=2, L={\tiny ${32\ \atop 541}\atop [01012]$}]{b32}

\Vertex[x=8.5, y=4, L={\tiny ${32\ \atop 541}\atop [00211]$}]{c31}
\Vertex[x=10, y=4, L={\tiny ${32\ \atop 541}\atop [01021]$}]{c32}
\Vertex[x=11.5, y=4, L={\tiny ${32\ \atop 541}\atop [10012]$}]{c33}

\Edge[label={$s_4$}](a3)(b31)
\Edge[label={$s_2$}](a3)(b32)
\Edge[label={$s_3$}](b31)(c31)
\Edge[label={$s_2$}](b31)(c32)
\Edge[label={$s_4$}](b32)(c32)
\Edge[label={$s_1$}](b32)(c33)

\Vertex[x=2.5, y=-7, L={\tiny ${31\ \atop 542}\atop [00112]$}]{a4}
\Vertex[x=7.5, y=-7, L={\tiny ${21\ \atop 543}\atop [00112]$}]{a5}

\Vertex[x=1.5, y=-5, L={\tiny ${31\ \atop 542}\atop [00121]$}]{b41}
\Vertex[x=3.5, y=-5, L={\tiny ${31\ \atop 542}\atop [01012]$}]{b42}

\Vertex[x=1, y=-3, L={\tiny ${31\ \atop 542}\atop [00211]$}]{c41}
\Vertex[x=2.5, y=-3, L={\tiny ${31\ \atop 542}\atop [01021]$}]{c42}
\Vertex[x=4, y=-3, L={\tiny ${31\ \atop 542}\atop [10012]$}]{c43}

\Vertex[x=5, y=-6,
 L={${12\ \atop001}$},style={color=red}
]{x2}

\Vertex[x=6.5, y=-5, L={\tiny ${21\ \atop 543}\atop [00121]$}]{b51}
\Vertex[x=8.5, y=-5, L={\tiny ${21\ \atop 543}\atop [01012]$}]{b52}

\Vertex[x=6, y=-3, L={\tiny ${21\ \atop 543}\atop [00211]$}]{c51}
\Vertex[x=7.5, y=-3, L={\tiny ${21\ \atop 543}\atop [01021]$}]{c52}
\Vertex[x=9, y=-3, L={\tiny ${21\ \atop 543}\atop [10012]$}]{c53}

\Edge[label={$s_4$}](a4)(b41)
\Edge[label={$s_2$}](a4)(b42)
\Edge[label={$s_3$}](b41)(c41)
\Edge[label={$s_2$}](b41)(c42)
\Edge[label={$s_4$}](b42)(c42)
\Edge[label={$s_1$}](b42)(c43)

\Edge[label={$s_4$}](a5)(b51)
\Edge[label={$s_2$}](a5)(b52)
\Edge[label={$s_3$}](b51)(c51)
\Edge[label={$s_2$}](b51)(c52)
\Edge[label={$s_4$}](b52)(c52)
\Edge[label={$s_1$}](b52)(c53)

\Edge[label={$s_3$},style={post,out=-10,in=190},color=blue](a4)(a5)
\Edge[label={$s_4$},style={post,out=90,in=90},color=blue](c41)(c51)

    \tikzstyle{VertexStyle}=[shape = rectangle
]
\Vertex[x=-0.4, y=5,L={$\dots$}]{aa13}
\Vertex[x=0.4, y=5,L={$\dots$}]{aa14}
\Vertex[x=-1.9, y=5,L={$\dots$}]{aa11}
\Vertex[x=-1.1, y=5,L={$\dots$}]{aa12}
\Vertex[x=1.1, y=5,L={$\dots$}]{aa15}
\Vertex[x=1.9, y=5,L={$\dots$}]{aa16}

\Vertex[x=4.6, y=5,L={$\dots$}]{aa23}
\Vertex[x=5.4, y=5,L={$\dots$}]{aa24}
\Vertex[x=3.1, y=5,L={$\dots$}]{aa21}
\Vertex[x=3.9, y=5,L={$\dots$}]{aa22}
\Vertex[x=6.1, y=5,L={$\dots$}]{aa25}
\Vertex[x=6.9, y=5,L={$\dots$}]{aa26}

\Vertex[x=9.6, y=5,L={$\dots$}]{aa33}
\Vertex[x=10.4, y=5,L={$\dots$}]{aa34}
\Vertex[x=8.1, y=5,L={$\dots$}]{aa31}
\Vertex[x=8.9, y=5,L={$\dots$}]{aa32}
\Vertex[x=11.1, y=5,L={$\dots$}]{aa35}
\Vertex[x=11.9, y=5,L={$\dots$}]{aa36}

\Vertex[x=2.1, y=-2,L={$\dots$}]{aa43}
\Vertex[x=2.9, y=-2,L={$\dots$}]{aa44}
\Vertex[x=0.6, y=-2,L={$\dots$}]{aa41}
\Vertex[x=1.4, y=-2,L={$\dots$}]{aa42}
\Vertex[x=3.6, y=-2,L={$\dots$}]{aa45}
\Vertex[x=4.4, y=-2,L={$\dots$}]{aa46}

\Vertex[x=7.1, y=-2,L={$\dots$}]{aa53}
\Vertex[x=7.9, y=-2,L={$\dots$}]{aa54}
\Vertex[x=5.6, y=-2,L={$\dots$}]{aa51}
\Vertex[x=6.4, y=-2,L={$\dots$}]{aa52}
\Vertex[x=8.6, y=-2,L={$\dots$}]{aa55}
\Vertex[x=9.4, y=-2,L={$\dots$}]{aa56}
\SetUpEdge[lw = 0.5pt,
color = black,
 labelstyle = {draw}
]
\Edge(c12)(aa13)
\Edge(c12)(aa14)
\Edge(c11)(aa11)
\Edge(c11)(aa12)
\Edge(c13)(aa15)
\Edge(c13)(aa16)

\Edge(c22)(aa23)
\Edge(c22)(aa24)
\Edge(c21)(aa21)
\Edge(c21)(aa22)
\Edge(c23)(aa25)
\Edge(c23)(aa26)

\Edge(c32)(aa33)
\Edge(c32)(aa34)
\Edge(c31)(aa31)
\Edge(c31)(aa32)
\Edge(c33)(aa35)
\Edge(c33)(aa36)

\Edge(c42)(aa43)
\Edge(c42)(aa44)
\Edge(c41)(aa41)
\Edge(c41)(aa42)
\Edge(c43)(aa45)
\Edge(c43)(aa46)

\Edge(c52)(aa53)
\Edge(c52)(aa54)
\Edge(c51)(aa51)
\Edge(c51)(aa52)
\Edge(c53)(aa55)
\Edge(c53)(aa56)
\end{tikzpicture}
\caption{\label{connected} Some connected components of $H_{32}$.}
\end{figure}
\end{Example}

We we use the following result in the sequel, its proof is easy and left
to the reader.
\begin{Proposition}\label{tabinHT}
Let $(\tau,\zeta,v,\sigma)$ be a vertex of $H_T$ such that
$(\tau,\zeta,v,\sigma)s_i=\varnothing$. One has
\begin{enumerate}\itemsep=0pt
\item[$1.$] If $H_T$ is $1$-compatible then $\sigma[i]$ and
$\sigma[i+1](=\sigma[i]+1)$ are in the same  row in $\tau$.
\item[$2.$] If $H_T$ is $(-1)$-compatible then $\sigma[i]$ and
$\sigma[i+1](=\sigma[i]+1)$ are in the same column in $\tau$.
\end{enumerate}
\end{Proposition}

\def\std#1{{\rm std}(#1)}

The following def\/inition is used to f\/ind a RST corresponding to a
f\/illing of a shape.
\begin{Definition}
Let $T$ be a f\/illing of shape $\lambda$, the {\it standardization} $\std
T$ of $T$ is the reverse standard tableau with shape $\lambda$ obtained
by the following process:
\begin{enumerate}\itemsep=0pt
 \item Denote by $|T|_i$ the number of occurrences of $i$ in $T$
 \item Read the tableau $T$ from the left to the right and the bottom to
the top and replace successively each occurrence of $i$ by the numbers
$N-|T|_0-\dots-|T|_{i-1}$, $N-|T|_0-\dots-|T|_{i-1}-1$, \dots,
$N-|T|_0-\dots-|T|_{i}$.
\end{enumerate}
Alternatively, one has
\begin{gather*}
\std T\ \left[ i,j\right]   :=\#\left\{ \left(
k,l\right) :T\left[ k,l\right] >T\left[ i,j\right] \right\} +\#\left\{
\left( k,l\right) :
l>j,T\left[ k,l\right] =T\left[ i,j\right] \right\}
\\
\phantom{\std T\ \left[ i,j\right]   :=}{}  +\#\left\{ \left( k,j\right) :k\geq i,T\left[ k,j\right] =T\left[
i,j\right] \right\} .
\end{gather*}
We will denote by $\lambda_T$ the unique partition obtained by sorting
in the decreasing order all the entries of $T$.
\end{Definition}

\begin{Example}
Pictorially, reading $01\ \atop 002$ one obtains
\[
 \begin{array}{ccc|cc}
 0&0&2&0&1\\\hline
 0&0&.&0&.\\
 .&.&.&.&1\\
 .&.&2&.&.
 \end{array}
\]
Renumbering in increasing order from the bottom to the top and the right
to the left, one reads
\[
 \begin{array}{ccc|cc}
 0&0&2&0&1\\\hline
 5&4&.&3&.\\
 .&.&.&.&2\\
 .&.&1&.&.
 \end{array}
\]
Hence,
we have
 ${\rm std}\left(01\ \atop 002\right)={32\ \atop 541}$ and $\lambda_{01\
\atop 002}=[21000]$.
\end{Example}

Note that each $H_T$ has a unique sink (that is a vertex with no
outward edge) and this vertex is labeled by $(\std T,\zeta_T,
\lambda_T,{\it Id})$ for a certain vector $\zeta_T$ and a unique root.

\begin{Example}
Consider the tableau $T={01\atop00}$. Its standardization is $\std
T={21\atop 43}$ and the graph~$H_T$ is:
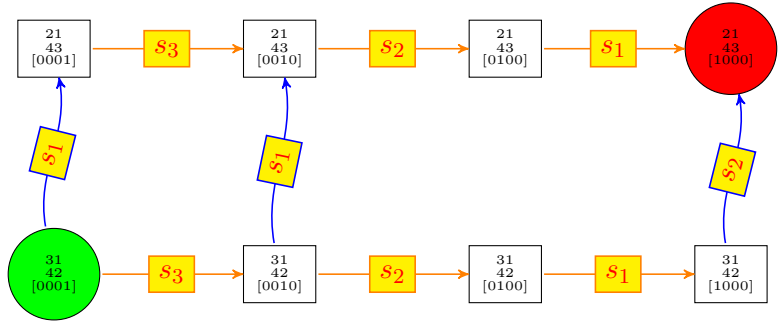
\begin{center}\begin{tikzpicture}%
\GraphInit[vstyle=Shade]
    \tikzstyle{VertexStyle}=[shape = rectangle,
                             draw
]

\SetUpEdge[lw = 1.5pt,
color = orange,
 labelcolor = gray!30,
 labelstyle = {draw,sloped},
 style={post}
]

\tikzset{LabelStyle/.style = {draw,
                                     fill = yellow,
                                     text = red}}

\Vertex[x=0, y=0, L={\tiny ${31\atop 42}\atop
[0001]$},style={shape=circle,fill = green}]{a1}
\Vertex[x=3, y=0, L={\tiny ${31\atop 42}\atop [0010]$}]{a2}
\Vertex[x=6, y=0, L={\tiny ${31\atop 42}\atop [0100]$}]{a3}
\Vertex[x=9, y=0, L={\tiny ${31\atop 42}\atop [1000]$}]{a4}

\Vertex[x=0, y=3, L={\tiny ${21\atop 43}\atop [0001]$}]{b1}
\Vertex[x=3, y=3, L={\tiny ${21\atop 43}\atop [0010]$}]{b2}
\Vertex[x=6, y=3, L={\tiny ${21\atop 43}\atop [0100]$}]{b3}
\Vertex[x=9, y=3, L={\tiny ${21\atop 43}\atop
[1000]$},style={shape=circle,fill = red}]{b4}

\Edge[label={$s_1$},style={post,in=-80,out=100},color=blue](a1)(b1)
\Edge[label={$s_1$},style={post,in=-80,out=100},color=blue](a2)(b2)
\Edge[label={$s_2$},style={post,in=-80,out=100},color=blue](a4)(b4)
\Edge[label={$s_3$}](a1)(a2)
\Edge[label={$s_2$}](a2)(a3)
\Edge[label={$s_1$}](a3)(a4)
\Edge[label={$s_3$}](b1)(b2)
\Edge[label={$s_2$}](b2)(b3)
\Edge[label={$s_1$}](b3)(b4)

\end{tikzpicture}
\end{center}
The sink is denoted by a red disk and the root by a green disk.
\end{Example}

\subsection{Symmetric and antisymmetric Jack polynomials}\label{ss5.2SymJack}
For convenience, let us def\/ine:
\[ (v,\tau)s_i=(v',\tau') \qquad \mbox{if} \quad
(\tau,\zeta,v,\sigma)s_i=(\tau',\zeta',v',\sigma')
\]
and
\[ (v,\tau)s_i=\varnothing \qquad \mbox{if}  \quad (\tau,\zeta,v,\sigma)s_i=\varnothing.\]
Denote also, $J_\varnothing:=0$.

Let $(\tau,\zeta,v,\sigma)$ be a vertex of $H_T$, set
$b_{v,\tau}[i]=\frac1{\zeta_{v,\tau}[i+1]-\zeta_{v,\tau}[i]}$ and
$c_{v,\tau}[i]=\frac{\zeta_{v,\tau}[i]-\zeta_{v,\tau}[i+1]}{\zeta_{v,\tau}[i]-\zeta_{v,\tau}[i+1]+1}$.

Note that
\begin{equation}\label{cbtau1}
1+c_{v,\tau}[i]b_{v,\tau}[i]=c_{v,\tau}[i]
\end{equation}
and
\begin{equation}\label{cbtau2}
c_{v,\tau}[i](1-b_{v,\tau}[i]^2)-b_{v,\tau}[i]=1.
\end{equation}
Let $H_T$ be a $1$-compatible component of $G_\lambda$. For each vertex
$(\tau,\zeta,v,\sigma)$ of $H_T$, we def\/ine the coef\/f\/icient ${\cal
E}_{v,\tau}$ by the following induction:
\begin{enumerate}\itemsep=0pt
\item ${\cal E}_{v,\tau}=1$ if there is no arrows of the form
\begin{center}\begin{tikzpicture}
\GraphInit[vstyle=Shade]
    \tikzstyle{VertexStyle}=[shape = rectangle,
draw
]
\SetUpEdge[lw = 1.5pt,
color = orange,
 labelcolor = gray!30,
 style={post},
labelstyle={sloped}
]
\tikzset{LabelStyle/.style = {draw,
                                     fill = white,
                                     text = black}}
\tikzset{EdgeStyle/.style={post}}
\Vertex[x=0, y=-1,
 L={\tiny$(\tau',\zeta',v',\sigma')$}]{x}
\Vertex[x=3, y=-1,
 L={\tiny $(\tau,\zeta,v,\sigma)$}]{y}
\Edge[label={\tiny$s_i$}](x)(y)
\end{tikzpicture}\end{center}

in $H_T$.
\item ${\cal E}_{v,\tau}={\zeta'[i+1]-\zeta'[i]\over
\zeta'[i+1]-\zeta'[i]-1}{\cal E}_{v',\tau'}=
{\zeta[i+1]-\zeta[i]\over \zeta[i+1]-\zeta[i]+1}{\cal
E}_{v',\tau'}=c_{v',\tau'}{\cal E}_{v',\tau'}$ if there is an arrow

\begin{center}\begin{tikzpicture}
\GraphInit[vstyle=Shade]
    \tikzstyle{VertexStyle}=[shape = rectangle,
draw
]
\SetUpEdge[lw = 1.5pt,
color = orange,
 labelcolor = gray!30,
 style={post},
labelstyle={sloped}
]
\tikzset{LabelStyle/.style = {draw,
                                     fill = white,
                                     text = black}}
\tikzset{EdgeStyle/.style={post}}
\Vertex[x=0, y=-1,
 L={\tiny$(\tau',\zeta',v',\sigma')$}]{x}
\Vertex[x=3, y=-1,
 L={\tiny $(\tau,\zeta,v,\sigma)$}]{y}
\Edge[label={\tiny$s_i$}](x)(y)
\end{tikzpicture}\end{center}
in $H_T$.
\end{enumerate}
The symmetric group acts on the spectral vectors $\zeta$ by permuting
their components. Hence the value of ${\cal E}_{v,\tau}$ does not depend
on the path used for its computation and the ${\cal E}_{v,\tau}$ are
well def\/ined. Indeed, it suf\/f\/ices to check that the def\/inition is
compatible with the commutations $s_is_j=s_js_i$ with $|i-j|>1$ and the
braid relations $s_is_{i+1}s_i=s_{i+1}s_is_{i+1}$.

Let us f\/irst prove the compatibility with the commutation relations.
Suppose
\[
(\tau_0,\zeta_0,v_0,\sigma_0)s_is_j=(\tau_1,\zeta_1,v_1,\sigma_1)s_j=(\tau_2,\zeta_2,v_2,\sigma_2),\]
with $|i-j|>1$ and
\[
(\tau_0,\zeta_0,v_0,\sigma_0)s_js_i=(\tau'_1,\zeta'_1,v'_1,\sigma'_1)s_i=(\tau'_2,\zeta'_2,v'_2,\sigma'_2).\]
Note that $\tau'_2=\tau_2$, $\zeta'_2=\zeta_2$, $v'_2=v_2$ and
$\sigma'_2=\sigma_2$.
 But, since the symmetric group acts on $\zeta$ by permuting its
components, one has
\begin{gather*}
 \zeta_2[j+1]=\zeta'_1[j+1],\qquad \zeta_2[j]=\zeta'_1[j],\qquad
\zeta_1[j+1]=\zeta'_2[j+1]\qquad \mbox{and} \qquad \zeta_1[j]=\zeta'_2[j].
\end{gather*}
Hence,
\begin{gather*}
{\zeta_2[j+1]-\zeta_2[j]\over
\zeta_2[j+1]-\zeta_2[j]+1}.{\zeta_1[i+1]-\zeta_1[i]\over
\zeta_1[i+1]-\zeta_1[i]+1} = {\zeta_1[i+1]-\zeta_1[i]\over
\zeta_1[i+1]-\zeta_1[i]+1}.{\zeta_2[j+1]-\zeta_2[j]\over
\zeta_2[j+1]-\zeta_2[j]+1}\\
\hphantom{{\zeta_2[j+1]-\zeta_2[j]\over
\zeta_2[j+1]-\zeta_2[j]+1}.{\zeta_1[i+1]-\zeta_1[i]\over
\zeta_1[i+1]-\zeta_1[i]+1}}{}
= {\zeta'_1[i+1]-\zeta'_1[i]\over
\zeta'_1[i+1]-\zeta'_1[i]+1}.{\zeta'_2[j+1]-\zeta'_2[j]\over
\zeta'_2[j+1]-\zeta'_2[j]+1},
\end{gather*}
and the def\/inition of ${\cal E}_{v,\tau}$ is compatible with the
commutations.

Now, let us show that the def\/inition is compatible with the braid
relations and set
\[
(\tau_0,\zeta_0,v_0,\sigma_0)s_is_{i+1}s_i=(\tau_1,\zeta_1,v_1,\sigma_1)s_
{i+1}s_i=(\tau_2,\zeta_2,v_2,\sigma_2)s_i=(\tau_3,\zeta_3,v_3,\sigma_3),
\]
and
\[
(\tau_0,\zeta_0,v_0,\sigma_0)s_{i+1}s_{i}s_{i+1}=(\tau'_1,\zeta'_1,v'_1,\sigma'_1)s_
{i}s_{i+1}=(\tau'_2,\zeta'_2,v'_2,\sigma'_2)s_{i+1}=(\tau'_3,\zeta'_3,v'_3,\sigma'_3).
\]
Note that $\tau'_3=\tau_3$, $\zeta'_3=\zeta_3$, $v'_3=v_3$ and
$\sigma'_3=\sigma_3$.
Since the symmetric group acts on $\zeta$ by permuting its components,
one has
\begin{gather*}
\zeta_3[i+1]=\zeta'_1[i+1],\qquad \zeta_3[i]=\zeta'_1[i+1],\\
\zeta_2[i+2]=\zeta'_2[i+1], \qquad \zeta_2[i+1]=\zeta'_2[i],\\
\zeta_1[i+1]=\zeta'_2[i+2]\qquad \mbox{and} \qquad \zeta_1[i]=\zeta'_3[i+1].
\end{gather*}
Hence,
\begin{gather*}
{\zeta_3[i+1]-\zeta_3[i]\over
\zeta_3[i+1]-\zeta_3[i]+1}.{\zeta_2[i+2]-\zeta_2[i+1]\over
\zeta_2[i+2]-\zeta_2[i+1]+1}.{\zeta_1[i+1]-\zeta_1[i]\over
\zeta_1[i+1]-\zeta_1[i]+1} \\
\qquad{} =
{\zeta_1[i+1]-\zeta_1[i]\over
\zeta_1[i+1]-\zeta_1[i]+1}.{\zeta_2[i+2]-\zeta_2[i+1]\over
\zeta_2[i+2]-\zeta_2[i+1]+1}.{\zeta_3[i+1]-\zeta_3[i]\over
\zeta_3[i+1]-\zeta_3[i]+1}\\
\qquad{} =
{\zeta'_3[i+2]-\zeta'_3[i+1]\over
\zeta'_3[i+2]-\zeta'_3[i+1]+1}.{\zeta'_2[i+1]-\zeta'_2[i]\over
\zeta'_2[i+1]-\zeta'_2[i]+1}.{\zeta'_1[i+2]-\zeta'_1[i+1]\over
\zeta'_1[i+2]-\zeta'_1[i+1]+1},
\end{gather*}
and the def\/inition is compatible with the braid relations.

Def\/ine the symmetrization operator
\[
{\cal S}:=\sum_{\omega\in\S_N} \omega\otimes\omega.
\]

We will say that a polynomial is symmetric if it is invariant under the
action of $s_i\otimes s_i$ for each $i<N$.

\begin{Theorem}\label{Sym}\qquad
\begin{enumerate}\itemsep=0pt
\item[$1.$] Let $H_T$ be a connected component of $G_\lambda$. For each vertex
$(\tau,\zeta,v,\sigma)$ of $H_T$, the polyno\-mial~$J_{v,\tau}{\cal S}$
equals $J_{\lambda_T,\std T}{\cal S}$ up to a multiplicative constant.
\item[$2.$] One has $J_{\lambda_T,\std T}{\cal S}\neq 0$ if and only if $H_T$
is $1$-compatible.
\item[$3.$] More precisely, when $H_T$ is $1$-compatible, the polynomial
\[J_{T}=\sum_{(\tau,\zeta,v,\sigma) \ {\rm vertex \ of} \ H_T}{\cal
E}_{v,\tau}J_{v,\tau}\]
is symmetric.
\end{enumerate}
\end{Theorem}

\begin{proof}
1.  Let us prove the f\/irst assertion by induction on the length of a
path from $(\tau,\zeta,v,\sigma)$ to $(\std T,\zeta_T,\lambda_T,\sigma)$
in $H_T$. Let $(\tau',\zeta',v',\sigma')$ such that
\begin{center}
\begin{tikzpicture}
\GraphInit[vstyle=Shade]
    \tikzstyle{VertexStyle}=[shape = rectangle,
draw
]
\SetUpEdge[lw = 1.5pt,
color = orange,
 labelcolor = gray!30,
 style={post},
labelstyle={sloped}
]
\tikzset{LabelStyle/.style = {draw,
                                     fill = white,
                                     text = black}}
\tikzset{EdgeStyle/.style={post}}
\Vertex[x=0, y=-1,
 L={\tiny$(\tau,\zeta,v,\sigma)$}]{x}
\Vertex[x=3, y=-1,
 L={\tiny $(\tau',\zeta',v',\sigma')$}]{y}
\Edge[label={\tiny$s_i$}](x)(y)
\end{tikzpicture}\end{center}
is not a jump in $H_T$ (hence, $-1<b_{v,\tau}[i]<1$).
It follows that
\begin{gather*}
 J_{v,\tau}\S = \frac1{1-b_{v,\tau}[i]^2}\left(J_{v',\tau'}(s_i\otimes
s_i)+b_{v,\tau}[i]J_{v',\tau'}\right)\S
 = \frac1{1-b_{v,\tau}[i]^2}\left(1+b_{v,\tau}[i]\right)J_{v',\tau'}\S.
\end{gather*}
By induction $J_{v',\tau'}\S$ is proportional to $J_{\lambda_T,\std T}$,
which ends the proof.

2. If $H_T$ is not $1$-compatible, then there exists $s_i$ such that
$J_{\lambda_T,\std T}(s_i\otimes s_i)=-J_{\lambda_T,\std T}$. Hence,
since $\S=(s_i\otimes s_i)\S$, one obtains $J_{\lambda_T,\std T}\S=0$.

3. Let us prove that, when $H_T$ is $1$-compatible,
  $J_T(s_i\otimes s_i)=J_T$ for any $i$.
Fix $i$ and decompose $J_T:=J^++J_0+J^-$ where
\[J^{+}={\sum_{\tau,v}}^+ {\cal E}_{\tau,v}J_{v,\tau}, \]
where $\sum^+$ means that the sum is over the pairs $(\tau,v)$ such that
there exists an arrow
\begin{center}
\begin{tikzpicture}
\GraphInit[vstyle=Shade]
    \tikzstyle{VertexStyle}=[shape = rectangle,
draw
]
\SetUpEdge[lw = 1.5pt,
color = orange,
 labelcolor = gray!30,
 style={post},
labelstyle={sloped}
]
\tikzset{LabelStyle/.style = {draw,
                                     fill = white,
                                     text = black}}
\tikzset{EdgeStyle/.style={post}}
\Vertex[x=0, y=0,
 L={\tiny$(\tau',\zeta',v',\sigma')$}]{x}
\Vertex[x=3, y=0,
 L={\tiny$(\tau,\zeta,v,\sigma)$}]{y}
\Edge[label={\tiny$s_i$}](x)(y)
\end{tikzpicture}
\end{center}
in $H_T$,
\[
J^{-}={\sum_{\tau,v}}^- {\cal E}_{\tau,v}J_{v,\tau}
\]
where $\sum^-$ means that the sum is over the pairs $(\tau,v)$ such that
there exists an arrow
\begin{center}
\begin{tikzpicture}
\GraphInit[vstyle=Shade]
    \tikzstyle{VertexStyle}=[shape = rectangle,
draw
]
\SetUpEdge[lw = 1.5pt,
color = orange,
 labelcolor = gray!30,
 style={post},
labelstyle={sloped}
]
\tikzset{LabelStyle/.style = {draw,
                                     fill = white,
                                     text = black}}
\tikzset{EdgeStyle/.style={post}}
\Vertex[x=0, y=0,
 L={\tiny$(\tau,\zeta,v,\sigma)$}]{x}
\Vertex[x=3, y=0,
 L={\tiny$(\tau',\zeta',v',\sigma')$}]{y}
\Edge[label={\tiny$s_i$}](x)(y)
\end{tikzpicture}
\end{center}
in $H_T$ and
\[
J_0={\sum}^0 {\cal E}_{\tau,v}J_{v,\tau},
\]
where $\sum^0$ means that the sum is over the pairs $(\tau,v)$ such that
there exists an arrow
\begin{center}
\begin{tikzpicture}
\GraphInit[vstyle=Shade]
    \tikzstyle{VertexStyle}=[shape = rectangle,
draw
]
\SetUpEdge[lw = 1.5pt,
color = orange,
 labelcolor = gray!30,
 style={post},
labelstyle={sloped}
]
\tikzset{LabelStyle/.style = {draw,
                                     fill = white,
                                     text = black}}
\tikzset{EdgeStyle/.style={post}}
\Vertex[x=0, y=0,
 L={\tiny$(\tau,\zeta,v,\sigma)$}]{x}
\Vertex[x=2, y=0,
 L={\tiny$\varnothing$}]{y}
\Edge[label={\tiny$s_i$}](x)(y)
\end{tikzpicture}
\end{center}
in $G_T$ (equivalently there is no arrow from $(\tau,\zeta,v,\sigma)$
labeled by $s_i$ in $H_T$).
Suppose that
\begin{center}
\begin{tikzpicture}
\GraphInit[vstyle=Shade]
    \tikzstyle{VertexStyle}=[shape = rectangle,
draw
]
\SetUpEdge[lw = 1.5pt,
color = orange,
 labelcolor = gray!30,
 style={post},
labelstyle={sloped}
]
\tikzset{LabelStyle/.style = {draw,
                                     fill = white,
                                     text = black}}
\tikzset{EdgeStyle/.style={post}}
\Vertex[x=0, y=0,
 L={\tiny$(\tau,\zeta,v,\sigma)$}]{x}
\Vertex[x=2, y=0,
 L={\tiny$\varnothing$}]{y}
\Edge[label={\tiny$s_i$}](x)(y)
\end{tikzpicture}
\end{center}
 is a fall in $G_T$, then
\[
J_{v,\tau}(s_i\otimes
s_i)=J_{(v,\tau)s_i}-b_{v,\tau}[i]J_{v,\tau}=-b_{v,\tau}[i]J_{v,\tau}.
\]
Since, $H_T$ is $1$-compatible Proposition~\ref{tabinHT} implies that
$i$ and $i+1$ are in the same  row. Hence,
$b_{v,\tau}[i]=-1$ and $J_{v,\tau}(s_i\otimes s_i)=J_{v,\tau}$. It
follows that $J_0(s_i\otimes s_i)=J_0$.

Now, let
\begin{center}
\begin{tikzpicture}
\GraphInit[vstyle=Shade]
    \tikzstyle{VertexStyle}=[shape = rectangle,
draw
]
\SetUpEdge[lw = 1.5pt,
color = orange,
 labelcolor = gray!30,
 style={post},
labelstyle={sloped}
]
\tikzset{LabelStyle/.style = {draw,
                                     fill = white,
                                     text = black}}
\tikzset{EdgeStyle/.style={post}}
\Vertex[x=0, y=0,
 L={\tiny$(\tau,\zeta,v,\sigma)$}]{x}
\Vertex[x=3, y=0,
 L={\tiny$(\tau',\zeta',v',\sigma')$}]{y}
\Edge[label={\tiny$s_i$}](x)(y)
\end{tikzpicture}
\end{center}
be an arrow in $H_T$, then
\begin{gather*}
J_{v,\tau}(s_i\otimes s_i)=J_{v',\tau'}-b_{v,\tau}[i]J_{v,\tau},
\\
J_{v',\tau'}(s_i\otimes
s_i)=b_{v,\tau}[i]J_{v',\tau'}+(1-b_{v,\tau}[i]^2)J_{v,\tau}
\qquad \mbox{and}
\qquad
{\goth E}_{v',\tau'}=c_{v,\tau}{\goth E}_{v,\tau}.
\end{gather*}
Hence, equalities (\ref{cbtau1}) and (\ref{cbtau2}) imply
\begin{gather*}
({\cal E}_{v,\tau}J_{v,\tau}+{\cal E}_{v',\tau'}J_{v',\tau'})(s_i\otimes
s_i) = {\cal
E}_{v,\tau}(J_{v,\tau}+c_{v,\tau}[i]J_{v',\tau'})(s_i\otimes s_i) \\
\qquad{}
 =  {\cal
E}_{v,\tau}\left(((c_{v,\tau}[i](1-b_{v,\tau}[i]^2)-b_{v,\tau}[i])J_{v,\tau} +(1+c_{v,\tau}[i]b_{v,\tau}[i])J_{v',\tau'}\right)\\
\qquad{} = {\cal E}_{v,\tau}(J_{v,\tau}+c_{v,\tau}[i]J_{v',\tau'}) = ({\cal
E}_{v,\tau}J_{v,\tau}+{\cal E}_{v',\tau'}J_{v',\tau'}).
\end{gather*}
This proves that $(J^++J^-)(s_i\otimes s_i)=J^++J^-$. Hence,
$J_T(s_i\otimes s_i)=J_T$ for each $i$ and $J_T$ is symmetric.
\end{proof}

\begin{Example}
Consider the graph $H_{11\atop 00}$
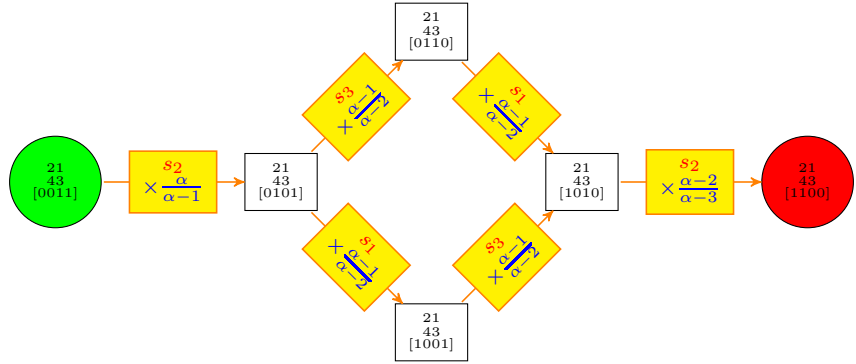
\begin{center}\begin{tikzpicture}%

\GraphInit[vstyle=Shade]
    \tikzstyle{VertexStyle}=[shape = rectangle,
                             draw
]

\SetUpEdge[lw = 1.5pt,
color = orange,
 labelcolor = gray!30,
 labelstyle = {draw,sloped},
 style={post}
]

\tikzset{LabelStyle/.style = {draw,
                                     fill = yellow,
                                     text = red}}

\Vertex[x=0, y=0, L={\tiny ${21\atop 43}\atop
[0011]$},style={shape=circle,fill = green}]{a1}
\Vertex[x=3, y=0, L={\tiny ${21\atop 43}\atop [0101]$}]{a2}
\Vertex[x=5, y=2, L={\tiny ${21\atop 43}\atop [0110]$}]{a3}
\Vertex[x=5, y=-2, L={\tiny ${21\atop 43}\atop [1001]$}]{a4}
\Vertex[x=7, y=0, L={\tiny ${21\atop 43}\atop [1010]$}]{a5}
\Vertex[x=10, y=0, L={\tiny ${21\atop 43}\atop
[1100]$},style={shape=circle,fill = red}]{a6}

\Edge[label={$s_2\atop \blue\times {\alpha\over \alpha-1}$}](a1)(a2)
\Edge[label={$s_1\atop \blue\times {\alpha-1\over \alpha-2}$}](a2)(a4)
\Edge[label={$s_3\atop \blue\times {\alpha-1\over \alpha-2}$}](a2)(a3)
\Edge[label={$s_1\atop \blue\times {\alpha-1\over \alpha-2}$}](a3)(a5)
\Edge[label={$s_3\atop \blue\times {\alpha-1\over \alpha-2}$}](a4)(a5)
\Edge[label={$s_2\atop \blue\times {\alpha-2\over \alpha-3}$}](a5)(a6)

\end{tikzpicture}
\end{center}
The polynomial
 \begin{gather*}
J_{11\atop00} = J_{0011,{21\atop 43}}+{\alpha\over
\alpha-1}J_{0101,{21\atop 43}}+{\alpha\over \alpha-2}J_{0110,{21\atop
43}}+{\alpha\over \alpha-2}J_{1001,{21\atop 43}}
+{\alpha(\alpha-1)\over (\alpha-2)^2}J_{1010,{21\atop
43}}\\
\phantom{J_{11\atop00} =}{} +{\alpha(\alpha-1)\over (\alpha-2)(\alpha-3)}J_{1100,{21\atop 43}}
\end{gather*}
is symmetric.
\end{Example}

\def\sink{{\rm sink}}
\def\root{{\rm root}}

Let $H_T$ be a connected component, denote by ${\rm root}(T)$ the only
vertex of $H_T$ without inward edge and by ${\sink}(T)=({\rm
std}(T),\zeta_T,\lambda_T,{\rm Id})$ the only vertex of $H_T$ without outward
edge.
Denote by $\# H_T$ the number of vertices of $H_T$.
The following proposition allows to compare the polynomial $J_T$ to the
symmetrization of $J_{\root(T)}$.
\begin{Proposition}\label{TtoRootSym}
 One has
\[
 J_T={\# H_T\over N!}{\goth E}_{\sink(T)}J_{\root(T)}{\cal S}.
\]
\end{Proposition}
\begin{proof}
It suf\/f\/ices to compare the coef\/f\/icient of $J_{\sink(T)}$ in $J_T$ and in
$J_{\root(T)}.{\cal S}$.
The coef\/f\/icient of $J_{\sink(T)}$ in $J_T$ equals ${\goth E}_{\sink(T)}$
while the coef\/f\/icient of
$J_{\sink(T)}$ in $J_{\root(T)}{\cal S}$ equals $N!\over \# H$. Indeed
$N!\over \# H$ is the order of the stabilizer of $\lambda_T$. The
leading monomial of $J_{\sink(T)}$ does not appear in any other $J_{v,\tau}$
so its coef\/f\/icient in the symmetrization of $J_{\root(T)}$ equals the
order of the stabilizer.
\end{proof}

Let $H_T$ be a $(-1)$-compatible component of $G_\lambda$. For each
vertex $(\tau,\zeta,v,\sigma)$ of $H_T$, we def\/ine the coef\/f\/icient
${\cal F}_{v,\tau}$ by the following induction:
\begin{enumerate}\itemsep=0pt
\item ${\cal F}_{v,\tau}=1$ if there is no arrow of the
form\begin{center}\begin{tikzpicture}
\GraphInit[vstyle=Shade]
    \tikzstyle{VertexStyle}=[shape = rectangle,
draw
]
\SetUpEdge[lw = 1.5pt,
color = orange,
 labelcolor = gray!30,
 style={post},
labelstyle={sloped}
]
\tikzset{LabelStyle/.style = {draw,
                                     fill = white,
                                     text = black}}
\tikzset{EdgeStyle/.style={post}}
\Vertex[x=0, y=0,
 L={\tiny$(\tau',\zeta',v',\sigma')$}]{x}
\Vertex[x=3, y=0,
 L={\tiny$(\tau,\zeta,v,\sigma)$}]{y}
\Edge[label={\tiny$s_i$}](x)(y)
\end{tikzpicture}\end{center} in $H_T$.
\item ${\cal F}_{v,\tau}=-{\zeta[i]-\zeta[i+1]\over
\zeta[i]-\zeta[i+1]+1}{\cal F}_{v',\tau'}
=-{\zeta'[i+1]-\zeta[i]\over \zeta'[i+1]-\zeta[i]+1}{\cal F}_{v',\tau'}$
if there is an arrow \begin{center}\begin{tikzpicture}
\GraphInit[vstyle=Shade]
    \tikzstyle{VertexStyle}=[shape = rectangle,
draw
]
\SetUpEdge[lw = 1.5pt,
color = orange,
 labelcolor = gray!30,
 style={post},
labelstyle={sloped}
]
\tikzset{LabelStyle/.style = {draw,
                                     fill = white,
                                     text = black}}
\tikzset{EdgeStyle/.style={post}}
\Vertex[x=0, y=0,
 L={\tiny$(\tau',\zeta',v',\sigma')$}]{x}
\Vertex[x=3, y=0,
 L={\tiny$(\tau,\zeta,v,\sigma)$}]{y}
\Edge[label={\tiny$s_i$}](x)(y)
\end{tikzpicture} \end{center}in $H_T$.
\end{enumerate}

Again the ${\cal F}_{v,\tau}$ are well def\/ined since the symmetric group
acts on the spectral vectors by permuting their components.
Def\/ine also the antisymmetrization operator
\[
{\cal A}:=\sum_{\omega\in\S_N}(-1)^{\ell(\omega)} (\omega\otimes\omega).
\]

We will say that a polynomial is antisymmetric if it vanishes under
the action of $1-s_i\otimes s_i$ for each $i<N$.

\begin{Theorem}\label{Antisym}\qquad
\begin{enumerate}\itemsep=0pt
\item[$1.$] Let $H_T$ be a connected component of $G_\lambda$. For each vertex
$(\tau,\zeta,v,\sigma)$ of $H_T$, the polyno\-mial~$J_{v,\tau}{\cal A}$
equals $J_{\lambda_T,\std T}{\cal A}$ up to a multiplicative constant.
\item[$2.$] One has $J_{\lambda_T,\std T}{\cal A}\neq 0$ if and only if $H_T$
is $(-1)$-compatible.
\item[$3.$] More precisely, when $H_T$ is $(-1)$-compatible, the polynomial
\[J'_{T}=\sum_{(\tau,\zeta,v,\sigma)\ {\rm  vertex \ of} \ H_T}{\cal
F}_{v,\tau}J_{v,\tau}\]
is antisymmetric.
\end{enumerate}
\end{Theorem}

\begin{Example}\rm
Consider the graph $H_{01\atop 01}$
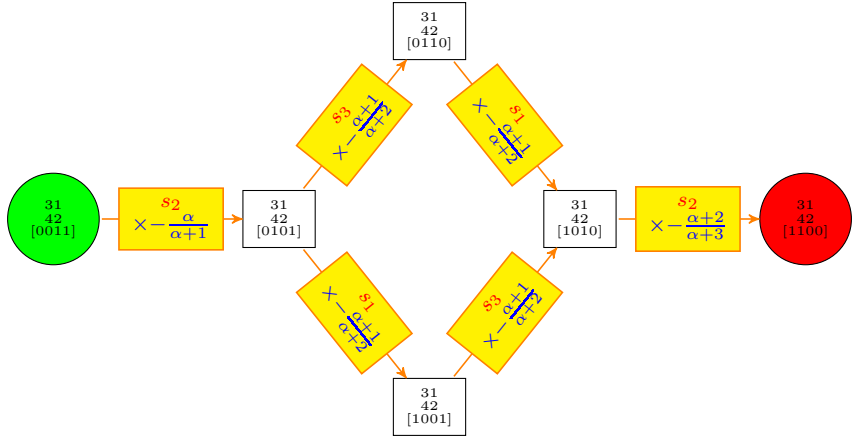
\begin{center}\begin{tikzpicture}%

\GraphInit[vstyle=Shade]
    \tikzstyle{VertexStyle}=[shape = rectangle,
                             draw
]

\SetUpEdge[lw = 1.5pt,
color = orange,
 labelcolor = gray!30,
 labelstyle = {draw,sloped},
 style={post}
]

\tikzset{LabelStyle/.style = {draw,
                                     fill = yellow,
                                     text = red}}

\Vertex[x=0, y=0, L={\tiny ${31\atop 42}\atop
[0011]$},style={shape=circle,fill = green}]{a1}
\Vertex[x=3, y=0, L={\tiny ${31\atop 42}\atop [0101]$}]{a2}
\Vertex[x=5, y=2.5, L={\tiny ${31\atop 42}\atop [0110]$}]{a3}
\Vertex[x=5, y=-2.5, L={\tiny ${31\atop 42}\atop [1001]$}]{a4}
\Vertex[x=7, y=0, L={\tiny ${31\atop 42}\atop [1010]$}]{a5}
\Vertex[x=10, y=0, L={\tiny ${31\atop 42}\atop
[1100]$},style={shape=circle,fill = red}]{a6}

\Edge[label={$s_2\atop \blue\times -{\alpha\over \alpha+1}$}](a1)(a2)
\Edge[label={$s_1\atop \blue\times -{\alpha+1\over \alpha+2}$}](a2)(a4)
\Edge[label={$s_3\atop \blue\times -{\alpha+1\over \alpha+2}$}](a2)(a3)
\Edge[label={$s_1\atop \blue\times -{\alpha+1\over \alpha+2}$}](a3)(a5)
\Edge[label={$s_3\atop \blue\times -{\alpha+1\over \alpha+2}$}](a4)(a5)
\Edge[label={$s_2\atop \blue\times -{\alpha+2\over \alpha+3}$}](a5)(a6)

\end{tikzpicture}
\end{center}
The polynomial
 \begin{gather*}
J'_{01\atop01} = J_{0011,{31\atop 42}}-{\alpha\over
\alpha+1}J_{0101,{31\atop 42}}+{\alpha\over \alpha+2}J_{0110,{31\atop
42}}+{\alpha\over \alpha+2}J_{1001,{31\atop 42}}
-{\alpha(\alpha+1)\over (\alpha+2)^2}J_{1010,{31\atop 42}}\\
\phantom{J'_{01\atop01} =}{} +
{\alpha(\alpha+1)\over (\alpha+2)(\alpha+3)}J_{1100,{31\atop 42}}
\end{gather*}
is antisymmetric.
\end{Example}

And, as in the symmetric case, one has:
\begin{Proposition}\label{TtoRootAntiSym}
 One has
\[
 J_T={\# H_T\over N!}{\goth F}_{\sink(T)}J_{\root(T)}.{\cal A}.
\]
\end{Proposition}

\subsection{Normalization}\label{s5.3NormSym}
As a consequence of Proposition~\ref{recnorm}, one deduces the following
result using Theorems~\ref{Sym} and~\ref{Antisym}.
\begin{Corollary}\label{CorNorm}
Let $H_T$ be a connected component and $(\tau,\zeta,v,\sigma)$ be a
vertex of $H_T$. Denote by $\ell_{\tau,v}^T$ the length of a path from
${\rm root}(T)$ to $(\tau,\zeta,v,\sigma)$. One has,
\[
||J_{v,\tau}||^2=(-1)^{\ell_{\tau,v}^T}{\goth E}_{v,\tau}^{-1}{\goth
F}_{v,\tau}^{-1}||J_{\rm root}(T)||^2.
\]
\end{Corollary}

From Theorems~\ref{Sym} and~\ref{Antisym}, vector-valued symmetric and
antisymmetric Jack polynomials are also pairwise orthogonal.
\begin{Proposition}\qquad\sloppy
\begin{enumerate}\itemsep=0pt
\item[$1.$] Let $H_{T_1}$ and $H_{T_2}$ be two $1$-compatible connected
components. If $T_1\neq T_2$ then $\langle J_{T_1},J_{T_2}\rangle=0$.
\item[$2.$] Let $H_{T_1}$ and $H_{T_2}$ be two $(-1)$-compatible connected
components. If $T_1\neq T_2$ then $\langle J'_{T_1},J'_{T_2}\rangle=0$.
\end{enumerate}
\end{Proposition}

\begin{proof} It suf\/f\/ices to remark that from Theorem~\ref{Sym} (resp.\
Theorem~\ref{Antisym}) each $J_{T}$ (resp. $J'_T$) is a linear
combination of $J_{v,\tau}$ for $(\tau,\zeta,v,\sigma)$ vertex in the
connected component $H_T$. \end{proof}
In the special cases when $H_T$ is $\pm1$-compatible, the value of
$||J_T||^2$ admits a remarkable equality.
\begin{Proposition}\label{NormT2Root}
 One has:
\begin{enumerate}\itemsep=0pt
 \item[$1.$] If $H_T$ is a $1$-compatible connected component then
\[ ||J_T||^2=\# H_T{\goth E}_{{\rm sink}(T)}||J_{{\rm root}(T)}||^2.\]
 \item[$2.$] If $H_T$ is a $(-1)$-compatible connected component then
\[ ||J'_T||^2=\# H_T{\goth F}_{{\rm sink}(T)}||J_{{\rm root}(T)}||^2.\]
\end{enumerate}
\end{Proposition}

\begin{proof}  The two cases being very
similar, let us only prove the symmetric case.
From Proposition~\ref{TtoRootSym}, one has:
 \begin{gather*}
||J_T||^2 = {\# H_T\over N!}{\goth E}_{\sink(T)}\langle
J_T,J_{\root(T)}.{\cal S}\rangle
 = {\# H_T\over N!}{\goth E}_{\sink(T)}\sum_{\sigma\in\S_N}\langle
J_T,J_{\root(T)}(\sigma\otimes\sigma)\rangle\\
\phantom{||J_T||^2}{}
 = \# H_T{\goth E}_{{\rm sink}(T)}||J_{{\rm root}(T)}||^2.\tag*{\qed}
\end{gather*}
\renewcommand{\qed}{}
\end{proof}

From Corollary~\ref{CorNorm} and Theorem~\ref{NormT2Root}, one obtains
the surprising equalities:
\begin{Corollary}
If $H_T$ is $1$-compatible, one has:
\begin{equation}\label{sumlemmanorm}\sum_{(\tau,\zeta,v,\sigma)\ {\rm vertex \ of} \ H_T} (-1)^{\ell_{v,\tau}^T}{{\goth E}_{v,\tau}\over {\goth
F}_{v,\tau}}=\# H_T{\goth E}_{{\rm sink}(T)}.
\end{equation}
If $H_T$ is $(-1)$-compatible, one has:
\begin{equation*}
\sum_{(\tau,\zeta,v,\sigma)\ {\rm vertex \
of} \ H_T} (-1)^{\ell_{v,\tau}^T}{{\goth F}_{v,\tau}\over {\goth
E}_{v,\tau}}=\# H_T{\goth F}_{{\rm sink}(T)}.
\end{equation*}
\end{Corollary}

\begin{Example}
Consider the graph $H_{11\atop 00}$, the sum (\ref{sumlemmanorm}) gives
\[
1+{\alpha+1\over\alpha-1}\left(1+{\alpha\over\alpha-2}\left(2+
{\alpha\over\alpha-2}\left(1+{\alpha-1\over\alpha-3}\right)\right)\right)=
6\,{\frac {\alpha \left( \alpha-1 \right) }{ \left( \alpha-2 \right)
\left( \alpha-3
 \right) }}
\]
as expected.

\end{Example}
\subsection{Symmetric and antisymmetric polynomials with minimal
degree\label{s5.4minim}}

Since the irreducible characters of $\mathfrak{S}_{N}$ are real it follows
that the tensor product of an irreducible module with itself contains the
trivial representation exactly once. The tensor product of the module
corresponding to a partition $\lambda$ with the module for $^{t}\lambda$
(the transpose) contains the sign representation exactly once. We
demonstrate these facts explicitly. Using the concepts from Section
\ref{ss4.1YB}
let
\[
\zeta_{1}=\sum_{\tau\in\mathrm{Tab}_{\lambda}}a\left( \tau\right)
(\tau\otimes\tau)\in V_{\lambda}\otimes V_{\lambda}%
\]
be symmetric with (rational) coef\/f\/icients $a\left( \tau\right) $ to be
determined. We impose the conditions $\zeta_{1}(s_{i}\otimes
s_{i})=\zeta_{1}$
for $i=1,\ldots,N-1$. Fix some $i$ and   split the sum as suggested by equation~(\ref{Murphy}) %
\begin{gather*}
\zeta_{1}   =\sum_{b_{\tau
}\left[ i\right] =\pm1}a\left( \tau\right) (\tau\otimes\tau)
  +\sum_{0<b_{\tau}\left[ i\right] \leq\frac{1}{2}}\left( a\left( \tau\right) (\tau\otimes\tau)+a\big( \tau^{\left(
i,i+1\right) }\big)\big( \tau^{\left( i,i+1\right) }\otimes\tau^{\left(
i,i+1\right) }\big)\right) .
\end{gather*}
In the f\/irst sum $\left( \tau\otimes\tau\right) \left( s_{i}\otimes
s_{i}\right) =b_{\tau}\left[ i\right]
^{2}(\tau\otimes\tau)=\tau\otimes\tau$.
For the second sum, note that $\tau^{\left( i,i+1\right) }s_{i}=\big(
1-b_{\tau}\left[ i\right] ^{2}\big) \tau-b_{\tau}\left[ i\right]
\tau^{\left( i,i+1\right) }$. Simple computations show that
\begin{gather*}
  \left( a\left( \tau\right) (\tau\otimes\tau)+a\big( \tau^{\left(
i,i+1\right) }\big) \big(\tau^{\left( i,i+1\right) }\otimes\tau^{\left(
i,i+1\right) }\big)\right) \left( s_{i}\otimes s_{i}\right) \\
\qquad{}  =a\left( \tau\right) (\tau\otimes\tau)+a\big( \tau^{\left( i,i+1\right)
}\big) \big(\tau^{\left( i,i+1\right) }\otimes\tau^{\left(
i,i+1\right) }\big)
\end{gather*}
exactly when $a\left( \tau\right) =\big( 1-b_{\tau}\left[ i\right]
^{2}\big) a\left( \tau^{\left( i,i+1\right) }\right) $. The unique (up
to a constant multiple) $\mathfrak{S}_{N}$-invariant norm on $V_{\tau}$
satisf\/ies $\big\Vert \tau^{\left( i,i+1\right) }\big\Vert ^{2}=\big(
1-b_{\tau}\left[ i\right] ^{2}\big) \left\Vert \tau\right\Vert ^{2}$ (see
Section~\ref{ss4.3Norm}); thus $a\left( \tau\right) =c/\left\Vert
\tau\right\Vert ^{2}$ for
some constant~$c$.

Consider the module $V_{^{t}\lambda}$. The transpose map takes each RST
$\tau$
with shape $\lambda$ to the RST $^{t}\tau$ of shape $^{t}\lambda$. Thus
$b_{^{t}\tau}\left[ i\right] =-b_{\tau}\left[ i\right] $ for $1\leq i\leq
N$. Suppose $0<b_{\tau}\left[ i\right] \leq\frac{1}{2}$ for some $\tau$ and
$i$, then $-\frac{1}{2}\leq b_{^{t}\tau}\left[ i\right] <0$ and the
following transformation rules apply:
\begin{gather*}
^{t}\tau s_{i}   =b_{^{t}\tau}\left[ i\right] \,{}^{t}\tau+\big(
1-b_{^{t}\tau}\left[ i\right] ^{2}\big) \big( ^{t}\tau^{\left(
i,i+1\right) }\big) ,\\
^{t}\tau^{\left( i,i+1\right) } s_{i}   = {}^{t}\tau-b_{^{t}\tau}\left[
i\right] \,{}^{t}\tau^{\left( i,i+1\right) }.
\end{gather*}
Let
\[
\zeta_{\det}=\sum_{\tau\in\mathrm{Tab}_{\lambda}}a\left( \tau\right) \left(
^{t}\tau\right) \otimes\tau\in V_{^{t}\lambda}\otimes V_{\lambda}%
\]
be antisymmetric with (rational) coef\/f\/icients $a\left( \tau\right) $ to be
determined. We impose the conditions $\zeta_{\det} (s_{i}\otimes s_{i})%
=-\zeta_{\det}$ for $i=1,\ldots,N-1$. Fix some $i$ and write%
\begin{gather*}
\zeta_{\det}   =\sum_{b_{\tau}\left[ i\right] =\pm1} a\left( \tau\right) \left( ^{t}\tau\right)
\otimes\tau
  +\sum_{0<b_{\tau}\left[
i\right] \leq\frac{1}{2}} \left(a\left( \tau\right) \left( ^{t}\tau\right) \otimes
\tau+a\big( \tau^{\left( i,i+1\right) }\big) \big( ^{t}\tau^{\left(
i,i+1\right) }\big) \otimes\tau^{\left( i,i+1\right)}\right).
\end{gather*}
In the f\/irst sum $\left( \tau\otimes\left( ^{t}\tau\right) \right)
 \left( s_{i}\otimes s_{i}\right) =b_{\tau}\left[ i\right] b_{^{t}\tau
}\left[ i\right] \left( ^{t}\tau\right) \otimes\tau=-\left( ^{t}%
\tau\right) \otimes\tau$. We f\/ind that
\begin{gather*}
  \left( a\left( \tau\right) \left( ^{t}\tau\right) \otimes\tau+a\big(
\tau^{\left( i,i+1\right) }\big) \big( ^{t}\tau^{\left( i,i+1\right)
}\big) \otimes\tau^{\left( i,i+1\right) }\right) \left( s_{i}\otimes
s_{i}\right) \\
\qquad =-\left( a\left( \tau\right) \left( ^{t}\tau\right) \otimes
\tau+a\big( \tau^{\left( i,i+1\right) }\big) \big( ^{t}\tau^{\left(
i,i+1\right) }\big) \otimes\tau^{\left( i,i+1\right) }\right)
\end{gather*}
exactly when $a\left( \tau\right) =-a\left( \tau^{\left( i,i+1\right)
}\right) $.

 Thus $a\left( \tau\right) =c\left( -1\right) ^{\mathrm{inv}%
\left( \tau\right) }$ (recall $\mathrm{inv}\left( \tau\right) =\#\left\{
\left( i,j\right) :1\leq i<j\leq N,\mathrm{rw}\left( i,\tau\right)
>\mathrm{rw}\left( j,\tau\right) \right\} $, and $0<b_{\tau}\left[
i\right] \leq\frac{1}{2}$ implies $\mathrm{inv}\left( \tau^{\left(
i,i+1\right) }\right) =\mathrm{inv}\left( \tau\right) +1$).

We can now write down the symmetric and antisymmetric Jack polynomials of
lowest degree, by replacing the f\/irst factors in $\zeta_{1}$ and
$\zeta_{\det
}$ by the corresponding polynomials $P_{\tau}\left( x\right) $ and
$P_{^{t}\tau}\left( x\right) $ (as constructed in Section~\ref{s3VectorValued}). Let
$l=\ell\left( \lambda\right) =~^{t}\lambda\left[ 1\right] $.

In the symmetric case let $v=\big[ \left( l-1\right) ^{\lambda\left[
l\right] },\left( l-2\right) ^{\lambda\left[ l-1\right] },\ldots
,1^{\lambda\left[ 2\right] },0^{\lambda\left[ 1\right] }\big] $ (using
exponents to indicate the multiplicity of an entry) The corresponding
tableau is
\[
T_{1}:= \begin{array}{ccccccc}
l-1&\dots&l-1&&&&(\lambda[l]\times)\\
\vdots&&&\ddots&&&\vdots\\
1&\dots&\dots&\dots&1&&(\lambda[2]\times)\\
0&\dots&\dots&\dots&\dots&0&(\lambda[1]\times)
\end{array}
\]
and $\mathrm{std}\left( T_{1}\right) $ contains the numbers $N,N-1,\ldots
,2,1$ entered row-by-row.

\begin{Example}
 If $\lambda=\left[
4,3,2\right] $ then $v=\left[ 221110000\right] $. The corresponding
tableau is
\[
T_1=\begin{array}{cccc}
2&2\\
1&1&1\\
0&0&0&0
\end{array}
\qquad\mbox{and}\qquad \std{T_1}=\begin{array}{cccc}
2&1\\
5&4&3\\
9&8&7&6
\end{array}.
\]
\end{Example}

In
\[
\zeta_{1}\left( x\right) =\sum_{\tau\in\mathrm{Tab}_{\lambda}}\frac
{c}{\left\Vert \tau\right\Vert ^{2}}P_{\tau}\left( x\right) \otimes\tau,
\]
the monomial $x^{v}$ occurs only when $\tau=\mathrm{std}\left( T_{1}\right)
$, with coef\/f\/icient $c/\left\Vert \mathrm{std}\left( T_{1}\right)
\right\Vert ^{2}$. This polynomial is a multiple of $J_{T_{1}}$ (see Theorem
\ref{Sym}).

For the antisymmetric case let%
\[
T_{\det}:=
 \begin{array}{cccccccc}
0 & 1 &\dots &\lambda_l-1\\
\vdots&\vdots & &\vdots &\ddots\\
0 &1 &\dots &\lambda_l-1&\dots&\lambda_2-1\\
0 &1 &\dots &\lambda_l-1&\dots&\lambda_2-1&\dots&\lambda_1-1
\end{array}
\]
Thus $\mathrm{std}\left( T_{\det}\right) =\tau_{\lambda}$ and $v=\big[
\left( \lambda\left[ 1\right] -1\right) ^{^{t}\lambda\left[
\lambda\left[ 1\right] \right] },\left( \lambda\left[ 1\right]
-2\right) ^{^{t}\lambda\left[ \lambda\left[ 1\right] -1\right] }%
,\ldots,0^{^{t}\lambda\left[ 1\right] }\big] $.

\begin{Example}   If
$\lambda=\left[ 4,3,2\right] $ then $^{t}\lambda=\left[ 3,3,2,1\right] $
and $v=\left[ 322111000\right] $. The corresponding tableau is
\[
T_{\rm det}=\begin{array}{cccc}0&1\\0&1&2\\0&1&2&3 \end{array}
\qquad \mbox{and} \qquad
 \std{T_{\rm det}}=\begin{array}{cccc}7&4\\8&5&2\\9&6&3&1
\end{array}=\tau_{[4,3,2]}
.
\]
\end{Example}

 Let
\[
\zeta_{\det}\left( x\right) =\sum_{\tau\in\mathrm{Tab}_{\lambda}}\left(
-1\right) ^{\mathrm{inv}\left( \tau\right) }P_{^{t}\tau}\left( x\right)
\otimes\tau.
\]
The monomial $x^{v}$ occurs only in the term $\tau=\tau_{\lambda}$ (see
Def\/inition~\ref{Glambda}). This polynomial is a~constant multiple of
$J_{T_{\det}%
}^{\prime}$ (see Theorem~\ref{Antisym}).

 We summarize the results of this section in the following theorem.
\begin{Theorem}\label{Minimal}
The subspace of $M_\lambda$ of the symmetric $($resp.\ antisymmetric$)$
polynomials with minimal degree is spanned by only one generator: the
symmetric $($resp.\ antisymmetric$)$ Jack polynomial $J_{T_1}$ $($resp.\
$J_{T_{\rm det}})$.
\end{Theorem}

As a consequence one observes a remarkable property.

\begin{Corollary}
The Jack polynomial $J_{T_1}$ $($resp.\ $J_{T_{\rm det}})$ is equal to a
polynomial which does not depend on the parameter $\alpha$ multiplied by
the global multiplicative constant ${\goth E}_{\sink(T_1)}$ $($resp.\
${\goth F}_{\sink(T_{\rm det})})$.
\end{Corollary}

\begin{proof} The f\/irst part of the sentence is a consequence of Theorem~\ref{Minimal} since the dimension of the space is~$1$. The values of the
multiplicative constants follow from Theorems~\ref{Sym} and~\ref{Antisym} together with the fact that the coef\/f\/icient of the leading
terms in a Jack polynomials $J_{v,\tau}$ is $1$ (see Theorem
\ref{thJack}).\end{proof}

Note also that $T_1$ (resp.\ $T_{\rm det}$) is not the only tableau for
which the corresponding symmetric (resp. antisymmetric) Jack does not
depend on $\alpha$ (up to a global multiplicative constant).

\begin{Example}
Consider the partition $\lambda=[221]$ together with the vector
$v=[2,1,1,0,0]$. The corresponding symmetric Jack $\frac1{{\goth
E}_{21100,{1\ \atop{32\atop 54}}}}{J_{2\ \atop{11\atop 00}}}$
does not depend on $\alpha$.

There are two symmetric Jack polynomials in degree $5$: $J_{2\
\atop{12\atop 00}}$ and $J_{3\ \atop{11\atop 00}}$. Note that the (non
minimal) polynomial $\frac1{{\goth E}_{22100,{1\ \atop{32\atop
54}}}}J_{2\ \atop{12\atop 00}}$ does not depend on $\alpha$ whilst the
parameter $\alpha$ appears in $\frac1{{\goth E}_{31100,{1\ \atop{32\atop
54}}}}J_{3\ \atop{11\atop 00}}$ even after simplifying the expression.
\end{Example}

\section{Restrictions}\label{s6Rest}

\subsection[Restrictions on Yang-Baxter graphs]{Restrictions on Yang--Baxter graphs}\label{ss6.1RestYang}

\def\rest#1{{\displaystyle\mathop\downarrow_#1}}

Consider the operator $\rest M$ acting on the Yang--Baxter graphs
$G_\lambda$ by producing a new graph $G_\lambda \rest M$ following the
rules below:
\begin{enumerate}\itemsep=0pt
\item Add all the possible edges of the form
 \begin{center}
\begin{tikzpicture}
\GraphInit[vstyle=Shade]
    \tikzstyle{VertexStyle}=[shape = rectangle,
draw
]
\SetUpEdge[lw = 1.5pt,
color = orange,
 labelcolor = gray!30,
 style={post},
labelstyle={sloped}
]
\tikzset{LabelStyle/.style = {draw,
                                     fill = white,
                                     text = black}}
\tikzset{EdgeStyle/.style={post}}
\Vertex[x=0, y=0,
 L={\tiny$(\tau,\zeta, [v[1],\dots,v[M],0,\dots, 0],\sigma)$}]{x}
\Vertex[x=8, y=0,
 L={\tiny$(\tau,\zeta', [v[2],\dots,v[M],v[1]+1,0,\dots, 0],\sigma')$}]{y}
\Edge[label={\tiny$\Psi'$ },style={post,in=182,out=2},color=green](x)(y)
\end{tikzpicture}
\end{center}
 More precisely, the action of $\Psi'$ on the $4$-tuples is given by
\[\Psi'=\Psi s_{N-1}\cdots s_M.\]
\item Suppress the vertices labeled by $(\tau,\zeta,v,\sigma)$ with
$v[i]\neq 0$ for some $i>M$, with the associated inward and outward edges.
\item Relabel the remained vertices $(\tau,\zeta,v,\sigma)\rest M:=
(\tau\rest M,\zeta\rest M,v\rest M,\sigma\rest M)$ with
\begin{enumerate}\itemsep=0pt
\item $\tau\rest M$ is obtained from $\tau$ by removing the nodes
labeled by $M+1,\dots, N$. Note that the shape of $\tau\rest M$ could be
a skew partition.
\item $v\rest M=[v[1],\dots,v[M]]$.
\item $\sigma\rest M=[\sigma[1],\dots,\sigma[M]]$.
\item $\zeta\rest M=[\zeta[1]-\CT_\tau[M],\dots,\zeta[M]-\CT_\tau[M]]$.
\end{enumerate}
\item Relabel by $\Psi$ the edges labeled by $\Psi'$.
\end{enumerate}
\begin{Example}\rm
Consider the partition $\lambda=21$ and $M=2$, the graph $G_{21}$ in
Fig.~\ref{G21} with edges~$\Psi'$ added. We obtain the graph
$G_{21}\rest M$ (Fig.~\ref{rhoG21}) applying the other rules.
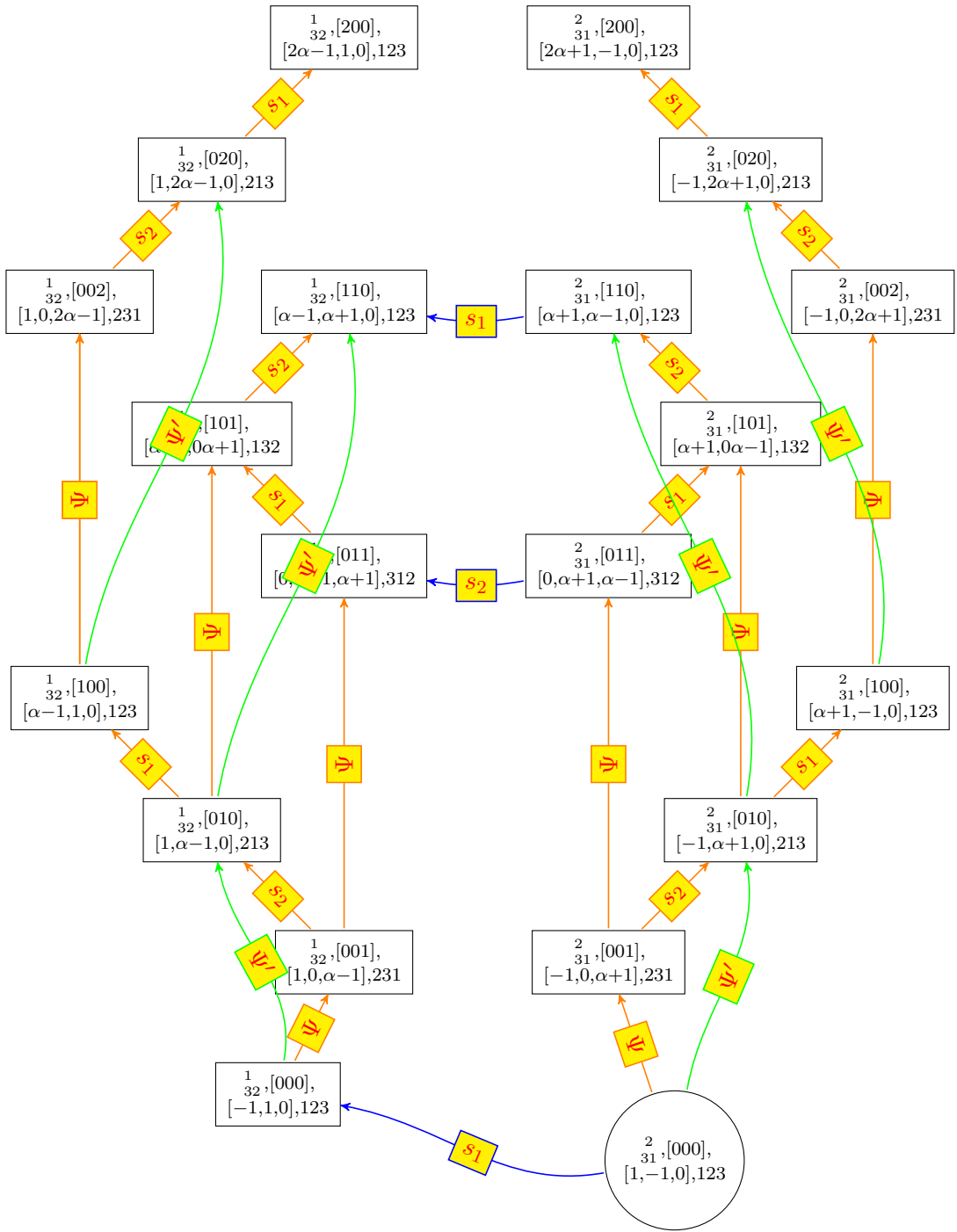
\begin{figure}[t]\centering
\begin{tikzpicture}%
\GraphInit[vstyle=Shade]
    \tikzstyle{VertexStyle}=[shape = rectangle,
                             draw
]

\Vertex[x=3, y=-1,
 L={${2\atop31}, [000],\atop
[1,-1,0],123$},style={shape=circle}
]{x2}
\SetUpEdge[lw = 1.5pt,
color = orange,
 labelcolor = gray!30,
 labelstyle = {draw,sloped},
 style={post}
]

\tikzset{LabelStyle/.style = {draw,
                                     fill = yellow,
                                     text = red}}
\Vertex[x=-3, y=0, L={${1\ \atop32}, [000],\atop [-1,1,0],123$}]{x1}
\Edge[label={$s_1$},style={post,in=-10,out=190},color=blue](x2)(x1)
\Vertex[x=-2, y=2, L={${1\ \atop32}, [001],\atop [1,0,\alpha -1],231$}]{y1}
\Vertex[x=-4, y=4, L={${1\ \atop32}, [010],\atop [1,\alpha -1,0],213$}]{y2}
\Vertex[x=-6, y=6, L={${1\ \atop32}, [100],\atop [\alpha-1,1,0],123$}]{y3}
\Vertex[x=2, y=2, L={${2\ \atop31}, [001],\atop [-1,0,\alpha +1],231$}]{z1}
\Vertex[x=4, y=4, L={${2\ \atop31}, [010],\atop [-1,\alpha +1,0],213$}]{z2}
\Vertex[x=6, y=6, L={${2\ \atop31}, [100],\atop [\alpha+1,-1,0],123$}]{z3}

\Edge[label={$s_2$}](y1)(y2)
\Edge[label={$s_1$}](y2)(y3)
\Edge[label={$s_2$}](z1)(z2)
\Edge[label={$s_1$}](z2)(z3)
\Edge[label={$\Psi$}](x1)(y1)
\Edge[label={$\Psi$}](x2)(z1)

\Vertex[x=-2, y=8, L={${1\ \atop32}, [011],\atop [0,\alpha -1,\alpha
+1],312$}]{xy1}
\Vertex[x=2, y=8, L={${2\ \atop31}, [011],\atop [0,\alpha +1,\alpha
-1],312$}]{xz1}
\Vertex[x=-4, y=10, L={${1\ \atop32}, [101],\atop [\alpha-1,0\alpha
+1],132$}]{xy2}
\Vertex[x=4, y=10, L={${2\ \atop31}, [101],\atop [\alpha+1,0\alpha
-1],132$}]{xz2}
\Vertex[x=-6, y=12, L={${1\ \atop32}, [002],\atop
[1,0,2\alpha-1],231$}]{xy3}
\Vertex[x=6, y=12, L={${2\ \atop31}, [002],\atop
[-1,0,2\alpha+1],231$}]{xz3}
\Edge[label={$\Psi$}](y1)(xy1)
\Edge[label={$\Psi$}](y2)(xy2)
\Edge[label={$\Psi$}](y3)(xy3)
\Edge[label={$\Psi$}](z1)(xz1)
\Edge[label={$\Psi$}](z2)(xz2)
\Edge[label={$\Psi$}](z3)(xz3)

\Edge[label={$s_2$},style={post,in=-10,out=190},color=blue](xz1)(xy1)
\Edge[label={$s_1$}](xy1)(xy2)
\Edge[label={$s_1$}](xz1)(xz2)
\Edge[label={$\Psi$}](y3)(xy3)

\Vertex[x=-2, y=12, L={${1\ \atop32}, [110],\atop [\alpha-1,\alpha
+1,0],123$}]{xy4}
\Vertex[x=2, y=12, L={${2\ \atop31}, [110],\atop [\alpha+1,\alpha
-1,0],123$}]{xz4}
\Edge[label={$s_2$}](xy2)(xy4)
\Edge[label={$s_2$}](xz2)(xz4)
\Edge[label={$s_1$},style={post,in=-10,out=190},color=blue](xz4)(xy4)

\Vertex[x=-4, y=14, L={${1\ \atop32}, [020],\atop
[1,2\alpha-1,0],213$}]{xy5}
\Vertex[x=-2, y=16, L={${1\ \atop32}, [200],\atop
[2\alpha-1,1,0],123$}]{xy6}
\Vertex[x=4, y=14, L={${2\ \atop31}, [020],\atop
[-1,2\alpha+1,0],213$}]{xz5}
\Vertex[x=2, y=16, L={${2\ \atop31}, [200],\atop
[2\alpha+1,-1,0],123$}]{xz6}

\Edge[label={$s_2$}](xz3)(xz5)
\Edge[label={$s_1$}](xz5)(xz6)
\Edge[label={$s_2$}](xy3)(xy5)
\Edge[label={$s_1$}](xy5)(xy6)

\Edge[label={$\Psi'$},style={post,in=-80,out=80},color=green](x1)(y2)
\Edge[label={$\Psi'$},style={post,in=-80,out=80},color=green](x2)(z2)
\Edge[label={$\Psi'$},style={post,in=-80,out=80},color=green](y2)(xy4)
\Edge[label={$\Psi'$},style={post,in=-80,out=80},color=green](y3)(xy5)
\Edge[label={$\Psi'$},style={post,in=-80,out=80},color=green](z2)(xz4)
\Edge[label={$\Psi'$},style={post,in=-80,out=80},color=green](z3)(xz5)

\end{tikzpicture}
\caption{\label{PsiG21} The f\/irst vertices of the graph $G_{21}$ with
edges $\Psi'$ for $M=2$.}\vspace{-1mm}
\end{figure}

\begin{figure}[t]\centering
\begin{tikzpicture}%
\GraphInit[vstyle=Shade]
    \tikzstyle{VertexStyle}=[shape = rectangle,
                             draw
]

\Vertex[x=3, y=-1,
 L={${2\atop\ 1}, [0],\atop
[0,2],12$},style={shape=circle}
]{x2}

\SetUpEdge[lw = 1.5pt,
color = orange,
 labelcolor = gray!30,
 labelstyle = {draw,sloped},
 style={post}
]

\tikzset{LabelStyle/.style = {draw,
                                     fill = yellow,
                                     text = red}}
\Vertex[x=-3, y=0, L={${1\atop\ 2}, [00],\atop [-2,0],12$}]{x1}
\Edge[label={$s_1$},style={post,in=-10,out=190},color=blue](x2)(x1)
\Vertex[x=-4, y=4, L={${1\atop\ 2}, [01],\atop [0,\alpha -2],21$}]{y2}
\Vertex[x=-6, y=6, L={${1\atop\ 2}, [10],\atop [\alpha-2,0],12$}]{y3}

\Vertex[x=4, y=4, L={${2\atop\ 1}, [01],\atop [0,\alpha +2],21$}]{z2}
\Vertex[x=6, y=6, L={${2\atop\ 1}, [10],\atop [\alpha+2,0],12$}]{z3}

\Edge[label={$s_1$}](y2)(y3)

\Edge[label={$s_1$}](z2)(z3)
\Vertex[x=-2, y=12, L={${1\atop\ 2}, [11],\atop [\alpha-2,\alpha
],12$}]{xy4}
\Vertex[x=2, y=12, L={${2\atop\ 1}, [11],\atop [\alpha+2,\alpha ],12$}]{xz4}
\Edge[label={$s_1$},style={post,in=-10,out=190},color=blue](xz4)(xy4)

\Vertex[x=-4, y=14, L={${1\atop\ 2}, [02],\atop [0,2\alpha-2],21$}]{xy5}
\Vertex[x=-2, y=16, L={${1\atop\ 2}, [20],\atop [2\alpha-2,0],12$}]{xy6}
\Vertex[x=4, y=14, L={${2\atop\ 1}, [02],\atop [0,2\alpha+2],21$}]{xz5}
\Vertex[x=2, y=16, L={${2\atop\ 1}, [20],\atop [2\alpha+2,0],12$}]{xz6}

\Edge[label={$s_1$}](xz5)(xz6)
\Edge[label={$s_1$}](xy5)(xy6)

\Edge[label={$\Psi$},style={post,in=-80,out=80},color=green](x1)(y2)
\Edge[label={$\Psi$},style={post,in=-80,out=80},color=green](x2)(z2)
\Edge[label={$\Psi$},style={post,in=-80,out=80},color=green](y2)(xy4)
\Edge[label={$\Psi$},style={post,in=-80,out=80},color=green](y3)(xy5)
\Edge[label={$\Psi$},style={post,in=-80,out=80},color=green](z2)(xz4)
\Edge[label={$\Psi$},style={post,in=-80,out=80},color=green](z3)(xz5)

\end{tikzpicture}
\caption{\label{rhoG21} The f\/irst vertices of the graph $G_{21}\rest 2$.}\vspace{-2mm}
\end{figure}
\end{Example}

\begin{Definition}
A RST $\tau$ has the property $R(M)$ if the removal of the nodes labeled
by $M+1,\dots,N$ in $\tau$ produces a RST whose Ferrers diagram is a
partition.
\end{Definition}

\begin{Example}
The RST
\[\begin{array}{ccc}
 5&2\\7&3&1\\8&6&4
\end{array}\]
has the property $R(3)$ while the RST $2\ \atop 31$ does not have
property $R(2)$.
\end{Example}
Denote by $\overline {G_\tau}$ the subgraph of $G_\lambda$ whose root is
$\tau$. In particular, one has

\begin{Proposition}\label{PropRest}\looseness=-1
Let $\tau$ have the property $R(M)$ and satisfy $\tau\rest M\!=\!\tau_{\lambda\rest
M}\!$ where $\lambda\rest M\!$ denotes the~Fer\-rers diagram of $\tau\rest M$.
The graph $G_{\lambda\rest M}\!$ is identical to the subgraph $\overline
{G_{\tau\rest M}}\!$ of $G_{\lambda}\rest M\!$ whose root is~$\tau\rest M$.
\end{Proposition}

\vspace{-3mm}

\begin{proof}
Obviously, since the Ferrers diagram of $\tau\rest M$ is a partition,
all the spectral vectors $\zeta$ labeling the vertices of $\overline
{G_{\tau\rest M}}$ are obtained by subtracting the same integer (that is
$\CT_\tau[M]$) from the corresponding spectral vector in $G_{\lambda}$.
It follows that the action of the $s_i$ permutes the components of the
spectral vectors in $\overline {G_{\tau\rest M}}$.

Let ${\bf v'}=(\tau',\zeta',[v'[1],\dots,v'[M],0,\dots,0],\sigma')$ be a
vertex of $\overline G_\tau$. Let us prove by induction on the length of
a path from the root to $\bf v'$ that
\begin{enumerate}\itemsep=0pt
\item There is a vertex labeled by ${\bf v'}\rest M:=(\tau'\rest
M,\zeta'\rest M,[v'[1],\dots,v'[M]],\sigma'\rest M)$
in $G_{\lambda\rest M}$.
\item If there is a non af\/f\/ine edge labeled by $s_i$ with $i<M$ from
\[{\bf v''}=(\tau'',\zeta'',[v''[1],\dots,v''[M],0,\dots,0],\sigma'')\] to
${\bf v'}$ in $\overline G_\tau$ then there is the same edge from ${\bf
v''}\rest M$ to ${\bf v'}\rest M$ in $G_{\lambda\rest M}$.
\item If there is an edge from ${\bf
v''}=(\tau'',\zeta'',[v''[1],\dots,v''[M],0,\dots,0],\sigma'')$ to ${\bf
v'}$ in $\overline {G_\tau}$ labeled by $\Psi'$ then there is an edge
labeled by $\Psi$ from ${\bf v''}\rest M$ to ${\bf v'}\rest M$ in
$G_{\lambda\rest M}$.
\end{enumerate}

First, observe that if $\tau'=\tau$ and $v'[i]=0$ for each $i$
(i.e.\ ${\bf v'}$ is a component of the label of the root of
$\overline {G_\tau}$) then the construction gives, straightforwardly the
result.

Suppose that there is a non af\/f\/ine edge
\begin{center}
\begin{tikzpicture}
\GraphInit[vstyle=Shade]
    \tikzstyle{VertexStyle}=[shape = rectangle,
draw
]
\SetUpEdge[lw = 1.5pt,
color = orange,
 labelcolor = gray!30,
 style={post},
labelstyle={sloped}
]
\tikzset{LabelStyle/.style = {draw,
                                     fill = white,
                                     text = black}}
\tikzset{EdgeStyle/.style={post}}
\Vertex[x=0, y=0,
 L={\tiny${\bf v''}$}]{x}
\Vertex[x=4, y=0,
 L={\tiny$\bf v'$}]{y}
\Edge[label={\tiny$s_i$ }](x)(y)
\end{tikzpicture}
\end{center}
in $\overline {G_\tau}$.
By induction ${\bf v''}\rest M$ labels a vertex of $G_{\lambda\rest M}$.
We verify that
\[({\bf v''}s_i)\rest M=({\bf v''})\rest Ms_i={\bf v'}\rest M\]
Hence, ${\bf v'}\rest M$ labels a vertex of $G_{\lambda\rest M}$ and the
assertion (2) is recovered.

Suppose now, that there is a af\/f\/ine edge \begin{center}
\begin{tikzpicture}
\GraphInit[vstyle=Shade]
    \tikzstyle{VertexStyle}=[shape = rectangle,
draw
]
\SetUpEdge[lw = 1.5pt,
color = orange,
 labelcolor = gray!30,
 style={post},
labelstyle={sloped}
]
\tikzset{LabelStyle/.style = {draw,
                                     fill = white,
                                     text = black}}
\tikzset{EdgeStyle/.style={post}}
\Vertex[x=0, y=0,
 L={\tiny${\bf v''}$}]{x}
\Vertex[x=4, y=0,
 L={\tiny$\bf v'$}]{y}
\Edge[label={\tiny$\Psi'$ }](x)(y)
\end{tikzpicture}
\end{center}
in $\overline {G_\tau}$.
By induction ${\bf v''}\rest M$ labels a vertex of $G_{\lambda\rest M}$.
We verify that
\[({\bf v''}\Psi')\rest M=({\bf v''}\Psi s_M\dots s_N)\rest M=({\bf
v''})\rest M\Psi={\bf v'}\rest M.\]
Hence, ${\bf v'}\rest M$ labels a vertex of $G_{\lambda\rest M}$ and the
assertion (3) is recovered.

Conversely, if ${\bf v'}$ labels a vertex of $G_{\lambda\rest M}$, there
exists a vertex labeled by ${\bf v'}^{(N)}$ in $\overline {G_\tau}$
verifying ${\bf v'}^{(N)}\rest M={\bf v'}$. Indeed, suppose ${\bf
v'}=(\tau',\zeta',v',\sigma')$ then ${\bf
v'}^{(N)}=({\tau'}^{(N)},{\zeta'}^{(N)},{v'
}^{(N)},{\sigma'}^{(N)})$, where ${\tau'}^{(N)}$ is obtained from
$\tau'$ by adding the nodes of $\tau$ labeled by $M+1,\dots,N$,
${v'}^{(N)}=[v'[1],\dots,v'[M],0,\dots,0]$,
${\zeta'}^{(N)}=\zeta_{{v'}^{(N)},{\tau'}^{(N))}}$ and
${\sigma'}^{(N)}=\sigma_{{v'}^{(N)}}$. Furthermore if ${\bf v'}s_i={\bf
v''}$ then ${\bf v'}^{(N)}s_i={\bf v''}^{(N)}$ and if ${\bf v'}\Psi={\bf
v''}$ then ${\bf v'}^{(N)}\Psi'={\bf v''}^{(N)}$.
This concludes the proof.
\end{proof}

\begin{Example}
Consider in Fig.~\ref{G332} the restriction problem for
\[\tau=\begin{array}{ccc}5&2\\7&3&1\\8&6&4 \end{array}\] and $M=3$.
\begin{figure}[t]\centering
\begin{tikzpicture}%
\GraphInit[vstyle=Shade]
    \tikzstyle{VertexStyle}=[shape = rectangle,
                             draw
]

\Vertex[x=3, y=-0.5,
 L={\tiny${52\ \atop {731\atop864}}, [000\dots],\atop
[1,-1,0,\dots],123\dots$},style={shape=circle}
]{x2}
\SetUpEdge[lw = 1.5pt,
color = orange,
 labelcolor = gray!30,
 labelstyle = {draw,sloped},
 style={post}
]

\tikzset{LabelStyle/.style = {draw,
                                     fill = yellow,
                                     text = red}}
\Vertex[x=-3, y=0, L={\tiny${51\ \atop {732\atop864}}, [000\dots],\atop
[-1,1,0,\dots],123\dots$}]{x1}
\Edge[label={\tiny$s_1$},style={post,in=-10,out=190},color=blue](x2)(x1)
\Vertex[x=-2, y=2, L={\tiny${51\ \atop {732\atop864}}, [001\dots],\atop
[1,0,\alpha -1,\dots],231\dots$}]{y1}
\Vertex[x=-4, y=4, L={\tiny${51\ \atop {732\atop864}}, [010\dots],\atop
[1,\alpha -1,0,\dots],213\dots$}]{y2}
\Vertex[x=-6, y=6, L={\tiny${51\ \atop {732\atop864}}, [100\dots],\atop
[\alpha-1,1,0,\dots],123\dots$}]{y3}
\Vertex[x=2, y=2, L={\tiny${52\ \atop {731\atop864}}, [001\dots],\atop
[-1,0,\alpha +1,\dots],231\dots$}]{z1}
\Vertex[x=4, y=4, L={\tiny${52\ \atop {731\atop864}}, [010\dots],\atop
[-1,\alpha +1,0,\dots],213\dots$}]{z2}
\Vertex[x=6, y=6, L={\tiny${52\ \atop {731\atop864}}, [100\dots],\atop
[\alpha+1,-1,0,\dots],123\dots$}]{z3}

\Edge[label={\tiny$s_2$}](y1)(y2)
\Edge[label={\tiny$s_1$}](y2)(y3)
\Edge[label={\tiny$s_2$}](z1)(z2)
\Edge[label={\tiny$s_1$}](z2)(z3)
\Edge[label={\tiny$\Psi'$},style={post,in=-80,out=80},color=green](x1)(y1)
\Edge[label={\tiny$\Psi'$},style={post,in=-80,out=80},color=green](x2)(z1)

\Vertex[x=-2, y=8, L={\tiny${51\ \atop {732\atop864}}, [011\dots],\atop
[0,\alpha -1,\alpha +1,\dots],312\dots$}]{xy1}
\Vertex[x=2, y=8, L={\tiny${52\ \atop {731\atop864}}, [011\dots],\atop
[0,\alpha +1,\alpha -1,\dots],312\dots$}]{xz1}
\Vertex[x=-4, y=10, L={\tiny${51\ \atop {732\atop864}}, [101\dots],\atop
[\alpha-1,0\alpha +1,\dots],132\dots$}]{xy2}
\Vertex[x=4, y=10, L={\tiny${52\ \atop {731\atop864}}, [101\dots],\atop
[\alpha+1,0\alpha -1,\dots],132\dots$}]{xz2}
\Vertex[x=-6, y=12, L={\tiny${51\ \atop {732\atop864}}, [002\dots],\atop
[1,0,2\alpha-1,\dots],231\dots$}]{xy3}
\Vertex[x=6, y=12, L={\tiny${52\ \atop {731\atop864}}, [002\dots],\atop
[-1,0,2\alpha+1,\dots],231\dots$}]{xz3}
\Edge[label={\tiny$\Psi'$},style={post,in=-80,out=80},color=green](y1)(xy1)
\Edge[label={\tiny$\Psi'$},style={post,in=-80,out=80},color=green](y2)(xy2)
\Edge[label={\tiny$\Psi'$},style={post,in=-80,out=80},color=green](y3)(xy3)
\Edge[label={\tiny$\Psi'$},style={post,in=-80,out=80},color=green](z1)(xz1)
\Edge[label={\tiny$\Psi'$},style={post,in=-80,out=80},color=green](z2)(xz2)
\Edge[label={\tiny$\Psi'$},style={post,in=-80,out=80},color=green](z3)(xz3)

\Edge[label={\tiny$s_2$},style={post,in=-10,out=190},color=blue](xz1)(xy1)
\Edge[label={\tiny$s_1$}](xy1)(xy2)
\Edge[label={\tiny$s_1$}](xz1)(xz2)
\Edge[label={\tiny$\Psi'$},style={post,in=-80,out=80},color=green](y3)(xy3)

\Vertex[x=-2, y=12, L={\tiny${51\ \atop {732\atop864}}, [110\dots],\atop
[\alpha-1,\alpha +1,0,\dots],123\dots$}]{xy4}
\Vertex[x=2, y=12, L={\tiny${52\ \atop {731\atop864}}, [110\dots],\atop
[\alpha+1,\alpha -1,0,\dots],123\dots$}]{xz4}
\Edge[label={\tiny$s_2$}](xy2)(xy4)
\Edge[label={\tiny$s_2$}](xz2)(xz4)
\Edge[label={\tiny$s_1$},style={post,in=-10,out=190},color=blue](xz4)(xy4)

\Vertex[x=-4, y=14, L={\tiny${51\ \atop {732\atop864}}, [020\dots],\atop
[1,2\alpha-1,0,\dots],213\dots$}]{xy5}
\Vertex[x=-2, y=16, L={\tiny${51\ \atop {732\atop864}}, [200\dots],\atop
[2\alpha-1,1,0,\dots],123\dots$}]{xy6}
\Vertex[x=4, y=14, L={\tiny${52\ \atop {731\atop864}}, [020\dots],\atop
[-1,2\alpha+1,0,\dots],213\dots$}]{xz5}
\Vertex[x=2, y=16, L={\tiny${52\ \atop {731\atop864}}, [200\dots],\atop
[2\alpha+1,-1,0,\dots],123\dots$}]{xz6}

\Edge[label={\tiny$s_2$}](xz3)(xz5)
\Edge[label={\tiny$s_1$}](xz5)(xz6)
\Edge[label={\tiny$s_2$}](xy3)(xy5)
\Edge[label={\tiny$s_1$}](xy5)(xy6)
\end{tikzpicture}
\caption{\label{G332} The f\/irst vertices of the subgraph of $G_{332}$
associated to the restriction for
$M=3$ whose root is $\mbox{\tiny $52\ \atop {731\atop864}$}$ with the edges $\Psi'$ added.}
\end{figure}
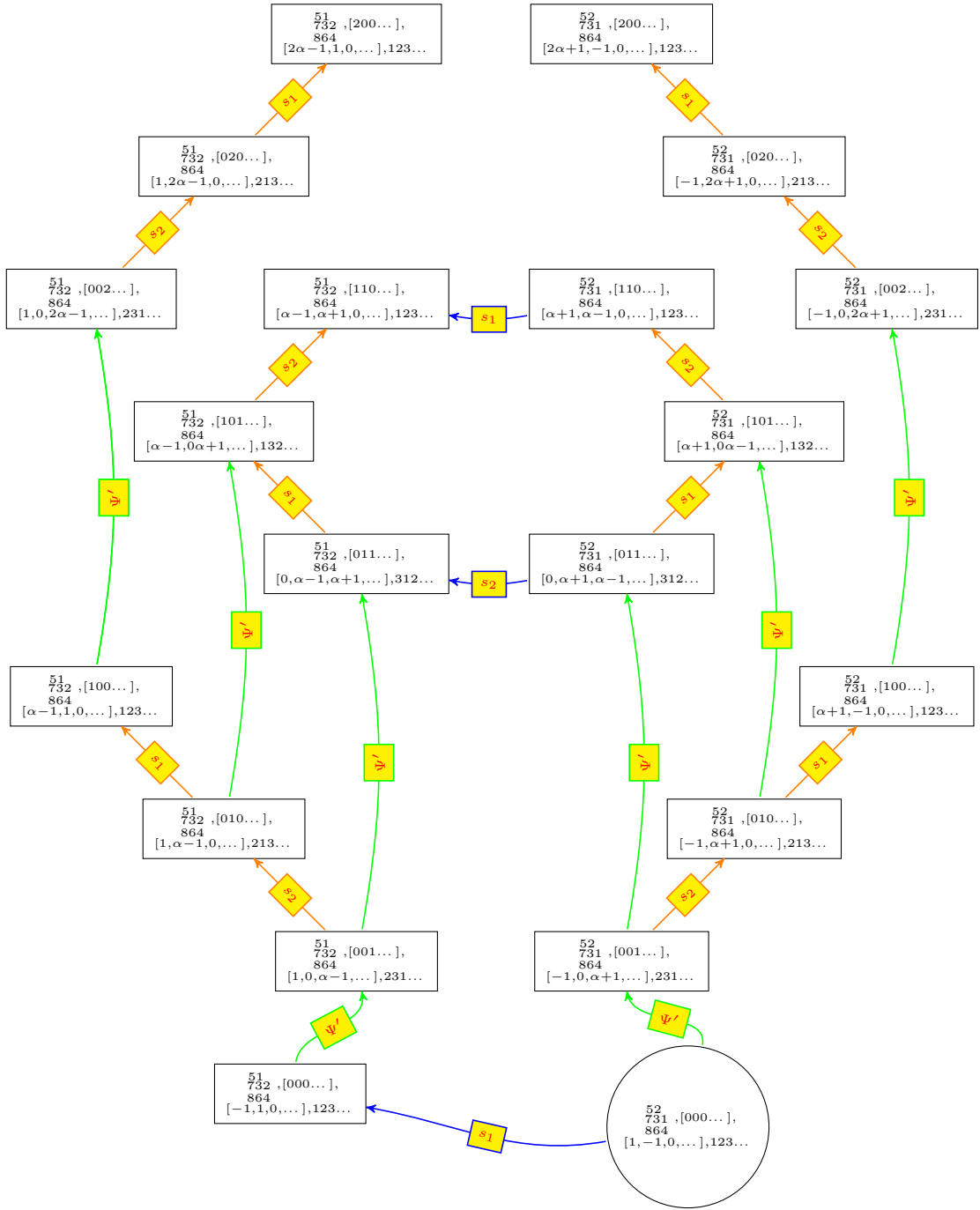
The subgraph of $\overline{G_{\mbox{\tiny$\begin{array}{ccc}5&2\\7&3&1\\8&6&4
\end{array}$}}}$ obtained using only the root and the arrows labeled by
$\Psi'$, $s_1$ and $s_2$ is isomorphic to the graph $\overline{G_{2\
\atop 31}}$ (see Fig.~\ref{G21}).
\end{Example}

\subsection{Restrictions on tableaux\label{ss6.2RestTab}}

In the sequel, as in \cite{Macdo}, we will denote a skew
partition by $\lambda/\mu$.

Let $\tau$ be a RST of shape
$\lambda=[\lambda_1,\dots,\lambda_k]$ and $M<N$. Consider $P_\tau$ as a
polynomial in $\C[t_1,\dots,t_M][t_{M+1},\dots ,t_N]$. Let $\mu$ the
sub-partition of $\lambda$ which is the shape of the RST obtained by
removing the nodes labeled by $1,\dots, M$ in $\tau$ and denote by
$\tau^{(M)}$ the associated RST. Consider also the skew-RST $\tau\rest
M$ of shape $\lambda/\mu$ obtained by removing the nodes labeled
$M+1,\dots, N$ in~$\tau$. Let $C(\tau,M):=\{\rho :
\rho^{(M)}=\tau^{(M)}\}$ and $T(\tau,M):=\{\rho\rest M : \rho\in
C(\rho,M) \} $. To each skew-RST $\rho\rest M$ in $T(\tau,M)$ we
associate the polynomial $P'_{\rho\rest M}$ which is the coef\/f\/icient of
the monomial $\prod\limits_{i=M+1}^Nt_i^{\lambda_{\cl(i,\rho)}-\rw(i,\rho)}$ in
$P_{\rho}$ (recall that $\cl(i,\tau)$ is the column of $\tau$
containing $i$ and $\rw(i,\tau)$ is the row of $\tau$ containing $i$).

\begin{Example}
Consider $\tau=\begin{array}{cc}\blue 2\\ \red4&\blue1 \\
\red5&\blue3\end{array}$ and $M=3$. Then $\tau^{(3)}={\red 4\atop 5}$,
$\tau\rest 3=\begin{array}{cc}\blue 2&\\&\blue1\\&\blue3\end{array}$,
\[C(\tau,3)=\left\{\begin{array}{cc}\blue 3&\\ \red4&\blue1
\\\red5&\blue2\end{array},\begin{array}{cc}\blue 2& \\ \red4&\blue1 \\
\red5&\blue3\end{array}, \begin{array}{cc}\blue 1&\\ \red4&\blue2 \\
\red5&\blue3\end{array} \right\}\] and
 \[T(\tau,3)=\left\{\begin{array}{cc}\blue 3&\\ &\blue1 \\
&\blue2\end{array},\begin{array}{cc}\blue 2& \\ &\blue1 \\
&\blue3\end{array}, \begin{array}{cc}\blue 1&\\ &\blue2 \\
&\blue3\end{array} \right\}.\]
The polynomial $P'_{\begin{array}{cc}\blue
2&\\&\blue1\\&\blue3\end{array}}$ is the coef\/f\/icient of $t_5^2t_4$ in
$P_{\begin{array}{cc}\blue 2\\ \red4&\blue1 \\ \red5&\blue3\end{array}}$.
\end{Example}

Suppose that $\gamma:=\lambda/\mu$ is a partition ($\tau\rest M$ is a RST).
When $i<M$, since $s_i$ does not act on the variables
$t_{M+1},\dots,t_N$, the Murphy rules (equation~\ref{Murphy}) give
\begin{equation}\label{MurphRest}
P'_{\rho\rest M}s_i=\left\{
\begin{array}{ll}
  b_\rho[i]P'_{\rho\rest M}&\mbox{if} \ b_\rho[i]^2=1,\\
b_\rho[i]P'_{\rho\rest M}+P'_{\rho^{(i,i+1)}\rest M}&\mbox{if} \ 0<b_\rho[i]\leq \frac12,\\
b_\rho[i]P'_{\rho\rest M}+(1-b_\rho[i]^2)P'_{\rho^{(i,i+1)}\rest
M}&\mbox{otherwise}.
\end{array}
\right.
\end{equation}
Since the action $s_i$ ($i<M$) on the RST commutes with the restriction
$\rest M$, one has $P'_{\rho^{(i,i+1)}\rest M}=P'_{(\rho\rest
M)^{(i,i+1)}}$ and Proposition~\ref{eigenMurph} and equation~(\ref{MurphRest})
imply that the polynomials $P'_{\rho\rest M}$ are simultaneous
eigenfunctions of the Jucys--Murphy operators
\[
 \omega_i^{(M)}=\sum_{i+1}^Ms_{ij},
\]
with eigenvalues $\CT_{\rho\rest M}[i]$. Since the $P'_{\rho\rest M}$
for $\rho\in C(\tau,M)$ span a polynomial representation of the
symmetric group $\S_M$ with minimal degree, the polynomials
$P'_{\rho\rest M}$ are equal up to a global multiplicative coef\/f\/icient
to the polynomials $P_{\rho\rest M}$.

{\samepage To summarize:
\begin{Proposition}
When $\tau\rest M$ is a RST, the coefficient of
$\prod\limits_{i=M+1}^Nt_i^{\lambda_{\cl(i,\tau)}-\rw(i,\tau)}$ in $P_{\tau}$
is proportional to $P_{\tau\rest M}$.
\end{Proposition}

}

\begin{Example}
\rm The coef\/f\/icient of $t_9^2t_8^2t_7t_6$ in $P_{51\ \ \atop{732\
\atop9864}}$ equals
\[
\frac16  t_1-\frac1{12} t_2 -\frac1{12} t_3 =\frac16
P_{1\ \atop 32}.
\]
\end{Example}

\subsection{Restrictions on Jack polynomials\label{ss6.3RestJack}}
Consider the linear map $\rest M$ which sends each $x_i$ to
$0$ when $i>M$ and each $\tau$ to $\tau\rest M$.
\begin{Theorem}
Let $\tau$ have property $R(M)$. Then
\[
 J_{[v[1],\dots,v[M],0,\dots,0],\tau}\rest
M=J_{[v[1],\dots,v[M]],\tau\rest M}.
\]
\end{Theorem}

\begin{proof} From Proposition~\ref{PropRest} the graphs
$\overline{G_{\tau\rest M}}$ and $G_{{\rm shape}(\tau\rest M)}$ are the
same.
Remark that:
\begin{enumerate}\itemsep=0pt
 \item The result is correct for the roots:
\[ J_{([0,\dots,0],\tau_{\rm root})}\rest M=J_{[0,\dots,0],\tau_{\rm root}\rest M}\]
where $\tau_{\rm root}$ denotes the RST obtained from $\tau$ by replacing
subtableau constituted with the nodes labeled $1,\dots,M$ by $\tau_{{\rm
shape}(\tau\rest M)}$.
 \item The action of the edges are compatible with the restriction:
\begin{enumerate}\itemsep=0pt
\item The non-af\/f\/ine edges: We use only the dif\/ference between the
content of two  boxes. Since, the values of
the spectral vector $\zeta$ in $\overline{G_{\tau_{\rm root}\rest M}}$ are
obtained from the values of the spectral vector in $G_{\tau_{\rm root}}$ by
adding to each component the same integer. Hence, the dif\/ferences are
the same and then the action of the non-af\/f\/ine edges are the same.
\item The af\/f\/ine edges: One has to verify that
 \begin{gather*}
 J_{[v[1],\dots,v[M],0,\dots,0],\tau}\Psi_N(s_{N-1}
\otimes s_{N-1}+(*))\cdots (s_M \otimes s_M+(*))\rest M
 \\
\qquad{} =J_{[v[1],\dots,v[M],0,\dots,0],\tau}\rest M\Psi_M,
\end{gather*}
where $(*)$ denote the correct rational numbers corresponding to the
edges of the Yang--Baxter graph and $\Psi_M$ means that one applies the
operator $\Psi$ for an alphabet of size $M$.
This identity is easy to obtained from the construction: since $\Psi$
gives a~polynomial whose a factor is $x_N$, the only non vanishing part of
\[J_{[v[1],\dots,v[M],0,\dots,0],\tau}\Psi_N(s_{N-1} \otimes
s_{N-1}+(*))\cdots (s_M \otimes s_M+(*))\rest M\]
is
 \[J_{[v[1],\dots,v[M],0,\dots,0],\tau}\Psi_N(s_{N-1} \cdots
s_M \otimes s_{N-1}\cdots s_M)\rest M=
J_{[v[1],\dots,v[M],0,\dots,0],\tau}\Psi_M \rest M,\]
the last part of the proof follows from the commutation $\Psi_M \rest
M=\rest M\Psi_M$.
\end{enumerate}
\end{enumerate}

This shows that the polynomials
$J_{[v[1],\dots,v[M],0,\dots,0],\tau}\rest M$ are inductively generated
following the same Yang--Baxter graph as the polynomials
$J_{[v[1],\dots,v[M]],\tau\rest M}$ with the same initial conditions.
Hence $J_{[v[1],\dots,v[M],0,\dots,0],\tau}\rest
M=J_{[v[1],\dots,v[M]],\tau\rest M}$ as expected.
\end{proof}

\section{Shifted vector-valued Jack polynomials}\label{s7shift}

\subsection{Knop and Sahi operators for vector-valued polynomials}\label{ss7.1Knop}
Let us def\/ine the following operators which are the vector-valued
versions of the operators def\/ined in \cite{KS}:
\begin{enumerate}\itemsep=0pt
\item $\varsigma_i:=\partial_i\otimes 1+s_i\otimes s_i,$ where
$\partial_i:=\partial_{i,i+1}=(1-s_i)\frac1{x_i-x_{i+1}}$ is a divided
dif\/ference.
\item Denote by $\Phi$ the operator sending each $x_i$ to $x_{i-1}$ for
$i>1$ and $x_1$ to $x_N-\alpha$ and $
T:=\Phi\otimes (s_1s_2\cdots s_{N-1})$.
\item
$\hat \Psi:=T(x_N+N-1).$
\end{enumerate}

\begin{Proposition}
The operators $\varsigma_i$ satisfy the braid relations
\[
\varsigma_i\varsigma_{i+1}\varsigma_i=\varsigma_{i+1}\varsigma_{i}\varsigma_{i+1},\qquad
\varsigma_i\varsigma_j=\varsigma_j\varsigma_i,\qquad |i-j|>1,
\]
and the relations between the $\varsigma_i$ and the multiplication by
the indeterminates are given by the Leibniz rules:
\[
 x_i\varsigma_i=\varsigma_{i+1}x_i+1,\qquad
x_{i+1}\varsigma_i=\varsigma_{i}x_i-1,\qquad
x_j\varsigma_i=\varsigma_ix_j,\qquad j\neq i,i+1.
\]
\end{Proposition}

\begin{proof}
The Leibniz rules are straightforward from the def\/inition while the
braid relations are a direct consequence of
 \begin{enumerate}\itemsep=0pt
\item[1)] the braid relations on the $s_i$ and the braid relations on the
$\partial_i$,

\item[2)]
$\partial_{i+1}\partial_is_i+s_{i+1}\partial_i\partial_{i+1}=\partial_is_{i+1}\partial_i$,

\item[3)]
$\partial_i\partial_{i+1}s_i+s_i\partial_{i+1}\partial_i=\partial_{i+1}s_i\partial_{i+1}$,

\item[4)] $s_i\partial_{i+1}s_i=s_{i+1}\partial_is_{i+1}$,

\item[5)] $\partial_is_{i+1}s_i=s_{i+1}s_i\partial_{i+1}$,

\item[6)]  $s_is_{i+1}\partial_i=\partial_{i+1}s_is_{i+1}$.\hfill \qed
\end{enumerate}
 \renewcommand{\qed}{}
\end{proof}

Since the $\varsigma_i$ verify the braid relations, they
 realize the braid group: given a
permutation $\omega\in\S_N$ and a reduced decomposition
$\omega=s_{i_1}\cdots s_{i_k}$, the product
$\varsigma_{i_1}\cdots\varsigma_{i_k}$ is independent of the choice of
the reduced decomposition. We will denote
$\varsigma_\omega:=\varsigma_{i_1}\cdots\varsigma_{i_k}$.

Furthermore, the algebra generated by the $\varsigma_i$ and the $x_j$ is
isomorphic to the \emph{degenerate Hecke
affine algebra}  generated by the operators $s_i+\partial_i$ and
the variables.

 Our goal is to f\/ind a basis of simultaneous eigenvectors of
the following operators
\[
\hat \xi_i:=x_i+N-1-\varsigma_i\cdots
\varsigma_{N-1}\hat\Psi\varsigma_{1}\cdots \varsigma_{i-1}.
\]

These operators commute and play the role of Cherednik elements for our
representation of the degenerate Hecke af\/f\/ine algebra. As a consequence,
one has the following relations:
\begin{Proposition}\label{CommShif}\qquad
\begin{enumerate}\itemsep=0pt
\item[$1.$] $\varsigma_i\hat\xi_{i+1}=\hat\xi_{i}\varsigma_i-1$,
\item[$2.$] $\varsigma_i\hat\xi_i=\hat\xi_{i+1}\varsigma_i+1$,
\item[$3.$] $\varsigma_i\hat\xi_{j}=\hat\xi_{j}\varsigma_i$ for $j\neq i,i+1$,
\item[$4.$] $\hat \Psi\hat\xi_i=\hat\xi_{i-1}\hat \Psi $ for $i\neq 1$,
\item[$5.$] $\hat \Psi\hat\xi_1=(\hat\xi_N-\alpha)\hat \Psi$.
\end{enumerate}
\end{Proposition}

Furthermore, the RST are simultaneous eigenfunctions of the operators
$\hat\xi_i$. More precisely
\begin{Proposition}\label{hat_itab}
\[\tau\hat\xi_i=\CT_\tau[i]\tau.
\]
\end{Proposition}

\begin{proof}\looseness=-1
The action of $\hat\xi_i$ on polynomials with degree $0$ in
the $x_i$ equals the action of the opera\-tors~$\tilde\xi_i$. Hence, the
result follows from the non-shifted version of the equality (Proposi\-tion~\ref{U_itab}).
\end{proof}

Straightforwardly, the operators $\varsigma_i$ and $\hat\Psi$ are
compatible with the  leading monomials
 in the following sense:
\begin{Proposition}\label{shiftdom}
Suppose that $P$ a polynomial such that its highest degree component has
 the leading monomial $ x^{v,
\tau}$ then
\begin{enumerate}\itemsep=0pt
\item[$1.$] If $v[i]< v[i+1]$ then the highest degree component of
$P\varsigma_i$ has the
leading monomial $x^{vs_{i},\tau}$.
\item[$2.$] The highest degree component of $P\hat\Psi$ has the
 leading monomial
 $x^{v\Psi,\tau}$.
\end{enumerate}
\end{Proposition}

\subsection[The Yang-Baxter graph]{The Yang--Baxter graph}\label{ss7.2YBShift}

Let $\lambda$ be a partition and $G_\lambda$ the associated graph. We
construct the set of the polynomials
$\big(\hat J_{\cal P}\big)_{\cal P \ {\rm path \ in} \ G_\lambda}$ using
the following recurrence rules:
\begin{enumerate}\itemsep=0pt
 \item $\hat J_{[]}:=1\otimes\tau_{\lambda}$.
 \item If ${\cal P}=[a_1,\dots,a_{k-1},s_i]$ then
\[\hat J_{\cal P}:=\hat
J_{[a_1,\dots,a_{k-1}]}\left(\varsigma_i+\frac1{\zeta[i+1]-\zeta[i]}\right),\]
where the vector $\zeta$ is def\/ined by
\[
\big(\tau_\lambda,\CT_{\tau_\lambda},0^N,[1,2,\dots,N]\big)a_1 \cdots
a_{k-1}=(\tau,\zeta,v,\sigma).
\]
 \item ${\cal P}=[a_1,\dots,a_{k-1},\Psi]$ then then
\[
\hat J_{\cal P}=\hat J_{[a_1,\dots,a_{k-1}]}\hat\Psi.
\]
\end{enumerate}

 As expected one obtains

\begin{Theorem}\label{thshifted}
Let ${\cal P}=[a_0,\dots,a_k]$ be a path in $G_\lambda$ from the root to
$(\tau,\zeta,v,\sigma)$.
The polyno\-mial~$\hat J_{\cal P}$ is a simultaneous eigenfunction of the
operators $\hat\xi_i$ whose leading monomial in the highest degree component
is $x^{v,\tau}$. Furthermore, the eigenvalue of $\hat\xi_i$ associated
to $\hat J_{\cal P}$ equals $\zeta[i]$.

 Consequently $\hat J_{\cal P}$ does not depend on the path, but
only on the end point $(\tau,\zeta,v,\sigma)$, and will be denoted by $\hat
J_{v,\tau}$.
 The
  family $(\hat J_{v,\tau})_{v,\tau}$ forms a basis of $M_\lambda$ of
simultaneous eigenfunctions of the Cherednik operators.

Furthermore, if $\cal P$ leads to $\varnothing$ then $\hat J_{\cal P}=0$.
\end{Theorem}

\begin{proof}
The proof goes as in Theorem~\ref{thJack} using respectively
Propositions~\ref{CommShif},~\ref{hat_itab} and~\ref{shiftdom} instead
of Propositions~\ref{Cherednik_com},~\ref{U_itab} and~\ref{dominance}.\end{proof}
In consequence, we will consider the family of polynomials
$(\hat J_{v,\tau})_{v,\tau}$ indexed by pairs $(v,\tau)$ where $v\in
\N^N$ is a weight and $\tau$ is a tableau.
\begin{Example}
\rm Let again $\tau=\begin{array}{cc}2\\3&1\end{array}$ and consider the
Yang--Baxter Graph $G_{2\atop 31}$ (see Fig.~\ref{hatG}).

\begin{figure}[t]
\centering
\begin{tikzpicture}
\GraphInit[vstyle=Shade]
    \tikzstyle{VertexStyle}=[shape = rectangle,
draw
]

\Vertex[style={
shape=circle},x=-6, y=0,
 L={$\hat J_{000,{1\ \atop 32}}$}]{x2}

\SetUpEdge[lw = 1.5pt,
color = orange,
 labelcolor = gray!30,
 style={post},
labelstyle={sloped}
]
 \tikzstyle{TempStyle}=[double = green,
                           double distance = 1pt]

\tikzset{LabelStyle/.style = {draw,
                                     fill = white,
                                     text = black}}
\tikzset{EdgeStyle/.style={post}}
\Vertex[x=-2, y=0, L={\tiny $\hat J_{001,{1\ \atop 32}}$}]{z1}
\Vertex[x=-2, y=-3, L={\tiny$\hat J_{010,{1\ \atop 32}}$}]{z2}
\Vertex[x=-2, y=-6, L={\tiny$\hat J_{100,{1\ \atop 32}}$}]{z3}
\Edge[label={\tiny $\varsigma_2+\frac1{1+\alpha}$}](z1)(z2)
\Edge[label={\tiny $\varsigma_1+\frac1{2+\alpha}$}](z2)(z3)
\Edge[,label={\tiny $\hat\Psi$}](x2)(z1)
\Vertex[x=2, y=0, L={\tiny$\hat J_{011,{1\ \atop 32}}$}]{xz1}
\Vertex[x=2, y=-3, L={\tiny$\hat J_{101{1\ \atop 32}}$}]{xz2}
\Vertex[x=2, y=-8, L={\tiny$\hat J_{002,{1\ \atop 32}}$}]{xz3}
\Edge[label={\tiny $\hat\Psi$}](z1)(xz1)
\Edge[label={\tiny $\hat\Psi$}](z2)(xz2)
\Edge[label={\tiny $\hat\Psi$}](z3)(xz3)
\Edge[label={\tiny $\varsigma_1+\frac1{1+\alpha}$}](xz1)(xz2)

\Vertex[x=2, y=-6, L={\tiny $\hat J_{110,{1\ \atop 32}}$}]{xz4}
\Edge[label={\tiny $\varsigma_2+\frac1{\alpha-1}$}](xz2)(xz4)
\Vertex[x=2, y=-11, L={\tiny $\hat J_{020,{1\ \atop 32}}$}]{xz5}
\Vertex[x=2, y=-14, L={\tiny $\hat J_{200,{1\ \atop 32}}$}]{xz6}

\Edge[label={\tiny $\varsigma_2+\frac1{1+2\alpha}$}](xz3)(xz5)
\Edge[label={\tiny$\varsigma_1+\frac1{2+2\alpha}$}](xz5)(xz6)

\end{tikzpicture}
\caption{\label{hatG} First values of $\hat J_{v,{2\ \atop 31}}$.}
\end{figure}
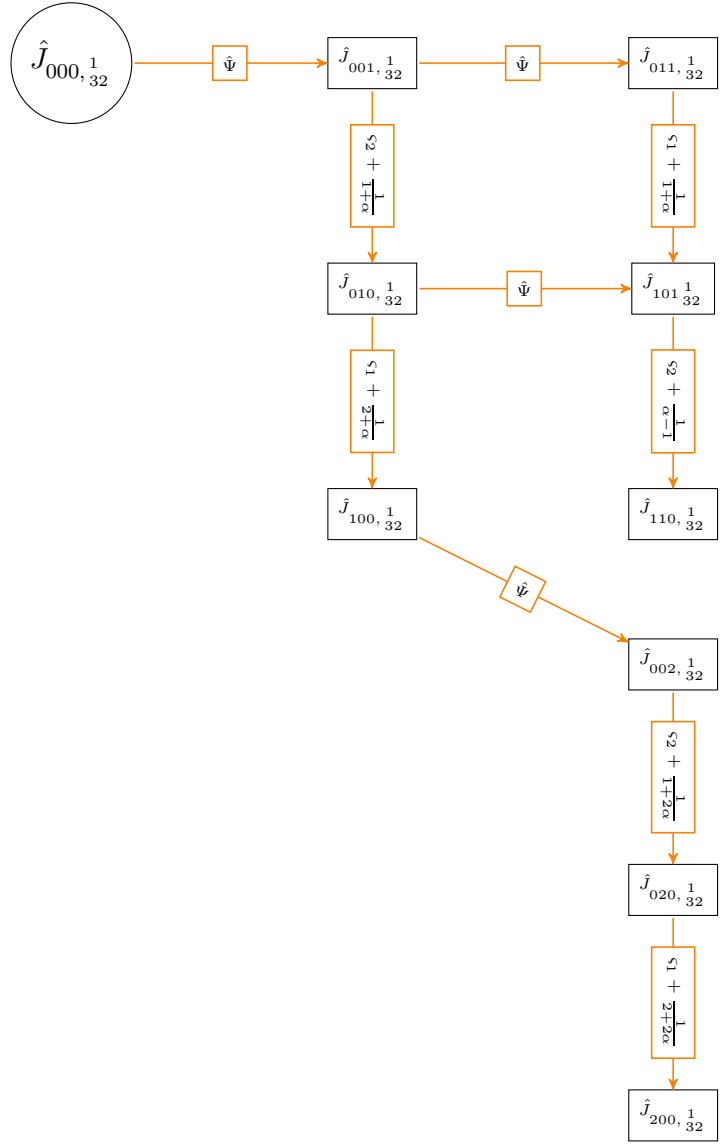

The eigenvalues of $\tau$ are $[1,-1,0]$. Let
$\tau'=\begin{array}{cc}1\\3&2\end{array}$, one has
\[
\hat J_{[001],\tau} = \tau\hat\Psi
 = \left( x_3 +2 \right) \otimes \left( -\frac12 \tau+\tau' \right)
\]
with associated eigenvalues $[-1,0,\alpha+1]$
and
\begin{gather*}
\hat J_{[010],\tau} = \hat
J_{[001],\tau}\left((s_2+\partial_{2,3})\otimes
s_2+\frac1{\alpha+1}\right) \\
\hphantom{\hat J_{[010],\tau}}{}
 = \frac12 \frac {\left( 3 \alpha+1+x_2+\alpha x_2-x_3 \right)}{\alpha+1}\tau+\frac {\left(x_3+3+x_2+\alpha x_2+\alpha
 \right)}{\alpha+1}\tau'
\end{gather*}
with eigenvalues $[-1,\alpha+1,0]$.
\end{Example}

The non-shifted vector-valued Jack polynomials can be recovered easily
from the shifted one:
\begin{Proposition}
The restriction of $\hat J_{v,\tau}$ to its component of top degree
equals $J_{v,\tau}$.
\end{Proposition}
\begin{proof}
It suf\/f\/ices to remark that $\varsigma_i=s_i\otimes s_i+\varsigma_i^-$
and $\hat\Psi=\Psi+\Psi^-$ where $\varsigma^-$ and $\Psi^-$ are
operators which decrease the degree in the $x_j$. Hence, one computes
the component of the top degree of $\hat J_{v,\tau}$ following the
Yang--Baxter graph replacing $\varsigma_i$ by $s_i\otimes s_i$ and
$\hat\Psi$ by $\Psi$, that is the Yang--Baxter graph which was used to
obtain the polynomials $J_{v,\tau}$. This ends the proof.\end{proof}

\subsection{Symmetrization, antisymmetrization} \label{ss7.3Sym}

We will say that a polynomial is symmetric if it is invariant under the
action of $\varsigma_i$ for each $i<N-1$.
Denote also $\hat {\cal S}=\sum_{\omega\in\S_N} \varsigma_\omega$.
As for non shifted Jack, one has:
\begin{Theorem}\qquad
\begin{enumerate}\itemsep=0pt
\item[$1.$] Let $H_T$ be a connected component of $G_\lambda$. For each vertex
$(\tau,\zeta,v,\sigma)$ of $H_T$, the polyno\-mial~$\hat J_{v,\tau}\hat
{\cal S}$ equals $\hat J_{\lambda_T,\std T}\hat{\cal S}$ up to a
multiplicative constant.
\item[$2.$] One has $\hat J_{\lambda_T,\std T}\hat{\cal S}\neq 0$ if and only
if $H_T$ is $1$-compatible.
\item[$3.$] More precisely, when $H_T$ is $1$-compatible, the polynomial
\[
\hat J_{T}=\sum_{(\tau,\zeta,v,\sigma)\ {\rm vertex \ of} \ H_T}{\cal
E}_{v,\tau}\hat J_{v,\tau}
\]
is symmetric.
\end{enumerate}
\end{Theorem}

\begin{proof} The proof is identical to the non-shifted case.\end{proof}

In the same way, for the antisymmetrization, denoting $\hat {\cal
A}:=\sum_{\omega\in\S_N}(-1)^{\ell(\omega)} \varsigma_\omega$, we have:
\begin{Theorem}\qquad
\begin{enumerate}\itemsep=0pt
\item[$1.$] Let $H_T$ be a connected component of $G_\lambda$. For each vertex
$(\tau,\zeta,v,\sigma)$ of $H_T$, the polyno\-mial~$\hat J_{v,\tau}\hat
{\cal A}$ equals $\hat J_{\lambda_T,\std T}\hat {\cal A}$ up to a
multiplicative constant.
\item[$2.$] One has $\hat J_{\lambda_T,\std T}\hat {\cal A}\neq 0$ if and only
if $H_T$ is $(-1)$-compatible.
\item[$3.$] More precisely, when $H_T$ is $(-1)$-compatible the polynomial
\[
\hat J'_{T}=\sum_{(\tau,\zeta,v,\sigma)\ {\rm vertex \ of} \ H_T}{\cal
F}_{v,\tau}\hat J_{v,\tau}
\]
is antisymmetric.
\end{enumerate}
\end{Theorem}

\subsection{Propagation of vanishing properties}\label{ss7.4vanish}

Some phenomena of propagation of vanishing properties can be deduced
from the classical case (see \cite{Las}).
\begin{Lemma}[Lascoux]\label{classical_vanish}
Let $f(x,y)$ be a function of two variables. Suppose that $f(b,a)=0$
with $a\neq b$ then
\[f\left(s_{x,y}+\gamma.\partial_{x,y}+{\gamma\over b-a}\right)(a,b)=0.\]
\end{Lemma}

Indeed, similar properties occur for vector-valued polynomials:
\begin{Lemma}\label{Vanish_si}
Let $f\in \C[x,y]\otimes V_\lambda$ then
\[
 f\left(s_{x,y}\otimes s_{x,y}+\gamma\partial_{x,y}\otimes
1+\frac\gamma{b-a}\right)(a,b)=0
\]
when $f(b,a)=0$.
\end{Lemma}

\begin{proof} Write $f(x,y)=\sum_\tau f^\tau(x,y)\otimes\tau$ and remark
that $f(b,a)=0$ implies $f^\tau(a,b)=0$ for each~$\tau$. But
$fs_{x,y}=\sum_\tau f^\tau s_{x,y}\otimes\tau s_{x,y}$ and by Lemma~\ref{classical_vanish} $f^\tau s_{x,y}(a,b)=0$. Hence $fs_{x,y}(a,b)=0$.
In the same way, Lemma~\ref{classical_vanish} implies
$f^\tau\big(\partial_{x,y}+\frac1{b-a}\big)(a,b)=0$ and
then $f\big(\partial_{x,y}\otimes 1+\frac1{b-a}\big)(a,b)=0$. This proves the
result.\end{proof}

\begin{Example}
Consider the polynomial $\hat J_{001,{2\ \atop 31}}$. This polynomial
vanishes for $x_3=-2$:
$\hat J_{001,{2\ \atop 31}}(x_1,x_2,-2)=0$. Since
\[
\hat J_{010,{2\ \atop 31}}=\hat J_{001,{2\ \atop
31}}\left(\varsigma_2+\frac1{\alpha+1}\right)
\]
one has $\hat J_{010,{2\ \atop 31}}(x_1,-2,\alpha-1)=0$. Finally,
\[
\hat J_{100,{2\ \atop 31}}=\hat J_{010,{2\ \atop
31}}\left(\varsigma_1+\frac1{\alpha+2}\right)
\] implies
$\hat J_{100,{2\ \atop 31}}(-2,\alpha,\alpha-1)=0$.
\end{Example}

Denote by $\Phi^+$ the operator def\/ined by $[x_1,\dots,
x_N]\Phi^+=[x_2,x_3,\dots,x_N,x_1+\alpha]$.
The action of the af\/f\/ine operator $\hat\Psi$ propagates information
about vanishing properties:
\begin{Lemma}\label{VanishPsi}
One has
\[f\hat\Psi\big([a_1,\dots,a_N]\Phi^+\big)=0\]
when $f(a_1,\dots,a_N)=0$.
\end{Lemma}
Hence, for each pair $(v,\tau)$ one can compute at least one $N$-tuple
$(a^{v,\tau}_1,\dots,a^{v,\tau}_N)$ such that
$J_{v,\tau}(a^{v,\tau}_1,\dots,a^{v,\tau}_N)=0$.

\begin{Proposition}\label{propVanish}
Denote by ${\bf V}_{v,\tau}$ the vector whose $i$-th component is
\[
{\bf
V}_{v,\tau}[i]:=\alpha(v^{+}[1]-v[i])+(\CT_\tau[1]-\CT_\tau[\sigma_v[i]])-N+1.
\]

Let $m$ be the smallest integer such that $v[m]=\max\{v[i] : 1\leq i\leq
N\}$.
 One has
\[
\hat J_{v,\tau}(x_1,\dots,x_{m-1},{\bf V}_{v,\tau}[m],\dots,{\bf
V}_{v,\tau}[N]
)=0
\]
\end{Proposition}

\begin{proof}
Denote by $\tilde\Phi$ the operator def\/ined by
\begin{gather*}
{} [x_1,\dots,x_N]\tilde\Phi= \begin{cases}
[x_1,\dots,x_N]\Phi&\mbox{if} \ x_1=a\alpha+b\mbox{ with } a>0, \\
[x_2+\alpha,\dots,x_N+\alpha,x_1]&\mbox{otherwise}.
\end{cases}
\end{gather*}

Observe f\/irst that the vectors ${\bf V}_{v,\tau}$ are obtained by
substituting $\Psi$ by $\tilde\Phi$ in $G_\tau$:

\begin{Lemma}\label{Psi2tildePhi}
Recall the notation $[v_1,\dots,v_N]\Psi=[v_2,\dots,v_N,v_1+1]$. One has
\[
 {\bf V}_{v\Psi,\tau}={\bf V}_{v,\tau}\tilde\Phi.
\]
\end{Lemma}

\begin{proof}
Let $m$ be the smallest integer such that $v[m]=\max\{v[i] : 1\leq i\leq
N\}$.
If $m>1$ then $v^{+}[1]={v\Psi}^{+}[1]$. It follows that ${\bf
V}_{v\Psi,\tau}[i]={\bf V}_{v,\tau}[i+1]$ if $i<N$ and, since
$v\Psi[N]=v[1]+1$ one has ${\bf V}_{v\Psi,\tau}[N]={\bf
V}_{v,\tau}[1]-\alpha$. That is
${\bf V}_{v\Psi,\tau}={\bf V}_{v,\tau}\Phi$. If $m=1$ then
$v^{+}[1]={(v\Psi)}^{+}[1]-1$. Hence,
${\bf V}_{v\Psi,\tau}[i]={\bf V}_{v,\tau}[i+1]+\alpha$ for $i<N$ and
${\bf V}_{v\Psi,\tau}[N]={\bf V}_{v,\tau}[1]$. This ends the proof.
\end{proof}

 Hence, using Lemma~\ref{Psi2tildePhi}, a straightforward induction
shows that
\begin{equation}\label{V[m]}
{\bf V}_{v,\tau}[m]=1-N.
\end{equation}
Let us prove the proposition by induction on the length of the path from
$ \hat J_{0\dots01,\tau}$ to $ \hat J_{v,\tau}$.
First note that
\[
{\bf V}_{0\dots01,\tau}[N]=-N+1,
\]
implies $\hat J_{0\dots01,\tau}(x_1,\dots,x_{N-1},{\bf
V}_{0\dots01,\tau}[N])=0$,
since
\[\hat J_{0\dots01,\tau}=\hat J_{0\dots0,\tau}\hat\Psi=\hat
J_{0\dots0,\tau}T(x_N-N+1).\]
Suppose that there is an arrow
\begin{center}\begin{tikzpicture}
\GraphInit[vstyle=Shade]
    \tikzstyle{VertexStyle}=[shape = rectangle,
draw
]
\SetUpEdge[lw = 1.5pt,
color = orange,
 labelcolor = gray!30,
 style={post},
labelstyle={sloped}
]
\tikzset{LabelStyle/.style = {draw,
                                     fill = white,
                                     text = black}}
\tikzset{EdgeStyle/.style={post}}
\Vertex[x=0, y=0,
 L={\tiny$(\tau,\zeta',v',\sigma')$}]{x}
\Vertex[x=3, y=0,
 L={\tiny$(\tau,\zeta,v,\sigma)$}]{y}
\Edge[label={\tiny$s_i$}](x)(y)
\end{tikzpicture} \end{center}in $G_\tau$.
Let $m'$ be the smallest integer such that $v'[m']=\max\{v'[i] : 1\leq
i\leq N\}$. We have by induction,
\begin{equation}\label{induVanishsi} \hat
J_{v',\tau}\big(x_1,\dots,x_{m'-1},{\bf V}_{v',\tau}[m'],\dots,{\bf
V}_{v',\tau}[N]
\big)=0.\end{equation}
By hypothesis, $i\neq m'$, otherwise $v'[i]\geq v'[i+1]$ and
the arrow is not in $G_\tau$.
Furthermore we have:
\begin{gather*}
\hat J_{v,\tau} =
\hat J_{v',\tau}\left(s_i\otimes s_i+\partial_i\otimes
1+\frac1{(v'[i+1]-v'[i])\alpha+(\CT_\tau[\sigma_{v'}[i+1]]-\CT_\tau[\sigma_{v'}[i]])}\right)\nonumber\\
\phantom{\hat J_{v,\tau}}{} = \hat
J_{v',\tau}\left(s_i\otimes s_i+\partial_i\otimes
1+\frac1{(v[i]-v[i+1])\alpha+(\CT_\tau[\sigma_v[i]]-\CT_\tau[\sigma_v[i+1]])}\right). 
\end{gather*}
Consider three cases
\begin{enumerate}\itemsep=0pt
\item If $i<m'-1$ then Lemma~\ref{Vanish_si} implies
\[ \hat J_{v,\tau}\big(x_1,\dots,x_{m-1},{\bf V}_{v',\tau}[m'],\dots,{\bf
V}_{v',\tau}[N]
\big)=0.
\]
But ${\bf V}_{v',\tau}[m',m'+1,\dots, N]={\bf V}_{v,\tau}[m',m'+1,\dots,
N]$ and $m'=m$. This proves the result.
 \item If $i=m'-1$ then $m=m'-1$ and, as a special case of~(\ref{induVanishsi}), one has
\[\hat J_{v',\tau}\big(x_1,\dots,x_{m'-2},{\bf V}_{v,\tau}[m],{\bf
V}_{v',\tau}[m'],\dots,{\bf V}_{v',\tau}[N]
\big)=0.
\]
But
\begin{gather*}
{\bf V}_{v',\tau}[m']- {\bf V}_{v,\tau}[m] = {\bf V}_{v,\tau}[m+1]- {\bf
V}_{v,\tau}[m]\\
\phantom{{\bf V}_{v',\tau}[m']- {\bf V}_{v,\tau}[m]}{}
 = (v[m]-v[m+1])\alpha+(\CT_\tau[\sigma_v[m]]-\CT_\tau[\sigma_v[m+1]]).
\end{gather*}
The result is, now, a direct consequence of Lemma~\ref{Vanish_si}.
\item If $i>m'$ then ${\bf V}_{v',\tau}[m',m'+1,\dots, N]={\bf
V}_{v,\tau}[m',m'+1,\dots, N].s_i$ and $m'=m$.
Remarking that,
 \begin{gather*}
{\bf V}_{v',\tau}[i]- {\bf V}_{v,\tau}[i+1] = {\bf V}_{v,\tau}[i+1]-
{\bf V}_{v,\tau}[i]\\
\phantom{{\bf V}_{v',\tau}[i]- {\bf V}_{v,\tau}[i+1]}{}
 = (v[i]-v[i+1])\alpha+(\CT_\tau[\sigma_v[i]]-\CT_\tau[\sigma_v[i+1]]).
\end{gather*}
Lemma~\ref{Vanish_si} gives the result.
\end{enumerate}

Suppose that there is an arrow
\begin{center}\begin{tikzpicture}
\GraphInit[vstyle=Shade]
    \tikzstyle{VertexStyle}=[shape = rectangle,
draw
]
\SetUpEdge[lw = 1.5pt,
color = orange,
 labelcolor = gray!30,
 style={post},
labelstyle={sloped}
]
\tikzset{LabelStyle/.style = {draw,
                                     fill = white,
                                     text = black}}
\tikzset{EdgeStyle/.style={post}}
\Vertex[x=0, y=0,
 L={\tiny$(\tau,\zeta',v',\sigma')$}]{x}
\Vertex[x=3, y=0,
 L={\tiny$(\tau,\zeta,v,\sigma)$}]{y}
\Edge[label={\tiny $\Psi$}](x)(y)
\end{tikzpicture} \end{center}in $G_\tau$.
Let $m'$ be the smallest integer such that $v'[m']=\max\{v'[i] : 1\leq
i\leq N\}$. We have by induction,
\begin{equation}\label{induVanishPsi} \hat
J_{v',\tau}\big(x_1,\dots,x_{m'-1},{\bf V}_{v',\tau}[m'],\dots,{\bf
V}_{v',\tau}[N]
\big)=0.\end{equation}
We need to consider two cases:
\begin{enumerate}\itemsep=0pt
\item If $m'>1$ then, since, $\hat J_{v,\tau}= \hat
J_{v',\tau}.\hat\Psi$ Lemma~\ref{VanishPsi} and equation~(\ref{induVanishPsi}) imply
 \begin{gather*}
\hat J_{v,\tau}\big(x_2,\dots,x_{m},{\bf V}_{v,\tau}[m],\dots,{\bf
V}_{v,\tau}[N-1],x_1
\big)\\
\qquad{} =
\hat J_{v,\tau}\big(x_2,\dots,x_{m'-1},{\bf V}_{v',\tau}[m'],\dots,{\bf
V}_{v',\tau}[N],x_1
\big) =
0.
\end{gather*}
In particular,
\[
\hat J_{v,\tau}\big(x_1,\dots,x_{m-1},{\bf V}_{v,\tau}[m],\dots,{\bf
V}_{v,\tau}[N]\big)=0.
 \]
\item If $m'=1$ then, from (\ref{V[m]}), one has ${\bf
V}_{v',\tau}[1]={\bf V}_{v,\tau}[N]=1-N$. But, since $\hat
J_{v,\tau}=\hat J_{v',\tau}T(x_N-N+1)$, one has
\[
\hat J_{v,\tau}\big(x_1,\dots,x_N-1,{\bf V}_{v,\tau}[N]\big)=0,
\]
and the result is just a special case obtained from this equality by
specializing the values of the $x_i$.\hfill \qed
\end{enumerate}
 \renewcommand{\qed}{}
\end{proof}

\begin{Example}
 One has:
\[
{\bf V}_{[0,2,2,1,0,3,5,1],{1\ \ \ \atop{3\ \ \ \atop{75\ \ \atop 8642}}}}
=[5\alpha-9,3\alpha-8,3\alpha-12,4\alpha-10,5\alpha-10,2\alpha-13,-7,4\alpha-11].
\]
Indeed,
\[
\sigma_{[0,2,2,1,0,3,5,1]}= [7, 3, 4, 5, 8, 2, 1, 6]\mbox{ and
}{[0,2,2,1,0,3,5,1]}^{+}=[5,3,2,2,1,1,0,0]
\]
and the values of $v^{+}[1]-v[i]$ and
$(\CT_\tau[1]-\CT_\tau[\sigma_v[i]])-N+1$ are computed by taking the
corresponding values in the RST
\[
\begin{array}{cccc}
5\\
3\\
0&5\\
4&2&4&3
\end{array} \mbox{ and }
\begin{array}{cccc}
-9\\
-12\\
-7&-10 \\
-11& -13& -10& -8\end{array}
\]

Hence,
\[
{\hat J}_{[0,2,2,1,0,3,5,1],{1\ \ \ \atop{3\ \ \ \atop{75\ \ \atop 8642}}}}
(x_1,x_2,x_3,x_4,x_5,x_6,-7,4\alpha-11)=0.
\]
\end{Example}

The vanishing properties described in Proposition~\ref{propVanish} are
obtained by combining the actions of the $s_i$ and $\tilde\Phi$ on the
initial vectors ${\bf V}_{0^N,\tau}$.

\begin{Example}
\rm Consider the propagation of vanishing properties described in
Fig.~\ref{VG21}.
\begin{figure}[t]
\centering
\begin{tikzpicture}%
\GraphInit[vstyle=Shade]
    \tikzstyle{VertexStyle}=[shape = rectangle,
                             draw
]

\SetUpEdge[lw = 1.5pt,
color = orange,
 labelcolor = gray!30,
 labelstyle = {sloped},
 style={post}
]

\tikzset{LabelStyle/.style = {draw,
                                     text = black}}

\Vertex[x=2, y=2, L={\tiny${\bf V}_{001,{2\ \atop
31}}=[\alpha,\alpha-1,-2] \atop [x_1,x_2,-2] $}]{z1}
\Vertex[x=4, y=4, L={\tiny$[\alpha,-2,\alpha-1] \atop
[x_1,-2,\alpha-1]$}]{z2}
\Vertex[x=6, y=6, L={\tiny$[-2,\alpha,\alpha-1] \atop
[-2,\alpha-1,\alpha-1]$}]{z3}

\Edge[label={\tiny $s_2$}](z1)(z2)
\Edge[label={\tiny $s_1$}](z2)(z3)

\Vertex[x=2, y=8, L={\tiny$[\alpha-1,-2,0] \atop [x_1,-2,0]$}]{xz1}
\Vertex[x=4, y=10, L={\tiny$[-2,\alpha-1,0] \atop [-2,\alpha-1,0]$}]{xz2}
\Vertex[x=6, y=12, L={\tiny$[2\alpha,2\alpha-1,-2] \atop
[x_1,x_2,-2]$}]{xz3}
\Edge[label={$\tilde\Phi$}](z1)(xz1)
\Edge[label={$\tilde\Phi$}](z2)(xz2)
\Edge[label={$\tilde\Phi$}](z3)(xz3)

\Edge[label={\tiny $s_1$}](xz1)(xz2)

\Vertex[x=2, y=12, L={\tiny$[-2,0,\alpha-1] \atop [-2,0,\alpha-1]$}]{xz4}
\Edge[label={\tiny $s_2$}](xz2)(xz4)

\Vertex[x=4, y=14, L={\tiny$[2\alpha,-2,2\alpha-1] \atop [x_1,-2,2
\alpha-1]$}]{xz5}
\Vertex[x=2, y=16, L={\tiny$[-2,2\alpha,2\alpha-1] \atop
[-2,2\alpha,2\alpha-1]$}]{xz6}

\Edge[label={\tiny $s_2$}](xz3)(xz5)
\Edge[label={\tiny$s_1$}](xz5)(xz6)

\Vertex[x=-2, y=2, L={\tiny${\bf V}_{001,{1\ \atop
32}}=[\alpha-4,\alpha-3,-2] \atop [x_1,x_2,-2]$}]{y1}
\Vertex[x=-4, y=4, L={\tiny$[\alpha-4,-2,\alpha-3] \atop
[x_1,-2,\alpha-3]$}]{y2}
\Vertex[x=-6, y=6, L={\tiny$[-2,\alpha-4,\alpha-3] \atop
[-2,\alpha-4,\alpha-3]$}]{y3}
\Vertex[x=-2, y=8, L={\tiny$[\alpha-3,-2,-4] \atop [x_1,-2,-4]$}]{xy1}
\Vertex[x=-4, y=10, L={\tiny$[-2,\alpha-3,-4] \atop [-2,\alpha-3,-4]$}]{xy2}
\Vertex[x=-6, y=12, L={\tiny$[2\alpha-4,2\alpha-3,-2] \atop
[x_1,x_2,-2]$}]{xy3}
\Vertex[x=-2, y=12, L={\tiny$[-2,-4,\alpha-3] \atop [-2,-4,\alpha-3]$}]{xy4}
\Vertex[x=-4, y=14, L={\tiny$[2\alpha-4,-2,2\alpha-3] \atop
[x_1,-2,2\alpha-3]$}]{xy5}
\Vertex[x=-2, y=16, L={\tiny$[-2,2\alpha-4,2\alpha-3] \atop
[-2,2\alpha-4,2\alpha-3]$}]{xy6}

\Edge[label={\tiny $s_2$}](y1)(y2)
\Edge[label={\tiny $s_1$}](y2)(y3)

\Edge[label={ $\tilde\Phi$}](y1)(xy1)
\Edge[label={$\tilde\Phi$}](y2)(xy2)
\Edge[label={$\tilde\Phi$}](y3)(xy3)
\Edge[label={\tiny $s_1$}](xy1)(xy2)
\Edge[label={$\tilde\Phi$}](y3)(xy3)
\Edge[label={\tiny $s_2$}](xy2)(xy4)
\Edge[label={\tiny $s_2$}](xy3)(xy5)
\Edge[label={\tiny $s_1$}](xy5)(xy6)

\end{tikzpicture}
\caption{\label{VG21} Propagation of some vanishing properties for
$\lambda=21$.}
\end{figure}
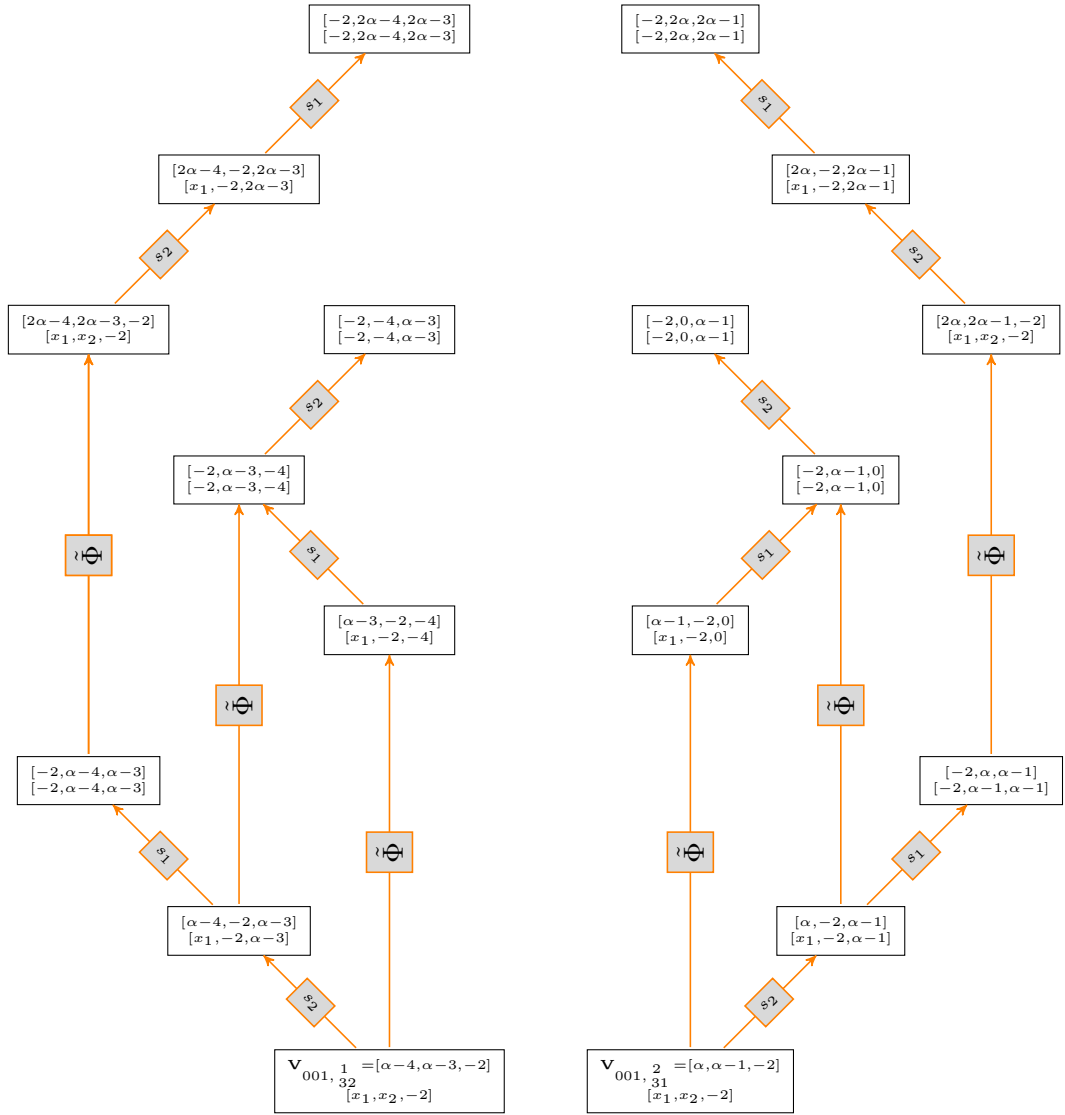
\end{Example}

Lemma~\ref{VanishPsi} and Proposition~\ref{propVanish} suggest that one
can compute other vanishing properties by combining the actions of the
$s_i$ and $\Phi^+$. A general closed formula remained to be found and,
unfortunately, the vanishing properties obtained by propagation from
${\bf V}_{0^N,\tau}$ are not suf\/f\/icient to characterize the shifted Jack
polynomials.

 \begin{Example}
  One has
\[\hat J_{[002],{2\ \atop31}}\big({\bf V}_{[001],{2\
\atop31}}s_2s_1\Phi^+\big)=\hat J_{[002],{2\
\atop31}}(\alpha,\alpha-1,\alpha-2)=0.
\]
Note that the only vanishing property obtained by propagation for $\hat
J_{[100],{2\ \atop31}}$ is
\[
\hat J_{[100],{2\ \atop31}}(-2,\alpha,\alpha-1)=0,
\]
and this is not suf\/f\/icient to characterize the polynomial $\hat
J_{[100],{2\ \atop31}}$.
\end{Example}

\section{Conclusion}

In this paper we used the Yang--Baxter graph technique to produce a structure
describing the nonsymmetric Jack polynomials whose values lie in an
irreducible $\mathfrak{S}_{N}$-module. The graph is directed with no
loops and
has exactly one root or base point. Any path joining the root to a vertex is
essentially an algorithm for constructing the Jack polynomial at that
vertex,
and the edges making up the path are the steps of the algorithm. The
edges are
labeled by the generators of the braid group or by an af\/f\/ine operation.

These techniques are used to analyze restriction to a subgroup
$\mathfrak{S}%
_{M}$ and also to construct symmetric and antisymmetric Jack polynomials.
These are associated with certain subgraphs.

Finally the graph technique is used to construct shifted, or inhomogeneous,
vector-valued Jack polynomials.

The theory is independent of the numerical value of the parameter $\alpha$
provided that the eigenspaces of $\widetilde{\xi}_{i}$ all have multiplicity
one, that is, that no two vertices of the graph have the same spectral
vector.
Future work is needed to analyze situations where this condition is
violated,
in particular when $\alpha$ has a singular value, $\left\{ \frac{n}{m}:2\leq
n\leq N,m\in\mathbb{Z},\frac{m}{n}\notin\mathbb{Z}\right\} $. There may not
be a basis of Jack polynomials for the space of all polynomials. For
particular choices of $\lambda$ there may exist symmetric Jack
polynomials of
highest weight, that is, those annihilated by
$\sum\limits_{i=1}^{N}\mathfrak{D}_{i}%
$. It seems plausible that any graph describing such a special case would be
signif\/icantly dif\/ferent from $G_{\lambda}$. As a f\/inal remark, note that
in the case of the trivial representation,
some families of highest weight Jack polynomials have been found (see
e.g.~\cite{BH,JL,BWF}) and related to the theory of the fractional quantum
Hall ef\/fect \cite{Laugh}.

\subsection*{Acknowledgments}
 The authors are grateful to A.~Lascoux
for fruitful discussions about the Yang--Baxter Graphs.
The authors acknowledge the referees for their meticulous reading and valuable comments.
This paper is partially supported by the ANR project PhysComb,
ANR-08-BLAN-0243-04.

\pdfbookmark[1]{References}{ref}
\LastPageEnding

 \end{document}